\documentclass[12pt,oneside,english]{amsart}
\usepackage[latin9]{inputenc}
\usepackage[letterpaper]{geometry}
\usepackage{babel}
\usepackage{textcomp}
\usepackage{mathrsfs}
\usepackage{amstext}
\usepackage{amsthm}
\usepackage{amssymb}
\usepackage{graphicx}
\usepackage[unicode=true,pdfusetitle,
 bookmarks=true,bookmarksnumbered=false,bookmarksopen=false,
 breaklinks=false,pdfborder={0 0 0},pdfborderstyle={},backref=false,colorlinks=false]
 {hyperref}

\makeatletter
\numberwithin{equation}{section}
\numberwithin{figure}{section}
\theoremstyle{plain}
\newtheorem*{thm*}{\protect\theoremname}
\theoremstyle{plain}
\newtheorem{thm}{\protect\theoremname}[section]
\theoremstyle{definition}
\newtheorem{defn}[thm]{\protect\definitionname}
\theoremstyle{plain}
\newtheorem{prop}[thm]{\protect\propositionname}
\theoremstyle{plain}
\newtheorem{lem}[thm]{\protect\lemmaname}
\theoremstyle{remark}
\newtheorem{rem}[thm]{\protect\remarkname}
\theoremstyle{plain}
\newtheorem{cor}[thm]{\protect\corollaryname}
\theoremstyle{remark}
\newtheorem*{rem*}{\protect\remarkname}

\usepackage{tikz}
\usepackage{tikz-cd}
\usepackage{pgfplots}

\usepackage{amsmath}
\usepackage{amssymb}

\usepackage{enumitem}
\setlist{leftmargin=0.5cm}

\usepackage[cmtip,all]{xy}
\usepackage{quiver}

\input xypic
\def\dotimes{\protect{\hspace{0.07cm}\otimes\hspace{-0.37cm}\bigcirc}}

\newcommand{\ootimes}{
  \mathbin{
    \mathchoice
      {\buildcircleotimes{\displaystyle}}
      {\buildcircleotimes{\textstyle}}
      {\buildcircleotimes{\scriptstyle}}
      {\buildcircleotimes{\scriptscriptstyle}}
  }
}

\newcommand\buildcircleotimes[1]{%
  \begin{tikzpicture}[baseline=(X.base), inner sep=0, outer sep=0]
    \node[draw,circle] (X)  {$#1\otimes$};
  \end{tikzpicture}%
}

\usepackage{mathabx,graphicx}

\makeatother

\providecommand{\corollaryname}{Corollary}
\providecommand{\definitionname}{Definition}
\providecommand{\lemmaname}{Lemma}
\providecommand{\propositionname}{Proposition}
\providecommand{\remarkname}{Remark}
\providecommand{\theoremname}{Theorem}

\begin{document}
\global\long\def\r{\mathbb{\mathbb{R}}}%

\global\long\def\n{\mathbb{\mathbb{N}}}%

\global\long\def\c{\mathbb{\mathbb{C}}}%

\global\long\def\b{\mathbb{\mathbb{B}}}%

\global\long\def\z{\mathbb{\mathbb{\mathbb{Z}}}}%

\global\long\def\q{\mathbb{\mathbb{\mathbb{Q}}}}%

\global\long\def\k{\textrm{\textbf{k}}}%

\global\long\def\h{\textrm{\textbf{H}}}%

\global\long\def\l{\ell}%

\global\long\def\o{\mathcal{O}}%

\global\long\def\nil{\emptyset}%

\global\long\def\longsquiggly{\xymatrix{{}\ar@{~>}[r]  &  {}}
 }%

\global\long\def\Hom{\mathrm{Hom}}%

\global\long\def\End{\mathrm{End}}%

\global\long\def\iso{\xrightarrow{{\sim}}}%

\global\long\def\idop{\textrm{Id}}%

\global\long\def\lra{\longrightarrow}%

\global\long\def\acts{\;\circlearrowright\;}%

\global\long\def\acted{\;\circlearrowleft\;}%

\global\long\def\dg{\mathsf{dg}}%

\global\long\def\op{\mathsf{op}}%

\global\long\def\vep{\varepsilon}%

\global\long\def\t{\tau}%

\global\long\def\xx{\tilde{x}}%

\global\long\def\tt{\tilde{\tau}}%

\global\long\def\ee{\tilde{E}}%

\global\long\def\veps{\varepsilon}%

\global\long\def\ci{\circ}%

\title{A Tensor 2-Product of 2-Representations of $\mathfrak{sl}_{2}^{+}$}
\author{Matthew McMillan}
\begin{abstract}
We construct an explicit abelian model for the operation of tensor
$2$-product for $2$-representations of $\mathfrak{sl}_{2}^{+}$,
specifically the product of a simple $2$-representation $\mathcal{L}(1)$
with a given abelian $2$-representation $\mathcal{V}$ taken from
the $2$-category of algebras. We study the case $\mathcal{V}=\mathcal{L}(1)$
in detail and identify the resulting $2$-representation. Our $2$-product
construction partially confirms a conjecture of Rouquier from 2008.

\noindent \smallskip{}

\noindent A first version of this paper appeared as arXiv:2209.06782
(2022).
\end{abstract}

\maketitle
\tableofcontents{}

\section{Introduction}

\subsection{Background and motivation}

The operation of tensor product is ubiquitous in representation theory
and its applications. It is a primary means of generating new representations
from old ones. In classical Lie theory this operation arises from
the Hopf structure of the enveloping algebra.

In \cite{craneFourDimensionalTopological1994}, Crane and Frenkel
outlined a program to build topological invariants using a higher
representation theory. The program was conceived as a way to formulate
invariants algebraically in $4d$ that upgrade known invariants in
$3d$ such as the TQFT of Witten-Reshetikhin-Turaev \cite{wittenQuantumFieldTheory1989,reshetikhinInvariants3manifoldsLink1991}.
The idea was to implement categorical versions of classical algebraic
structures. Crane and Frenkel proposed a concept of \textquoteleft Hopf
category\textquoteright{} to upgrade the Hopf structure of quantum
groups that was central to the WRT invariant.

A fully developed Hopf categorical representation theory will have
good definitions of categorical algebra, categorical representation,
and categorical Hopf structure. The notion of $2$-representation
was provided with a good definition for $\mathfrak{sl}_{2}$ in work
of Chuang-Rouquier \cite{chuangDerivedEquivalencesSymmetric2008a},
and the graded case descending to quantized structures in work of
Lauda \cite{laudaCategorificationQuantumMathfraksl2010}. The definitions
were generalized to all Kac-Moody algebras in \cite{rouquier2KacMoodyAlgebras2008}
and \cite{khovanovDiagrammaticApproachCategorification2009,khovanovDiagrammaticApproachCategorification2011}.

A tensor $2$-product operation would give the higher analog of an
aspect of Hopf structure, or at least of the expression of Hopf structure
on the category of representations of the algebra. A $2$-product
is defined in an $\mathcal{A}_{\infty}$ setting by Rouquier \cite{rouquierHopfCategories},
but no explicit formulas are known for the product action in that
setting, and the setting itself brings significant technical complications.
Rouquier has conjectured \cite{rouquierHigherRepresentationTheory2008}
that a subcategory affording an abelian $2$-representation should
exist. The main construction of this paper partially verifies his
conjecture by identifying an abelian $2$-product when one factor
is $\mathcal{L}(1)$ and the other factor $\mathcal{V}$ is taken
from the $2$-category of algebras. In addition, our construction
takes a step toward defining a practicable $2$-product by providing
explicit formulas for the $2$-representation component structures.

In early work of Bernstein-Frenkel-Khovanov \cite{bernsteinCategorificationTemperleyLiebAlgebra1999},
the authors consider a category whose Grothendieck group is the tensor
product of fundamental representations. Their methods were extended
by Stroppel \cite{stroppelCategorificationTemperleyLiebCategory2005}
and others (cf.~\cite{frenkelCategorificationFinitedimensionalIrreducible2007,mazorchukCombinatorialApproachFunctorial2009,sartoriCategorificationTensorProduct2015,sussanCategoryMathfrakslLink2007})
to find a category with Grothendieck group isomorphic to any given
tensor product of finite dimensional simples in type $A$. Graphical
methods were developed by Webster \cite{websterKnotInvariantsHigher2017,websterTensorProductAlgebras2016}
to produce categories for tensor products of simples for general Kac-Moody
algebras. We expect these categories to be equivalent to tensor $2$-products.
Ebert and Rouquier have established one such equivalence \cite{ebertHigherTensorProduct2025}.

The Crane-Frenkel program for building TQFTs gives perhaps the most
compelling motivation to find a categorical product. In the case of
$\mathfrak{sl}_{2}$, a $2$-product will play a central role in a
prospective $4d$ TQFT that extends Khovanov homology. Glimmers of
this $4d$ theory have been seen by physicists \cite{gukovFivebranes3manifoldHomology2017},
and some aspects are defined rigorously in some cases \cite{gukovTwovariableSeriesKnot2021}.
Along these lines, recent work of Manion-Rouquier \cite{manionHigherRepresentationsCornered2020}
on the case of the super Lie algebra $\mathfrak{gl}(1|1)^{+}$ shows
that a $2$-product can be used to describe Bordered Heegaard-Floer
theory for surfaces \cite{lipshitzBorderedHeegaardFloer2018a}. Lastly,
it is now known that $4d$ invariants rooted in quantum group representations
are able to detect smooth structure \cite{renKhovanovHomologyExotic2025}.

\subsection{\label{subsec:Results}Result}

Let $\mathcal{U}^{+}$ denote the monoidal category associated to
the positive half of the enveloping algebra of $\mathfrak{sl}_{2}$.
Let us be given a field $k$ and the data of a $k$-algebra $A$ and
a triple $(E,x,\tau)$ as follows. Let $E$ be an $(A,A)$-bimodule,
let $x\in\End(E)$ and $\tau\in\End(E^{2})$ be bimodule endomorphisms,
and suppose that $x$ and $\tau$ generate an action of the nil affine
Hecke algebra, that is, that they satisfy the following relations:
\begin{gather*}
\tau^{2}=0,\\
\tau E\circ E\tau\circ\tau E=E\tau\circ\tau E\circ E\tau,\\
\tau\circ Ex=xE\circ\tau+1,\;Ex\circ\tau=\tau\circ xE+1.
\end{gather*}
(Here we write $xE$ for the endomorphism $x\otimes\idop_{E}$ in
$\End(E^{2})$, and similarly for the others.) This data determines
a $2$-representation $\mathcal{V}$ of $\mathcal{U}^{+}$.

We can give such data for a simple $2$-representation $\mathcal{L}(1)$
of $\mathcal{U}^{+}$ that categorifies the fundamental representation
$L(1)$ of $\mathfrak{sl}_{2}$. The $k$-algebra is $k[y]_{+1}\times k[y]_{-1}$
(decomposed into weight algebras), and the triple is $(k[y],y,0)$.
Here $y\in k[y]_{-1}$ acts on $k[y]$ on the right by multiplication,
and $y\in k[y]_{+1}$ acts by zero. These roles are reversed for the
left action. The endomorphism $x$ acts by multiplication by $y$.

Let $P_{n}=k[x_{1},\dots,x_{n}]$ be the polynomial algebra. Then
$P_{n}$ acts on $E^{n}$ with $x_{i}\in P_{n}$ acting by the endomorphism
$E^{n-i}xE^{i-1}$.

This paper is organized around the following construction of $\mathcal{L}(1)\ootimes\mathcal{V}$:
\begin{thm*}[Main result]
 Suppose $x$ and $\tau$ satisfy the nil affine Hecke relations,
so $(E,x,\tau)$ gives a $2$-representation of $\mathcal{U}^{+}$
for the algebra $A$, denoted $\mathcal{V}$, and suppose the bimodule
$E$ has the following additional properties: 
\begin{itemize}
\item $_{A}E$ is finitely generated and projective,
\item $E^{n}$ is free as a $P_{n}$-module.
\end{itemize}
Then we define explicitly: 
\begin{itemize}
\item a $k$-algebra $C$ (Def.~\ref{def:C}),
\item a bimodule $\tilde{E}$ (Def.~\ref{def:tilde-E}),
\item endomorphisms $\tilde{x}$ and $\tilde{\tau}$ (Def.~\ref{def:tilde-x-tilde-tau}),
\end{itemize}
such that $\tilde{x}$ and $\tilde{\tau}$ satisfy the nil affine
Hecke relations, so $(\tilde{E},\tilde{x},\tilde{\tau})$ gives the
data of a $2$-representation of $\mathcal{U}^{+}$ for $C$ that
we call $\mathcal{L}(1)\ootimes\mathcal{V}$.
\end{thm*}
At the end of this paper we describe the case $\mathcal{L}(1)\ootimes\mathcal{L}(1)$
in detail, showing that it is the expected categorification of $L(1)\otimes L(1)$.

In a separate paper \cite{mcmillanTensor2product22024} we consider
the extension of the construction given in this paper to actions of
the full $2$-category $\mathcal{U}$ associated to the enveloping
algebra of $\mathfrak{sl}_{2}$, and not only its positive half. When
the functor $E\otimes_{A}-$ has a right adjoint given by tensor product
with a bimodule $F$, and the pair of them satisfies some additional
relations that categorify the commutator identities, the action is
said to give a $2$-representation of $\mathcal{U}$. We show that
if the original bimodule $E$ has such an adjoint $F$ giving with
$E$ an action of the full $\mathcal{U}$ on $\mathcal{V}$, then
there is also a bimodule $\tilde{F}$, given as the right-dual of
$\tilde{E}$, which together with $\tilde{E}$ provides an action
of the full $\mathcal{U}$ on $\mathcal{L}(1)\ootimes\mathcal{V}$.

In unpublished work, Rouquier defines an $\infty$-category for tensor
products including $\mathcal{L}(1)\ootimes\mathcal{V}$ that are provided
with an $\infty$-categorical action of $\mathcal{U}$ on the derived
category $D^{b}(\mathcal{L}(1)\ootimes\mathcal{V})$. Our algebra
$C$ is the endomorphism ring of a compact generator $X$ in Rouquier\textquoteright s
$D^{b}(\mathcal{L}(1)\ootimes\mathcal{V})$. This $X$ has vanishing
shifted self-homs, so it is tilting and defines a $t$-structure as
a projective generator of the heart.

In subsequent work, Ebert and Rouquier have considered the special
case where $\mathcal{V}$ is the regular $2$-representation of $\mathcal{U}$
\cite{ebertHigherTensorProduct2025}. They produce an infinite family
of objects that generate $D^{b}(\mathcal{L}(1)\ootimes\mathcal{U})$
as a thick subcategory and have vanishing shifted homs, thereby giving
a tilting subcategory $\mathcal{T}$ determining the same $t$-structure.
Moreover, this $\mathcal{T}$ has a $(\mathcal{U},\mathcal{U})$-bimodule
structure with left $\mathcal{U}$-action inherited from the $\infty$-categorical
$\mathcal{U}$-action on $D^{b}(\mathcal{L}(1)\ootimes\mathcal{U})$.
They show that $\mathcal{T}$ is a tensor product category defined
by Webster \cite[§4.2]{websterKnotInvariantsHigher2017}.

\subsection{\label{subsec:Technique}Technique }

Let us be given $\mathcal{V}$ as described above. Write $E_{y}$
for the $(A[y],A[y])$-bimodule $E[y]\Big/(x-y)E[y]$. We begin with
a \textquoteleft naive\textquoteright{} algebra $B$ formed from the
underlying data of $\mathcal{L}(1)$ and $\mathcal{V}$: 
\[
B=\begin{pmatrix}A[y] & E_{y}\\
0 & A[y]
\end{pmatrix}.
\]
 There is a natural candidate $E'$ for the diagonal action of $\mathcal{U}^{+}$,
but it is a complex of $(B,B)$-bimodules, not a bimodule. It is given
as a complex in degrees $0$ and $1$ by 
\[
E'=\begin{pmatrix}E[y] & E[y]E_{y}\\
0 & E[y]
\end{pmatrix}\overset{d}{\to}\begin{pmatrix}E_{y} & E_{y}E_{y}\\
A[y] & E_{y}
\end{pmatrix}.
\]
(The differential and action data are described in Definition \ref{def:E'}.)
There is also a natural candidate for $x\in\End(E')$ arising from
the data of $\mathcal{L}(1)$ and $\mathcal{V}$, but that $x$ is
not equivariant over the action of generators in $E_{y}$ in $B$. 

Let $e_{1}=$\scalebox{0.7}{$\begin{pmatrix}1 & 0\\
0 & 0
\end{pmatrix}$}$\,\in B$. Our technique in this paper is to define a new algebra
\[
C=\End_{K^{b}(B)}(Be_{1}\oplus E'e_{1})
\]
 that is derived-equivalent to $B$. The bimodule complex $E'$ may
be transported through the equivalence, and the result is quasi-isomorphic
to a complex $\tilde{E}$ of $(C,C)$-bimodules that is concentrated
in degree $0$ and projective on the left. We consider $\tilde{E}$
to be a $(C,C)$-bimodule, and we construct explicit bimodule endomorphisms
$\tilde{x}\in\End(\tilde{E})$ (compatible with $x$) and $\tilde{\tau}\in\End(\tilde{E}^{2})$
that satisfy the nil affine Hecke relations. The data $(C,\tilde{E},\tilde{x},\tilde{\tau})$
determines a $2$-representation that we call $\mathcal{L}(1)\ootimes\mathcal{V}$.

In order to define $\tilde{x}$ and $\tilde{\tau}$ and verify the
relations, we study the tensor powers $\tilde{E}^{n}$. These powers
can be parametrized by explicit models containing $\Hom_{K^{b}(B)}(E'e_{1},E'^{n}e_{1})$.
We give presentations of these modules by generators and relations
for $n=1,2,3,4$.

\subsection{Explanation }

Suppose $_{H}M$ and $_{H}N$ are two representations of a Hopf $k$-algebra
$H$ with coproduct $\Delta:H\to H^{2}$ and antipode $S:H\to H$.
There is a large outer product $M\otimes_{k}N$ with two commuting
actions of $H$ on the two factors, and a third, diagonal, action
given by first applying $\Delta$. There is a smaller product $M\otimes_{H}N$
using $S$ to view $M$ as a right $H$-module. The smaller product
is related to the larger one as follows: $M\otimes_{H}N$ is the largest
quotient of $M\otimes_{k}N$ on which $\Delta(H)$ acts by $0$. This
can be seen using the formulas $\Delta(h)=h\otimes1+1\otimes h$ and
$S(h)=-h$ for enveloping algebras of Lie algebras, with which the
condition $\Delta(h).(m\otimes n)=0$ may be written $m.h\otimes n=m\otimes h.n$. 

Now let $\mathcal{V}_{i}$ be an abelian category of $A_{i}$-modules
for $i=1,2$, where $\mathcal{V}_{i}$ is a $2$-representation of
$\mathcal{U}^{+}$ with data $(E_{i},x^{i},\tau^{i})$. We can easily
define a large outer product category $\mathcal{V}_{1}\boxtimes_{k}\mathcal{V}_{2}$
that has two commuting actions of $\mathcal{U}^{+}$. We seek a kind
of diagonal action of $\mathcal{U}^{+}$ on $\mathcal{V}_{1}\boxtimes_{k}\mathcal{V}_{2}$.
One can also describe a smaller product without diagonal $\mathcal{U}^{+}$-symmetry.
Objects should be generated by pairs of modules $M\in\mathcal{V}_{1}$,
$N\in\mathcal{V}_{2}$ together with functorial isomorphisms $E_{1}(M)\otimes_{k}N\iso M\otimes_{k}E_{2}(N)$
that are equivariant over the actions of $x^{i}$ on $E_{i}$ and
$\tau^{i}$ on $E_{i}^{2}$. These isomorphisms categorify the conditions
$\Delta(e).(m\otimes n)=0$.

At this point we make three conceptual moves. First, we expand the
larger product category by including with each pair $M\in\mathcal{V}_{1}$,
$N\in\mathcal{V}_{2}$ a morphism $\alpha_{M}^{N}:E_{1}(M)\otimes_{k}N\to M\otimes_{k}E_{2}(N)$,
functorial in $M$ and $N$, that is $x$- and $\tau$-equivariant.
So we define objects of $\mathcal{V}_{1}\ootimes\mathcal{V}_{2}$
to be triples of the form $(M,N;\alpha_{M}^{N})$. Second, we consider
morphisms $\alpha_{M}^{N}$ as two-term chain complexes, in particular
mapping cones, and move to a derived context. Third, for the new diagonal
action of $E$ on $(M,N;\alpha_{M}^{N})$ we take the cone complex
$C=Cone(\alpha_{M}^{N})$ itself. In the derived category, this complex
is zero precisely when $\alpha_{M}^{N}$ is an isomorphism, which
is the correspondence we sought.

To complete the idea, it is necessary to supply natural $x$- and
$\tau$-equivariant morphisms $\alpha_{(E_{1}\otimes\idop)C}^{(\idop\otimes E_{2})C}$
in order to make $C$ an object in $\mathcal{V}_{1}\ootimes\mathcal{V}_{2}$,
and to supply endomorphisms $x$ and $\tau$ of $Cone(\alpha_{M}^{N})$
and $Cone(Cone(\alpha_{M}^{N}))$ satisfying Hecke-type relations
in order to make a $2$-representation of $\mathcal{U}$ using $Cone(\alpha_{M}^{N})$
for the image of $E$. Here one encounters further technical difficulties.
In \cite{rouquierHopfCategories}, Rouquier is expected to give a
general definition of tensor $2$-product by working in an $\mathcal{A}_{\infty}$
setting that encodes the technical difficulties as higher homotopies.
For example, the failure of equivariance of the natural $x\in\End(E')$
mentioned in §\ref{subsec:Technique} can be expressed as a homotopy.

In our setting for $\mathcal{L}(1)\ootimes\mathcal{V}$, we have $\mathcal{L}(1)$
given by the data $(A^{\ci},k[y],y,0)$ with $A^{\ci}=k[y]_{+1}\times k[y]_{-1}$,
and $\mathcal{V}$ given by the data $(A,E,x,\tau)$. One can define
a tensor algebra $B'$: 

\[
B'=T_{A^{\ci}\otimes_{k}A}\bigl(\phantom{}^{\vee}k[y]\otimes_{k}E\bigr).
\]
 There is a canonical isomorphism $\phantom{}^{\vee}k[y]\otimes_{k}E\iso E[y]$,
and another $A^{\ci}\otimes_{k}A\iso A[y]\times A[y]$. The data of
a $B'$-module is equivalent to the data of a triple $(M,N,\alpha_{M}^{N})$
where $M,N\in A[y]\text{-mod}$ and $\alpha_{M}^{N}:E[y]\otimes_{A[y]}M\to N$.
Since $\tau^{1}=0$ in this case, $\alpha$ is automatically $\tau$-equivariant.
We can enforce $x$-equivariance of $\alpha$ by taking a quotient
by $I=\mathrm{Im}(x-y)$, where $x-y$ is understood in $\End_{A[y]}(E[y])$.
Then the algebra $B'/I$ is isomorphic to the algebra $B$ in §\ref{subsec:Technique}.

\subsection{Outline summary}

The paper is organized as follows:
\begin{itemize}
\item In §2 we describe some conventions and background theory. We are working
in the setting of monoidal categories of the form $\mathsf{Bim}_{k}(A)$
for a $k$-algebra $A$: objects are $(A,A)$-bimodules, morphisms
are bimodule maps. The data of a $2$-representation of $\mathcal{U}^{+}$
consists of an algebra $A$, a bimodule $_{A}E_{A}$, and endomorphisms
$x\in\End(E)$ and $\tau\in\End(E^{2})$ satisfying nil affine Hecke
relations.
\item In §3 we begin with a naive product algebra $B$ and complex of bimodules
$_{B}E'_{B}$. We construct a derived-equivalent algebra $C$. We
define a $(C,C)$-bimodule $\tilde{E}$ and study a new class of bimodules
we call $G_{n}$ that arise inside the tensor powers of $\tilde{E}$.
This study has a technical and computational flavor.
\item In §4 we construct the new nil affine Hecke action, with generators
$\tilde{x}$ and $\tilde{\tau}$, on powers of the new bimodule $\tilde{E}$.
More computations are required to establish the properties we need.
They rely on results about $G_{n}$ proved in §3.
\item In §5 we give explicit details for the most basic example of our construction:
$\mathcal{L}(1)\ootimes\mathcal{L}(1)$. This product agrees with
a well-known categorification of $L(1)\otimes L(1)$, where $L(1)$
is the fundamental representation of $\mathfrak{sl}_{2}$.
\end{itemize}

\subsection{Acknowledgments}

I thank Raphaël Rouquier for sharing his ideas, and for advice along
the way with many aspects of this project. I am indebted for his time
and attention. I thank Laurent Vera for working with me in early stages
of the project. I thank You Qi and Weiqiang Wang for helpful comments
on the manuscript. I thank David Reutter for opening my interest in
categorical mathematics while at Cambridge some years ago.

This work was supported by an NDSEG Research Fellowship from the US
Department of Defense and by the NSF through grant DMS-1702305.

\section{Background structures}

Let $k$ be a field.

\subsection{Nil affine Hecke algebras}

The nil affine Hecke algebra $^{0}H_{n}$ is the $k$-algebra with
generators $x_{1},\dots,x_{n},\tau_{1},\dots,\tau_{n-1}$ and relations:
\begin{gather*}
x_{i}x_{j}=x_{j}x_{i},\tau_{i}^{2}=0,\\
\tau_{i}\tau_{i+1}\tau_{i}=\tau_{i+1}\tau_{i}\tau_{i+1},\\
\tau_{i}\tau_{j}=\tau_{j}\tau_{i}\text{ if }|i-j|>1,\\
\tau_{i}x_{j}=x_{j}\tau_{i}\text{ if }j-i\notin\{0,1\},\\
\tau_{i}x_{i}=x_{i+1}\tau_{i}+1,x_{i}\tau_{i}=\tau_{i}x_{i+1}+1.
\end{gather*}

Define $s_{i}=\tau_{i}(x_{i}-x_{i+1})-1$. Observe that $s_{i}^{2}=1$
and $s_{i}\circ\tau_{i}=\tau_{i}$.

\subsection{$\mathcal{U}^{+}(\mathfrak{sl}_{2})$ and its 2-representations}

\subsubsection{Monoidal category $\mathcal{U}^{+}$}
\begin{defn}
\label{def:u(sl_2)} Let $\mathcal{U}^{+}(\mathfrak{sl}_{2})$ (hereafter
\textquoteleft $\mathcal{U}^{+}$\textquoteright ) be the strict monoidal
$k$-linear category generated by an object $E$ and maps $x:E\to E$
and $\tau:E^{2}\to E^{2}$ subject to the relations: 
\begin{gather}
\tau^{2}=0,\label{eq:U-relns}\\
\tau E\circ E\tau\circ\tau E=E\tau\circ\tau E\circ E\tau,\\
\tau\circ Ex=xE\circ\tau+1,\;Ex\circ\tau=\tau\circ xE+1.\label{eq:U-relns-last}
\end{gather}
\end{defn}

We write $s=\tau\circ(Ex-xE)-1$. Observe that $s^{2}=1$ and $s\circ\tau=\tau$.

One easily checks that non-trivial $\Hom$ spaces of $\mathcal{U}^{+}$
are Hecke algebras:
\begin{prop}
The objects of $\mathcal{U}^{+}$ are the $E^{n}$ for $n\in\z{}_{\geq0}$,
and 
\[
\Hom(E^{n},E^{m})\cong\begin{cases}
^{0}H_{n} & n=m\\
0 & n\neq m
\end{cases}
\]
with the isomorphism from $^{0}H_{n}$ given by $x_{i}\mapsto E^{n-i}xE^{i-1}$,
$\tau_{i}\mapsto E^{n-i-1}\tau E^{i-1}$. Using the obvious morphism
$^{0}H_{n}\otimes\phantom{}^{0}H_{m}\to\phantom{}^{0}H_{n+m}$, the
diagram commutes: 

\noindent\begin{minipage}[t]{1\columnwidth}%
\[\begin{tikzcd}[ampersand replacement=\&] 	{{^0 H_n}\otimes {^0H_m}} \&\& {^0 H_{n+m}} \\ 	\\ 	{\mathrm{End}(E^n)\otimes \mathrm{End}(E^m)} \&\& {\mathrm{End}(E^{n+m}).} 	\arrow[from=1-1, to=1-3] 	\arrow["\cong"', from=1-1, to=3-1] 	\arrow["\cong"', from=1-3, to=3-3] 	\arrow["\otimes", from=3-1, to=3-3] \end{tikzcd}\]%
\end{minipage}
\end{prop}

\subsubsection{$2$-representations of $\mathcal{U}^{+}$}
\begin{defn}
A $2$-representation of $\mathcal{U}^{+}$ on a category $\mathcal{V}$
is a strict monoidal functor $\mathcal{U}^{+}\to\End(\mathcal{V})$.
The data of such a functor consists of an endofunctor $E$ of $\mathcal{V}$
and natural transformations $x\in\End(E)$, $\tau\in\End(E^{2})$
satisfying (\ref{eq:U-relns})--(\ref{eq:U-relns-last}). A morphism
of $2$-representations $(\mathcal{V},E,x,\tau)\to(\mathcal{V}',E',x',\tau')$
consists of a functor $\Phi:\mathcal{V}\to\mathcal{V}'$ and an isomorphism
of functors $\varphi:\Phi E\iso E'\Phi$ such that: 
\begin{align*}
\varphi\circ\Phi x=x'\Phi\circ\varphi & :\Phi E\to E'\Phi,\\
E'\varphi\circ\varphi E\circ\Phi\tau=\tau'\Phi\circ E'\varphi\circ\varphi E & :\Phi E^{2}\to E'^{2}\Phi.
\end{align*}
\end{defn}

Note that $\End(\mathcal{V})$ is the full sub-$2$-category of the
$2$-category of categories $\mathsf{Cat}$ generated by the object
$\mathcal{V}$. One can define $\mathcal{U}^{+}$ as a $2$-category
with a single object, so that the data of $2$-representation is the
data of $2$-functor from $\mathcal{U}^{+}$ to $\mathsf{Cat}$. This
justifies our \textquoteleft $2$\textquoteright{} prefixes.

In this paper we study monoidal functors from $\mathcal{U}^{+}$ to
monoidal categories of the form $\mathsf{Bim}_{k}(A)$ which are defined
for $k$-algebras $A$ as follows: the objects of $\mathsf{Bim}_{k}(A)$
are $(A,A$)-bimodules, and the morphisms of $\mathsf{Bim}_{k}(A)$
are bimodule maps. The monoidal structure on $\mathsf{Bim}_{k}(A)$
is given by tensor product of bimodules over $A$.

Note that there is a $2$-category $\mathsf{Alg}_{k}$ with $k$-algebras,
bimodules, and bimodule maps as the objects, $1$-morphisms, and $2$-morphisms.
Then $\mathsf{Bim}_{k}(A)$ is the full sub-$2$-category of $\mathsf{Alg}_{k}$
generated by the object $A$.
\begin{prop}
\label{prop:Def-2-rep-U} The data of a $2$-representation $\mathcal{U}^{+}\to\mathsf{Bim}_{k}(A)$
for a $k$-algebra $A$ consists of a bimodule $_{A}E_{A}$ and bimodule
maps $x\in\End(E)$, $\tau\in\End(E^{2})$ that satisfy (strictly)
the relations of $\mathcal{U}^{+}$.
\end{prop}

We will use \textquoteleft $x_{i}$\textquoteright{} and \textquoteleft $\tau_{i}$\textquoteright{}
to denote the generators in any $^{0}H_{n}$ (where $i\leq n$ for
$x_{i}$ and $i<n$ for $\tau_{i}$ are assumed). Given a $2$-representation
for a $k$-algebra $A$ with bimodule $E$, these symbols are also
used to denote the corresponding elements in each $\End(E^{n})$.

\subsubsection{\label{subsec:Minimal-rep} The $2$-representation $\mathcal{L}(1)$}

A simple $2$-representation of $\mathcal{U}^{+}$ is given for the
algebra $A=A_{+1}\times A_{-1}$, $A_{i}=k[y]$, by the bimodule $E=k[y]$,
where $y\in A_{-1}$ acts on the left by $0$ and on the right by
multiplication by $y$, and $y\in A_{+1}$ acts on the right by $0$
and the left by $y$. The Hecke actions are generated by $x\in\End(E)$
acting by multiplication by $y$, and $\tau\in\End(E^{2})$ satisfies
$\tau=0$ because $E^{2}=0$.

\subsection{Further conventions}

Assume we are given data $(A,E,x,\tau)$ determining a $2$-representation,
and fix these through §4. Assume that $_{A}E$ is f.g.~projective
and that $E^{n}$ is free as a $P_{n}$-module.

Consider the endomorphism $x-y$ of the $(A[y],A[y])$-bimodule $E[y]$.
Its image $(x-y)E[y]$ is a sub-bimodule of $E[y]$. Write $E_{y}$
for the quotient $E[y]\big/(x-y)E[y]$. (Alternatively: $E_{y}$ is
$E$ extended to an $(A[y],A[y])$-bimodule by specifying that $y$
acts on both sides by $x$.) The projection 
\begin{align*}
\pi:E[y] & \to E_{y}\\
ey^{n} & \mapsto x^{n}(e)
\end{align*}
 is a surjection of bimodules.

We simplify notation for tensor products by adopting a convention
that concatenation indicates the tensor product over an algebra that
is clear from the context. Sometimes it will be unclear whether a
tensor product is meant over $A$ or over $A[y]$, so we further stipulate
that if the expression for a module contains \textquoteleft $y$\textquoteright ,
it will be understood as an $A[y]$-module, and if the expression
lacks \textquoteleft $y$\textquoteright , it will be understood as
an $A$-module. Concatenation will indicate tensor product over $A[y]$
if both are $A[y]$-modules, otherwise it will indicate tensor product
over $A$.

We will tacitly use canonical isomorphisms such as 
\[
M[y]\otimes_{A[y]}N[y]\iso M[y]\otimes_{A}N\iso(MN)[y]
\]
 for $M$ a right $A$-module and $N$ a left $A$-module. For example,
$EE_{y}$ denotes $E\otimes_{A}E_{y}$ according to our convention,
but this is canonically isomorphic to $E[y]\otimes_{A[y]}E_{y}$,
and the latter may be written $E[y]E_{y}$. So we may write either
$EE_{y}$ or $E[y]E_{y}$ with equivalent meanings.

Extend $x$ to $\End(E[y])$ by $x:ey^{n}\mapsto x(e)y^{n}$ and $\tau$
to $\End(E^{2}[y])$ by $\tau:eey^{n}\mapsto\tau(ee)y^{n}$. The map
$s$ defined above in terms of $x$ and $\tau$ extends in the same
manner to a map in $\End(E^{2}[y])$. Note that we denote an arbitrary
element of $E[y]$ by the single letter \textquoteleft $e$\textquoteright .
Similarly an arbitrary element of $E^{2}[y]$ is denoted by the doubled
symbol \textquoteleft $ee$\textquoteright , which may well not be
a simple tensor of the form $e\otimes e$. Later we will use \textquoteleft $eee$\textquoteright{}
or \textquoteleft $eee_{i}$\textquoteright{} as suggestive notation
for elements of $E^{3}[y]$, and so on.

Define $\delta=\tau\circ(Ex-y)\in\End(E^{2}[y])$. We also consider
the extensions of $x_{i}$ and $\tau_{i}$ to $E^{n}[y]$, and then
$s_{i}$ and $\delta_{i}$ defined by their same formulas but replacing
$x$ with $x_{i}$ and $\tau$ with $\tau_{i}$. Some important identities
are quickly verified: 
\begin{lem}
We have: 
\begin{itemize}
\item $s^{2}=1$, so $s$ is an isomorphism
\item $\delta^{2}=\delta$, so $\delta$ is an idempotent,
\end{itemize}
and we also have $s_{i}^{2}=1$ and $\delta_{i}^{2}=\delta_{i}$.
\end{lem}

We adopt a flexible notation $y_{i}=x_{i}-y$ until §\ref{sec:Comparison}.
Here $y_{i}$ indicates $\bigl(E^{j}xE^{i-1}-y\bigr)$ for some $j$,
and context will determine the value of $j$. Note that $\delta_{i}=\tau_{i}y_{i}$.

One may check that $s\circ x_{2}=x_{1}\circ s$ and $s\ci x_{1}=x_{2}\ci s$.
It follows that $s$ exchanges $y_{2}$ and $y_{1}$ and descends
to a map: 
\[
s:E_{y}\otimes_{A[y]}E[y]\to E[y]\otimes_{A[y]}E_{y}.
\]
 So we have $s:E^{2}\to E^{2}$ a map of $(A,A)$-bimodules, and this
induces $s:E^{2}[y]\to E^{2}[y]$ as well as $s:E_{y}E\to EE_{y}$,
maps of $(A[y],A[y])$-bimodules. Context will determine the domain
and codomain for the symbol $s$.
\begin{lem}
We also have: 
\begin{itemize}
\item $\pi_{1}\circ\delta=s\circ\pi_{2}:E^{2}[y]\to EE_{y}$.
\end{itemize}
\end{lem}

We define projections $\pi_{i}:E^{n}[y]\to E^{n-i}E_{y}E^{i-1}=E^{n}[y]\big/(y_{i}$)
by $\pi_{i}=E^{n-i}\pi E^{i-1}$. The same names may be used for maps
between products with $E_{y}$ factors, for example $\pi_{2}:EE_{y}\to E_{y}E_{y}$.

Given a module $_{A}M$, its algebra of endomorphisms $\End_{A}(_{A}M)$
will use the traditional order of composition for multiplication:
$(f\circ g)(m)=f(g(m))$. Typically, but not always, \textquoteleft $\circ$\textquoteright{}
is written to emphasize this convention. A consequence is that for
a ring $A$, the algebra $\End_{A}(_{A}A)$ is identified with $A^{\op}$.

Given two complexes $M$, $N$ of $A$-modules, we will write $\mathscr{H}om_{A}(M,N)$
for the complex generated by homogeneous $A$-module homomorphisms
from $M$ to $N$. In degree $n$ it is given by homogeneous maps
of degree $n$, and the differential is $d(f)=d\circ f-(-1)^{|f|}f\circ d$
for $f$ a homogeneous map of degree $|f|$. The notation $Z^{i}M$
refers to the degree $i$ part of the kernel of $d$.

Given an algebra $R$, we write $D^{b}(R)$ for the derived category
of bounded complexes of left $R$-modules. A strictly perfect complex
of left $R$-modules is a bounded complex of finitely generated projective
$R$-modules. The category $\text{per }R\subset D^{b}(R)$ is the
full subcategory of complexes quasi-isomorphic to strictly perfect
complexes. Given $M\in D^{b}(R)$, we write $\langle M\rangle_{\Delta}$
for the smallest triangulated strictly full subcategory of $D^{b}(R)$
closed under direct summands and containing $M$.
\begin{lem}
We have $\langle R\rangle_{\Delta}=\text{per }R$.
\end{lem}

\subsection{Generalized matrix algebras and tensor product\label{subsec:matrix-tensor-products}}

Suppose we are given $k$-algebras $A$ and $D$, bimodules $_{A}B_{D}$
and $_{D}C_{A}$, and bimodule maps 
\begin{align*}
_{A}B\otimes_{D}C_{A} & \xrightarrow{\gamma_{1}}A\\
_{D}C\otimes_{A}B_{D} & \xrightarrow{\gamma_{2}}D.
\end{align*}
 With this data we can define a new $k$-algebra $R$: 
\[
R=\begin{pmatrix}A & B\\
C & D
\end{pmatrix},
\]
where multiplication of matrices is defined with the customary formulas
using the above bimodule structures and maps.

A right $R$-module consists of the data of $M_{1}$ a right $A$-module,
$M_{2}$ a right $D$-module, a map $M_{1}\otimes_{A}B\xrightarrow{\alpha}M_{2}$
of right $D$-modules, and a map $M_{2}\otimes_{D}C\xrightarrow{\beta}M_{1}$
of right $A$-modules, such that the latter two maps are compatible
with $\gamma_{1}$ and $\gamma_{2}$. Here compatibility with $\gamma_{1}$,
for example, means that the following compositions agree: 
\begin{align*}
M_{1}\otimes_{A}(B\otimes_{D}C) & \xrightarrow{\idop_{M_{1}}\otimes\gamma_{1}}M_{1}\otimes_{A}A\iso M_{1}\\
(M_{1}\otimes_{A}B)\otimes_{D}C & \xrightarrow{\alpha\otimes\idop_{C}}M_{2}\otimes_{D}C\xrightarrow{\beta}M_{1}.
\end{align*}
 The data of a left $R$-module may be given in a similar form.

Let 
\[
M=\begin{pmatrix}M_{1} & M_{2}\end{pmatrix}
\]
be a right $R$-module, and 
\[
N=\begin{pmatrix}N_{1}\\
N_{2}
\end{pmatrix}
\]
a left $R$-module. Their tensor product $M\otimes_{R}N$ may be formed
as follows. Consider the pair of maps given by the $R$ action data:
\begin{align*}
M_{1}\otimes_{A}B\otimes_{D}N_{2} & \overset{I_{B}}{\longrightarrow}M_{1}\otimes_{A}N_{1}\oplus M_{2}\otimes_{D}N_{2}\\
M_{2}\otimes_{D}C\otimes_{A}N_{1} & \overset{I_{C}}{\longrightarrow}M_{1}\otimes_{A}N_{1}\oplus M_{2}\otimes_{D}N_{2}
\end{align*}
by $I_{B}(m\otimes b\otimes n)=m\otimes b.n-m.b\otimes n$ and likewise
for $I_{C}$. Then we have an isomorphism: 
\[
\left(M_{1}\otimes_{A}N_{1}\oplus M_{2}\otimes_{D}N_{2}\right)\big/(I_{B}+I_{C})\iso M\otimes_{R}N.
\]

Now let $F\in\End_{R}(N)$ be an endomorphism of left $R$-modules.
It determines an endomorphism $\idop_{M}\otimes_{R}F\in\End_{k}(M\otimes_{R}N)$
which will be denoted $MF$. We can study this on components as follows.
There are induced endomorphisms $F_{1}\in\End_{A}(N_{1})$ and $F_{2}\in\End_{D}(N_{2})$
given by restriction of $F$. These determine endomorphisms $M_{1}F_{1}\in\End_{k}(M_{1}\otimes_{A}N_{1})$
and $M_{2}F_{2}\in\End_{k}(M_{2}\otimes_{D}N_{2})$, and these in
turn provide together an endomorphism $\begin{pmatrix}M_{1}F_{1} & 0\\
0 & M_{2}F_{2}
\end{pmatrix}$ of $M_{1}\otimes_{A}N_{1}\oplus M_{2}\otimes_{D}N_{2}$. The property
of full $R$-linearity of $F$ implies that this morphism preserves
the submodules $I_{B}$ and $I_{C}$, and descends to the quotient
$M\otimes_{R}N$ where it agrees with $MF$.
\begin{lem}
\label{lem:tensor-components-induced} In the notations used above,
an element of $\End_{k}(M\otimes_{R}N)$ of the form $MF$ for $F\in\End_{R}(N)$
is uniquely determined by the induced maps $M_{1}F_{1}$ and $M_{2}F_{2}$.
\end{lem}

\section{\label{sec:Product-category} Product category}

Given a $2$-representation $\mathcal{V}$ for $A$ with $\mathcal{U}^{+}$-action
data $(E,x,\tau)$, we seek a $2$-representation for $C$ with data
$(\tilde{E},\tilde{x},\tilde{\tau})$ to serve as the tensor $2$-product
$\mathcal{L}(1)\ootimes\mathcal{V}$. In this section we describe
our proposal for the algebra $C$ and data $(\tilde{E},\tilde{x},\tilde{\tau})$,
and in the next section we study this data and verify that the nil
affine Hecke relations hold for $\tilde{x}$ and $\tilde{\tau}$.

\subsection{Naive product category}

\subsubsection{Naive product algebra $B$}
\begin{defn}
Let $B$ be the $k$-algebra:
\[
B=\begin{pmatrix}A[y] & E_{y}\\
0 & A[y]
\end{pmatrix}.
\]
Here the algebra structure of $B$ is given by matrix multiplication,
with the $(A[y],A[y])$-bimodule structure of $E_{y}$ contributing
for products with generators in $B_{12}$.
\end{defn}

A left $B$-module consists of a pair $\left(\begin{smallmatrix}M_{1}\\
M_{2}
\end{smallmatrix}\right)$ of left $A[y]$-modules, together with a morphism $\alpha:E_{y}\otimes_{A[y]}M_{2}\to M_{1}$
of left $A[y]$-modules. A right $B$-module is the data of a pair
$\left(\begin{smallmatrix}N_{1} & N_{2}\end{smallmatrix}\right)$
of right $A[y]$-modules, together with a morphism $\beta:N_{1}\otimes_{A[y]}E_{y}\to N_{2}$
of right $A[y]$-modules. It follows that a $(B,B)$-bimodule can
be written as a matrix of $(A[y],A[y])$-bimodules with accompanying
maps $\alpha$ and $\beta$ giving left and right actions of $E_{y}$.
Such a matrix with $\alpha,\beta$ determines a $(B,B)$-bimodule
only if the actions commute. Usually this commutativity is obvious
and we do not bother to check it.

A complex of left $B$-modules is the same data as a pair of complexes
of $A[y]$-modules together with a morphism $\alpha$ of complexes;
note that the differential of $E_{y}\otimes M_{2}$ for a complex
$(M_{2},d)$ is just $E_{y}\otimes d$. Similarly for right $B$-module
complexes. Complexes in this paper have a cohomological grading.

\subsubsection{\label{subsec:endofunctor-E'} Endofunctor $E'$ of $B\text{-cplx}$}
\begin{defn}
\label{def:E'} Let $E'$ be the complex of $(B,B)$-bimodules that
is nonzero in degrees $0$ and $1$, where it is given by: 
\[
E'_{0}=\begin{pmatrix}E[y] & E[y]E_{y}\\
0 & E[y]
\end{pmatrix},\;E'_{1}=\begin{pmatrix}E_{y} & E_{y}E_{y}\\
A[y] & E_{y}
\end{pmatrix}.
\]
The left action of a generator in $E_{y}\subset B$ is specified on
vertical columns of $E'_{0}$ by the maps $0:E_{y}\otimes0\to E[y]$
and $s:E_{y}E[y]\to E[y]E_{y}$. The left action on $E'_{1}$ is specified
by the identity map on vertical columns. The right action on $E'_{0}$
is specified by the identity on the top row and $0$ on the bottom
row. The right action on $E'_{1}$ is specified by the identity on
both rows. The differential $E'_{0}\to E'_{1}$ is given \emph{componentwise}
by $\left(\begin{smallmatrix}\pi & \pi E_{y}\\
0 & \pi
\end{smallmatrix}\right)$.
\end{defn}

Tensoring by $E'$ on the left gives an endofunctor $_{B}E'\otimes_{B}-$
of the category of complexes of $B$-modules. It is convenient to
have a formula for the action of this endofunctor on an arbitrary
complex of modules: 
\begin{lem}
Let $M=\left(\left(\begin{smallmatrix}M_{1}\\
M_{2}
\end{smallmatrix}\right),\alpha\right)$ be a complex of left $B$-modules (written as a pair of complexes),
where $\alpha:E_{y}\otimes_{A[y]}M_{2}\to M_{1}$ specifies the action
for generators in $E_{y}$. The functor $E'\otimes_{B}-$ on $M$
may be given by: 
\[
\left(\begin{pmatrix}M_{1}\\
M_{2}
\end{pmatrix},\alpha\right)\overset{E'}{\longmapsto}\left(\begin{pmatrix}E[y]M_{1}\overset{\overset{\pi M_{1}}{\curvearrowright}}{\oplus}E_{y}M_{1}[-1]\\
E[y]M_{2}\overset{\overset{\alpha\circ\pi M_{2}}{\curvearrowright}}{\oplus}M_{1}[-1]
\end{pmatrix},\begin{pmatrix}E[y]\alpha\circ sM_{2} & 0\\
0 & \idop_{E_{y}M_{1}}
\end{pmatrix}\right).
\]
Here the top and bottom rows express cocones of the maps $\pi M_{1}$
and $\alpha\circ\pi M_{2}$.
\end{lem}

\begin{rem}
It may help motivation to consider the effect of $E'$ at the level
of the Grothendieck group when $M_{1}$ and $M_{2}$ are just modules,
not complexes. The following discussion is not intended to be precise
or complete. 

Suppose $M_{1}'$ and $M_{2}'$ are projective left $A$-modules,
and $R_{1}$ and $R_{2}$ are projective left $k[y]$-modules. Consider
the projective left $A[y]$-modules $M_{1}=R_{1}\otimes_{k}M_{1}'$
and $M_{2}=R_{2}\otimes_{k}M_{2}'$. These are elements of the outer
product of categories $(k[y]\text{-proj})\boxtimes_{k}(A\text{-proj})$.
Suppose $\alpha:E_{y}M_{2}\to M_{1}$ is given. Apply $E'$ to $\left(\left(\begin{smallmatrix}M_{1}\\
M_{2}
\end{smallmatrix}\right),\alpha\right)$. The upper row is quasi-isomorphic to: 
\begin{align*}
\ker\bigl(E[y]M_{1}\xrightarrow{\pi M_{1}}E_{y}M_{1}\bigr) & \iso(y_{1}E[y])M_{1}\iso E[y]M_{1}\iso R_{1}\otimes_{k}(E\otimes_{A}M_{1}'),
\end{align*}
 where the first isomorphism follows by flatness of $M_{1}$. Letting
$e$ denote the action of $E$ on the Grothendieck group, we have
$(1\otimes e)\bigl([R_{1}]\otimes_{k}[M_{1}']\bigr)$ for the upper
row in the Grothendieck group. The lower row is the cocone of $\alpha$,
which contributes $[E[y]M_{2}]+[M_{1}]$ in the Grothendieck group.
Now recall that the raising functor for $\mathcal{L}(1)$ is just
$k[y]$. So: 
\[
M_{1}\iso(k[y]\otimes1)\bigl(R_{1}\otimes_{k}M_{1}'\bigr),\quad[M_{1}]=(e\otimes1)\bigl([R_{1}]\otimes_{k}[M_{1}']\bigr),
\]
 and we should interpret the copy of $M_{1}$ coming from the lower
row in this way, since the factor of $k[y]$ in the $A[y]\cong k[y]\otimes_{k}A$
of the lower left corner of $B$ is the higher weight copy. We also
have $[E[y]M_{2}]=(e\otimes1)\bigl([R_{2}]\otimes_{k}[M_{2}']\bigr)$.
Finally, it is a fact that $(e\otimes1)\bigl([R_{2}]\otimes_{k}[M_{2}']\bigr)=0$
because $\mathcal{L}(1)$ has only two weight categories. It follows
from these calculations that the action of $e'=[E']$ on the Grothendieck
group of the derived category has the form: 
\begin{align*}
e'[\left(\left(\begin{smallmatrix}M_{1}\\
M_{2}
\end{smallmatrix}\right),\alpha\right)] & :=[E'\left(\left(\begin{smallmatrix}M_{1}\\
M_{2}
\end{smallmatrix}\right),\alpha\right)]\\
 & =(e\otimes1+1\otimes e)\bigl([M_{1}']\otimes_{k}[R_{1}]+[M_{2}']\otimes_{k}[R_{2}]\bigr).
\end{align*}
 This agrees with the Hopf coproduct formula $\Delta(e)=e\otimes1+1\otimes e$.
\end{rem}

\begin{proof}[Proof of the lemma]
We first check that the matrix specifying the new $E_{y}$ action
gives a morphism of complexes. The diagonal coefficients of the matrix
give morphisms of the separate summands, and these commute with the
differentials on the separate summands. It remains to see that $\pi M_{1}\circ E[y]\alpha\circ sM_{2}=\idop_{E_{y}}M_{1}\circ E_{y}(\alpha\circ\pi M_{2})$,
and these agree because $\pi E_{y}\circ s=E_{y}\pi$.

Now we compute the tensor product following the recipe of §\ref{subsec:matrix-tensor-products}.
We have:
\[
E'\otimes_{B}M=\begin{pmatrix}\Bigl(E[y]M_{1}\oplus E[y]E_{y}M_{2}\Bigr)\big/I_{1}\overset{\overset{\pi M_{1}}{\curvearrowright}}{\oplus}\biggl(\Bigl(E_{y}M_{1}\oplus E_{y}E_{y}M_{2}\Bigr)\big/I_{1}'\biggr)[-1]\\
\Bigl(0\oplus E[y]M_{2}\Bigr)\big/I_{2}\overset{\overset{\alpha\circ\pi M_{2}}{\curvearrowright}}{\oplus}\biggl(\Bigl(A[y]M_{1}\oplus E_{y}M_{2}\Bigr)\big/I_{2}'\biggr)[-1]
\end{pmatrix}.
\]
Here the submodule $I_{1}$ is generated by all terms of the form
$e\otimes\alpha(e',m_{2})-e\otimes e'\otimes m_{2}$ for $e\in E[y],e'\in E_{y},m_{2}\in M_{2}$.
So every element of the quotient has a canonical representative in
$E[y]M_{1}$, and the quotient is isomorphic to $E[y]M_{1}$. With
analogous reasoning we see that the quotient by $I_{1}'$ is isomorphic
to $E_{y}M_{1}$, that by $I_{2}$ is isomorphic to $E[y]M_{2}$,
and that by $I_{2}'$ is isomorphic to $M_{1}$. The differential
may be written before taking quotients as $dM_{1}$ on the top and
$dM_{2}$ on the bottom. The images of $dM_{2}$ in $E_{y}M_{2}$
represent elements in $M_{1}$ by way of $\alpha$, and this determines
the differential component $\alpha\circ\pi M_{2}$ between summands
of the bottom row.

Now we calculate the new $E_{y}$ action in order to view this as
a complex of $B$-modules. Using the description of the left $B$-action
on $E'$, one sees that the action on the left summand is by $sM_{2}$,
which is represented in $E[y]M_{1}$ through $\alpha$, so the action
written on the quotients as described above is given by $E[y]\alpha\circ sM_{2}$.
The action is obvious on the right summand.
\end{proof}

\subsubsection{Category $\text{per }B$ and generator $X$}
\begin{defn}
Let $X$ be the following complex of $B$-modules:
\begin{align*}
X & =X_{1}\oplus X_{2}\\
X_{1} & =\begin{pmatrix}A[y]\\
0
\end{pmatrix}\\
X_{2} & =E'(X_{1})=\begin{pmatrix}E[y] & \overset{\pi}{\lra} & E_{y}\\
0 & \lra & A[y]
\end{pmatrix}
\end{align*}
where $X_{1}$ lies in degree $0$ and $X_{2}$ in degrees $0$ and
$1$. The $E_{y}$ action on $X_{2}$ is given by $E_{y}\otimes_{A[y]}A[y]\xrightarrow{\sim}E_{y}$,
$e\otimes1\mapsto e$.
\end{defn}

One can see that $X_{1}=Be_{1}$ and $X_{2}=E'e_{1}$, with $e_{i}\in B$
the standard matrix idempotent. Observe that there is a canonical
right $A[y]$ action on $Be_{i}$ and on $X_{i}$ given componentwise.
\begin{prop}
The complex $X$ is strictly perfect and generates $\text{per }B$.
\end{prop}

\begin{proof}
We can write $X$ in terms of $B$: 
\begin{align*}
X_{1} & =Be_{1}\\
X_{2} & =Be_{1}\otimes_{A}E\to Be_{2},
\end{align*}
 where the differential is by $\pi$ on the upper row. This is a complex
of finitely generated projective $B$-modules because $_{A}E$ is
finitely generated and projective. So $X$ is strictly perfect. To
see that $X$ generates $\text{per }B$, first note that $Be_{1}=X_{1}\in\langle X\rangle_{\Delta}$.
Now consider $Be_{1}\otimes_{A}E$ as a complex in degree $0$. There
is a map of complexes $X_{2}\to Be_{1}\otimes_{A}E$ given by the
identity in degree $0$ and by $0$ in degree $1$. Then $Be_{2}[-1]$
(a complex in degree $1$) is quasi-isomorphic to the cocone of this
map. So $Be_{2}\in\langle X\rangle_{\Delta}$.
\end{proof}
Recall our notation $\pi_{i}=E^{n-i}\pi E^{i-1}:E^{n}[y]\to E^{n-i}E_{y}E^{i-1}$.
\begin{lem}
The kernel of $\varphi:E^{n}[y]\xrightarrow{(\pi_{i})_{i}}\bigoplus_{i=1}^{n}E^{n-i}E_{y}E^{i-1}$
is $(y_{1}\dots y_{n})E^{n}[y]$.\label{lem:kernels-lemma}
\end{lem}

\begin{proof}
We have assumed that $E^{n}$ is free as a $P_{n}$-module. It follows
that $E^{n}[y]$ is free as a $P_{n}[y]$-module. Let $e\in\ker\varphi$.
So $\pi_{i}(e)=0$ and therefore $e\in y_{i}E^{n}[y]$ for each $i\in\{1,\dots,n\}$.
Let $B$ be a basis of $E^{n}[y]$ over $P_{n}[y]$. Write 
\[
e=y_{i}\sum_{j=1}^{\ell}f_{j}^{i}(x_{1},\dots,x_{n},y)\cdot b_{j}
\]
 for $b_{j}\in B$ distinct and $f_{j}^{i}\in P_{n}[y]$. It follows
that $y_{i}f_{j}^{i}=y_{k}f_{j}^{k}$ in $P_{n}[y]$ for each $(i,k)\in\{1,\dots,n\}^{\times2}$
and $j\in\{1,\dots,\ell\}$. Then $e=y_{1}\dots y_{n}e^{\circ}$ for
some $e^{\circ}\in E^{n}[y]$ because $P_{n}[y]$ is a unique factorization
domain and each $y_{i}$ is irreducible.
\end{proof}
\begin{lem}
\label{lem:E'(X_2)} The complex $E'X_{2}$ is concentrated in degrees
$0$, $1$, and $2$: 
\[
E'X_{2}=\left(\begin{pmatrix}E^{2}[y]\xrightarrow{(\pi_{2},\pi_{1})}E_{y}E\oplus EE_{y}\xrightarrow{(-\pi_{1},\pi_{2})}E_{y}E_{y}\\
0\longrightarrow E[y]\oplus E[y]\xrightarrow{(-\pi,\pi)}E_{y}
\end{pmatrix},\alpha\right),
\]
where 
\begin{align*}
\alpha_{0} & =0\\
\alpha_{1} & =\left(\begin{smallmatrix}\idop_{E_{y}E} & 0\\
0 & s
\end{smallmatrix}\right)\\
\alpha_{2} & =\idop_{E_{y}E_{y}}.
\end{align*}
\end{lem}

\begin{proof}
Computation. The minus signs arise from shifting differentials.
\end{proof}
\begin{prop}
\label{prop:E'(X)-perfect} The complex $E'X$ is quasi-isomorphic
to a finite direct sum of summands of $X$.
\end{prop}

We define two complexes of $B$-modules before proving the proposition.
\begin{defn}
Let $R,X_{2}'\in B\text{-cplx}$ be given by: 
\[
R=\begin{pmatrix}E^{2}[y]\xrightarrow{\left(\begin{smallmatrix}\pi E\\
\pi E\ci\tau
\end{smallmatrix}\right)}E_{y}E\oplus E_{y}E\\
0\to E[y]\oplus E[y]
\end{pmatrix},
\]
 
\[
X_{2}'=\begin{pmatrix}\tau y_{1}E^{2}[y] & \overset{\pi E}{\lra} & E_{y}E\\
0 & \lra & E[y]
\end{pmatrix},
\]
 both lying in degrees $0$ and $1$ (cohomological grading), and
the $E_{y}$ action on $R$ specified by $0$ and the canonical map
\[
E_{y}\otimes(E[y]\oplus E[y])\to E_{y}E\oplus E_{y}E,
\]
 and on $X_{2}'$ specified by $0$ and $\idop_{E_{y}E[y]}$.
\end{defn}

\begin{lem}
We have that $X_{2}'$ is a finite direct sum of summands of $X_{2}$,
and hence of $X$.
\end{lem}

\begin{proof}
Observe first that $X_{2}\otimes_{A}E$ is a finite direct sum of
summands of $X$ because $_{A}E$ is finitely generated projective.
(Here we use the componentwise right $A$-action on $X_{2}$.) Using
the formulas 
\begin{align*}
\pi_{2}\circ\delta & =\pi_{2},\\
\pi_{2}\circ(1-\delta) & =0,
\end{align*}
and $\delta\cdot(1-\delta)=0$, one has the decomposition of $X_{2}\otimes_{A}E$:
\begin{align*}
X_{2}\otimes_{A}E & =\begin{pmatrix}E^{2}[y] & \overset{\pi_{2}}{\lra} & E_{y}E\\
0 & \lra & E[y]
\end{pmatrix}\\
 & =\begin{pmatrix}\delta\cdot E^{2}[y] & \overset{\pi_{2}}{\lra} & E_{y}E\\
0 & \lra & E[y]
\end{pmatrix}\oplus\begin{pmatrix}(1-\delta)\cdot E^{2}[y]\\
0
\end{pmatrix}.
\end{align*}
\end{proof}
The nil-affine Hecke algebra is an $n!\times n!$ matrix algebra over
the symmetric polynomials $P_{n}^{S_{n}}$ (cf.~Prop.~3.4 of \cite{rouquier2KacMoodyAlgebras2008}).
This provides a decomposition of $E^{n}$ into submodules called \textquoteleft divided
powers\textquoteright : 
\[
E^{n}\iso\overset{n\text{ copies}}{\overbrace{E^{(n)}\oplus\dots\oplus E^{(n)}}}.
\]
 When $n=2$, the isomorphism is given explicitly as follows (by extension
to left $A[y]$-modules): 
\begin{equation}
E^{2}[y]\xrightarrow[\left(\begin{smallmatrix}\tau y_{1}\\
\tau
\end{smallmatrix}\right)]{\sim}\tau y_{1}E^{2}[y]\oplus\tau y_{1}E^{2}[y].\label{eq:divided-squares-iso}
\end{equation}
 The inverse map is $(\iota,-y_{2})$, where $\iota:\tau y_{1}E^{2}[y]\hookrightarrow E^{2}[y]$
is the inclusion. The elements $\tau y_{1}$ and $-y_{2}\tau$ are
orthogonal idempotents summing to $\idop$, and $\tau$ becomes an
isomorphism from $-y_{2}\tau E^{2}[y]$ to $\tau y_{1}E^{2}[y]$.
(To verify this: $\tau(-y_{2}\tau)ee=\tau ee=\tau y_{1}(\tau ee)\in\tau y_{1}E^{2}[y]$,
$\tau(-y_{2}\tau)y_{1}ee=\tau y_{1}ee$, and if $\tau(-y_{2}\tau)ee=0$
then $\tau ee=0$ so $-y_{2}\tau ee=0$.)
\begin{lem}
There is an isomorphism $R\iso X_{2}'\oplus X_{2}'$ in $B\text{-cplx}$
given by the above isomorphism on the degree $0$ term of the upper
row, and the identity on all other terms. So $R$ is a finite direct
sum of summands of $X_{2}$, and hence of $X$. In particular, $R$
is strictly perfect.
\end{lem}

\begin{lem}
\label{lem:R-quasi-iso-E'X} There is a quasi-isomorphism $R\xrightarrow{q.i.}E'X_{2}$
determined by $\idop_{E^{2}[y]}$ on the degree $0$ term of the upper
row and $\left(\begin{smallmatrix}1 & 0\\
1 & -y_{1}
\end{smallmatrix}\right)$ on the degree $1$ term of the lower row.
\end{lem}

\begin{proof}
We first check that the map is a morphism in $B\text{-cplx}$. The
matrix of the morphism on the degree $1$ part of the upper row, as
determined by equivariance over generators of $B$ in $E_{y}$, is
given by $\left(\begin{smallmatrix}\idop & 0\\
s & s\ci(x_{2}-x_{1})
\end{smallmatrix}\right)$. Observe that: 
\begin{align*}
 & \idop\ci\pi_{2}+0\ci\pi_{2}\ci\tau=\pi_{2};\\
 & s\ci\pi_{2}+s\ci(x_{2}-x_{1})\ci\pi_{2}\ci\tau\\
 & =\pi_{1}\ci s+(x_{1}-x_{2})\ci s\ci\pi_{2}\ci\tau\\
 & =\pi_{1}\ci s+\pi_{1}\ci(x_{1}-x_{2})\ci s\ci\tau\\
 & =\pi_{1}\ci\Bigl((x_{2}-x_{1})\ci\tau+\idop\\
 & \qquad+(x_{1}-x_{2})\ci\bigl((x_{2}-x_{1})\ci\tau+\idop\bigr)\ci\tau\Bigr)=\pi_{1}.
\end{align*}
 This shows compatibility with the differential from degree $0$ in
the upper row. The other compatibility checks are easier.

Now we show that the map is a quasi-isomorphism. The lower row of
$E'X_{2}$ has $H^{1}$ given by: 
\[
\{(e_{1},e_{2})\in E[y]^{\oplus2}\mid e_{1}-e_{2}=y_{1}e\text{ for some }e\in E[y]\}.
\]
 This is also the image of the (injective) map from $R$ in degree
$1$ of the lower row. The upper row of $E'X_{2}$ has $H^{0}=\ker(d^{0})=y_{1}y_{2}E^{2}[y]$
by Lemma \ref{lem:kernels-lemma}. The cohomology of the upper row
of $R$ is computed as follows. We have an isomorphism: 
\[
E^{2}[y]\iso\tau y_{1}E^{2}[y]\oplus-y_{2}\tau E^{2}[y].
\]
 Notice that $\pi_{2}\ci\tau$ vanishes on the first summand, and
$\pi_{2}$ vanishes on the second. Then one may compute: 
\[
\ker\bigl(\tau y_{1}E^{2}[y]\xrightarrow{\pi_{2}}E_{y}E\bigr)=\tau y_{1}y_{2}E^{2}[y]\subset y_{1}y_{2}E^{2}[y]
\]
 and 
\[
\ker\bigl(-y_{2}\tau E^{2}[y]\xrightarrow[\sim]{\tau}\tau y_{1}E^{2}[y]\xrightarrow{\pi_{2}}E_{y}E\bigr)=-y_{2}\tau y_{1}y_{2}E^{2}[y]\subset y_{1}y_{2}E^{2}[y].
\]
 So 
\[
\ker\bigl(\left(\begin{smallmatrix}\pi_{2}\\
\pi_{2}\ci\tau
\end{smallmatrix}\right)\bigr)\subset y_{1}y_{2}E^{2}[y].
\]
 The reverse inclusion is obvious, so $H^{0}$ of the upper row is
$y_{1}y_{2}E^{2}[y]$. This shows that $\idop_{E^{2}[y]}$ induces
an isomorphism on homology in degree $0$ of the upper row. Using
the decomposition and inspecting the maps above, we also see that
$d^{0}$ on the upper row of $R$ is surjective. Finally we consider
$H^{1}$ of the upper row of $E'X_{2}$ and show it is zero. (Clearly
the $H^{2}$ is zero.) Let $(ee_{1},ee_{2})\in E_{y}E\oplus EE_{y}$
be in $\ker(d^{1})$, i.e.~such that $\pi_{1}(ee_{1})=\pi_{2}(ee_{2})$.
Then $ee_{1}=ee_{2}+(Ex-xE)ee^{\circ}$ for some $ee^{\circ}\in E^{2}$.
(Note that $E_{y}E_{y}\cong E^{2}\big/(Ex-xE)$ where $y$ acts by
$Ex$ or $xE$.) Then consider $ee_{2}+(Ex-y)ee^{\circ}\in E^{2}[y]$.
The differential $d^{0}$ sends this to $ee_{1}$ in $E_{y}E$ and
to $ee_{2}$ in $EE_{y}$.
\end{proof}
\begin{proof}[Proof of Proposition \ref{prop:E'(X)-perfect}]
 The proposition follows from the preceding three lemmas.
\end{proof}
\begin{cor}
Tensoring with $_{B}E'_{B}$ gives an endofunctor $E'\otimes_{B}-$
of $\text{per }B$.
\end{cor}

\begin{proof}
We know that $X\in\text{per }B$, and it follows from Prop.~\ref{prop:E'(X)-perfect}
that $E'\otimes_{B}X\in\text{per }B$. The corollary follows because
$X$ generates $\text{per }B$.
\end{proof}
\begin{rem}
We do not know that $E'\otimes_{B}-$ on $K^{b}(B)$ is exact, so
we do not know that it descends to an endofunctor defined on all of
$D^{b}(B)$.
\end{rem}

\subsection{Bimodules $G_{n}$}

The constructions of this paper make use of certain bimodules that
we describe next.
\begin{defn}
Let $G_{n}$ denote $\Hom_{K^{b}(B)}(X_{2},E'^{n}X_{1})$.
\end{defn}

Every $G_{n}$ has the structure of $(G_{1}^{\op},A[y])$-bimodule
by pre- and post-composition. Here we understand $A[y]\cong\End_{K^{b}(B)}(X_{1})^{\op}$
and use functoriality of $E'$ for the action. Note that $G_{1}=\Hom_{K^{b}(B)}(X_{2},X_{2})$
has an algebra structure, and the right regular action of $G_{1}^{\op}$
on $G_{1}$ extends the right $A[y]$ action.

In this section we gather some facts regarding these bimodules and
give concrete presentations in small cases that are easier to handle.
Given $n\in\{1,2,3,4\}$, we define $\bar{G}_{n}$ as an $(A[y],A[y])$-sub-bimodule
of 
\[
E^{n-1}[y]^{\oplus n}\oplus\Hom_{A}(_{A}E,E^{n})[y].
\]
 (By $E^{0}[y]$ we mean $A[y]$.) We give isomorphisms $\bar{G}_{n}\iso G_{n}$
for such $n$. These isomorphisms induce left $G_{1}^{\op}$-actions
on $\bar{G}_{n}$ that extend the left $A[y]$-actions. In future
sections we do not distinguish $G_{n}$ from $\bar{G}_{n}$ and write
only the former. For example, we will write $(\theta,\varphi)$ to
indicate an element of $G_{1}$ or $G_{1}^{\op}$.
\begin{defn}
Define the following $(A[y],A[y])$-sub-bimodule of $A^{\op}[y]\oplus\End_{A}(_{A}E)[y]$:
\begin{gather*}
\bar{G}_{1}=\biggl\langle(\theta,\varphi)\in A^{\op}[y]\oplus\End_{A}(_{A}E)[y]\biggr|\\
\begin{split}\varphi & =\_.\theta+y_{1}\varphi_{1}\\
 & \text{ for some }\varphi_{1}\in\End_{A}(_{A}E)[y]\biggr\rangle.
\end{split}
\end{gather*}
Here $\_.\theta\in\mathrm{End}_{A}(_{A}E)[y]$ is the morphism sending
$e$ to $e.\theta$. One easily checks that $\bar{G}_{1}$ is closed
under the bimodule actions. This bimodule also has a $k$-algebra
structure with componentwise multiplication (using the opposite multiplication
on generators in $A[y]$).
\end{defn}

Note that $\bar{G}_{1}$ contains a copy of $A^{\op}[y]$, namely
the subspace of elements of the form $(\theta,\_.\theta)$.
\begin{prop}
\label{prop:G_1-bar} There is an isomorphism of $(A[y],A[y])$-bimodules
$\bar{G}_{1}\iso G_{1}$ determined by equivariance over $E_{y}\subset B$
with the following data: 
\[
(\theta,\varphi)\mapsto\left(\begin{pmatrix}\left(e,0\right)\\
(0,1)
\end{pmatrix}\mapsto\begin{pmatrix}\left(\varphi(e),0\right)\\
(0,\theta)
\end{pmatrix}\right).
\]
 Here $(e,0)\in E[y]\oplus E_{y}$ is an element of the upper row
of $X_{2}$, with $e$ in degree $0$ and $0$ in degree $1$. Analogously
with the lower row. This isomorphism respects the $k$-algebra structure.
\end{prop}

\begin{proof}
The condition $\varphi=\_.\theta+y_{1}\varphi_{1}$ in the definition
of $\bar{G}_{1}$ is equivalent to the statement that the morphism
given as the image of $(\theta,\varphi)$ defined in the proposition
has zero differential.
\end{proof}
\begin{defn}
\label{def:G_2} Define the following $(A[y],A[y])$-sub-bimodule
of $E[y]^{\oplus2}\oplus\Hom_{A}(_{A}E,E^{2})[y]$: 
\begin{gather*}
\bar{G}_{2}=\biggl\langle(e_{1},e_{2},\xi)\in E[y]^{\oplus2}\oplus\Hom_{A}(_{A}E,E^{2})[y]\biggr|\\
\begin{split}\qquad e_{1}-e_{2} & =y_{1}e'\\
\xi & =\_\otimes e_{1}+y_{2}\xi_{1}\\
 & =\delta(\_\otimes e_{2})+y_{1}\xi_{2}\\
 & \text{ for some }e'\in E[y]\text{ and }\xi_{\ell}\in\Hom_{A}(_{A}E,E^{2})[y]\biggr\rangle.
\end{split}
\end{gather*}
 Here $\_\otimes e_{\ell}\in\mathrm{Hom}_{A}(_{A}E,E^{2})[y]$ sends
$e$ to $e\otimes e_{\ell}$. One quickly checks that the conditions
are preserved by the bimodule actions.
\end{defn}

\begin{prop}
There is an isomorphism of $(A[y],A[y])$-bimodules $\bar{G}_{2}\iso G_{2}$
determined by equivariance over $E_{y}\subset B$ with the following
data:
\[
(e_{1},e_{2},\xi)\mapsto\left(\begin{pmatrix}(e,0)\\
(0,1)
\end{pmatrix}\mapsto\begin{pmatrix}(\xi(e),0,0)\\
(0,\left(\begin{smallmatrix}e_{1}\\
e_{2}
\end{smallmatrix}\right),0)
\end{pmatrix}\right).
\]
\end{prop}

\begin{proof}
Use the description of $E'X_{2}$ in Lemma \ref{lem:E'(X_2)}. As
in Prop.~\ref{prop:G_1-bar}, the condition of the definition of
$\bar{G}_{2}$ is equivalent to the statement that the image of $(e_{1},e_{2},\xi)$
has zero differential.
\end{proof}
In order to parametrize $G_{3}$, we compute the components of $E'^{2}X_{2}=E'^{3}X_{1}$
in degrees $0$, $1$, and $2$: 
\[
\begin{pmatrix}E^{3}[y] & \to & E_{y}EE\oplus EE_{y}E\oplus EEE_{y} & \to & E_{y}E_{y}E\oplus E_{y}EE_{y}\oplus EE_{y}E_{y}\to\dots\\
0 & \to & E^{2}[y]\oplus E^{2}[y]\oplus E^{2}[y] & \to & E_{y}E\oplus EE_{y}\oplus EE_{y}\to\dots
\end{pmatrix}.
\]
 The upper left differential map is $(\pi_{3},\pi_{2},\pi_{1})$.
We don\textquoteright t make use of the upper right. The bottom right
differential map is given by the matrix: 
\[
\begin{pmatrix}-\pi_{2} & \pi_{2} & 0\\
-\pi_{1} & 0 & \pi_{1}\circ\delta\\
0 & -\pi_{1} & \pi_{1}
\end{pmatrix}.
\]

\begin{defn}
\label{def:G_3} Define the following $(A[y],A[y])$-sub-bimodule
of $E^{2}[y]^{\oplus3}\oplus\Hom_{A}(_{A}E,E^{3})[y]$: 
\begin{gather*}
\bar{G}_{3}=\biggl\langle(ee_{1},ee_{2},ee_{3},\chi)\in E^{2}[y]^{\oplus3}\oplus\Hom_{A}(_{A}E,E^{3})[y]\biggr|\\
\begin{split}\qquad ee_{1}-ee_{2} & =y_{2}ee'\\
ee_{3}-ee_{2} & =y_{1}ee''\\
\delta(ee_{3})-ee_{1} & =y_{1}ee''',\\
\chi & =\_\otimes ee_{1}+y_{3}\chi_{1}\\
 & =\delta E(\_\otimes ee_{2})+y_{2}\chi_{2}\\
 & =E\delta\circ\delta E(\_\otimes ee_{3})+y_{1}\chi_{3}\\
 & \text{ for some }ee^{k}\in E^{2}[y]\text{ and }\chi_{\ell}\in\Hom_{A}(_{A}E,E^{3})[y]\biggr\rangle.
\end{split}
\end{gather*}
Here $\_\otimes ee_{\ell}\in\mathrm{Hom}_{A}(_{A}E,E^{3})[y]$ sends
$e$ to $e\otimes ee_{\ell}$. One quickly checks that the conditions
are preserved by the bimodule actions.
\end{defn}

\begin{prop}
There is an isomorphism of $(A[y],A[y])$-bimodules $\bar{G}_{3}\iso G_{3}$
determined by equivariance over $E_{y}\subset B$ with the following
data:
\[
(ee_{1},ee_{2},ee_{3},\chi)\mapsto\left(\begin{pmatrix}(e,0)\\
(0,1)
\end{pmatrix}\mapsto\begin{pmatrix}(\chi(e),0,\dots)\\
(0,\left(\begin{smallmatrix}ee_{1}\\
ee_{2}\\
ee_{3}
\end{smallmatrix}\right),\dots)
\end{pmatrix}\right).
\]
\end{prop}

\begin{proof}
The condition of the definition of $\bar{G}_{3}$ is equivalent to
the statement that the image of $(ee_{1},ee_{2},ee_{3},\chi)$ has
zero differential.
\end{proof}
\begin{defn}
\label{def:G_4} Define the following $(A[y],A[y])$-sub-bimodule
of $E^{3}[y]^{\oplus4}\oplus\Hom_{A}(_{A}E,E^{4})[y]$: 
\begin{gather*}
\bar{G}_{4}=\biggl\langle(eee_{1},eee_{2},eee_{3},eee_{4},\psi)\in E^{3}[y]^{\oplus4}\oplus\Hom_{A}(_{A}E,E^{4})[y]\biggr|\\
\begin{split}\qquad eee_{3}-eee_{4} & =y_{1}eee^{(1)}\\
eee_{2}-eee_{3} & =y_{2}eee^{(2)}\\
E\delta(eee_{4})-eee_{2} & =y_{1}eee^{(3)}\\
eee_{1}-eee_{2} & =y_{3}eee^{(4)}\\
eee_{1}-\delta E(eee_{3}) & =y_{2}eee^{(5)}\\
eee_{1}-\delta E\circ E\delta(eee_{1}) & =y_{1}eee^{(6)}\\
\psi & =\_\otimes eee_{1}+y_{4}\psi_{1}\\
 & =\delta E^{2}(\_\otimes eee_{2})+y_{3}\psi_{2}\\
 & =E\delta E\circ\delta E^{2}(\_\otimes eee_{3})+y_{2}\psi_{3}\\
 & =E^{2}\delta\circ E\delta E\circ\delta E^{2}(\_\otimes eee_{4})+y_{1}\psi_{4}\\
 & \text{ for some }eee^{k}\in E^{3}[y]\text{ and }\psi_{\ell}\in\Hom_{A}(_{A}E,E^{4})[y]\biggr\rangle.
\end{split}
\end{gather*}
 Here $\_\otimes eee_{\ell}\in\mathrm{Hom}_{A}(_{A}E,E^{4})[y]$ sends
$e$ to $e\otimes eee_{\ell}$. One quickly checks that the conditions
are preserved by the bimodule actions.
\end{defn}

\begin{lem}
Under the conditions on $eee_{i}$ in the definition, there is a unique
$\overline{eee}\in E^{3}[y]$ such that: 
\begin{align*}
eee^{(5)}-eee^{(2)} & =y_{3}\overline{eee},\\
eee^{(4)}-\tau E(eee_{3}) & =y_{2}\overline{eee}.
\end{align*}
\end{lem}

\begin{proof}
Subtracting two equations from those conditions: 
\begin{align*}
y_{2}\left(eee^{(5)}-eee^{(2)}\right) & =eee_{1}-eee_{2}-y_{3}\tau E(eee_{3})\\
 & =y_{3}\left(eee^{(4)}-\tau E(eee_{3})\right)
\end{align*}
 By Lemma \ref{lem:kernels-lemma} we know there is some $\overline{eee}$
satisfying the claim. It is unique because the $y_{i}$ are injective.
\end{proof}
\begin{prop}
There is an isomorphism of $(A[y],A[y])$-bimodules $\bar{G}_{4}\iso G_{4}$
determined by equivariance over $E_{y}\subset B$ with the following
data:
\[
(eee_{1},eee_{2},eee_{3},eee_{4},\psi)\mapsto\left(\begin{pmatrix}(e,0)\\
(0,1)
\end{pmatrix}\mapsto\begin{pmatrix}(\psi(e),0,\dots)\\
(0,\left(\begin{smallmatrix}eee_{1}\\
eee_{2}\\
eee_{3}\\
eee_{4}
\end{smallmatrix}\right),\dots)
\end{pmatrix}\right).
\]
\end{prop}

\begin{proof}
The reader may compute the first terms of $E'^{4}X_{1}$ and show
that the condition of the definition of $\bar{G}_{4}$ is equivalent
to the statement that the image of $(eee_{1},eee_{2},eee_{3},eee_{4},\psi)$
defined in the proposition has zero differential. There is some ambiguity
in the order of summands in degree $1$ of the lower row. The convention
we have used is that the first summand arises from the latest application
of $E'$ which moves a term from degree $0$ of the upper row to degree
$1$ of the lower (and increments the exponents on existing terms
in the lower row).
\end{proof}
It will be useful to describe alternative, equivalent, conditions
defining $\bar{G}_{2}$ and $\bar{G}_{3}$. It is sometimes easier
to work with them.
\begin{prop}
\label{prop:K-conds} Given $(e_{1},e_{2},\xi)\in E[y]^{\oplus2}\oplus\Hom_{A}(_{A}E,E^{2})[y]$
with $e_{1}-e_{2}=y_{1}e'$, the following conditions are equivalent:
\[
\begin{split}\xi & =\_\otimes e_{1}+y_{2}\xi_{1}\\
 & =\delta(\_\otimes e_{2})+y_{1}\xi_{2}\\
 & \text{ for some }\xi_{\ell}\in\Hom_{A}(_{A}E,E^{2})[y]
\end{split}
\]
and 
\[
\begin{split}\xi & =\_\otimes e_{1}+y_{2}\xi_{1}\\
\xi_{1} & =\tau(\_\otimes e_{2})+y_{1}\xi'\\
 & \text{ for some }\xi'\in\Hom_{A}(_{A}E,E^{2})[y].
\end{split}
\]
 When these conditions hold, $\xi_{1}$, $\xi_{2}$, and $\xi'$ are
uniquely determined by the data $(e_{1},e_{2},\xi)$, and $\xi_{2}=\_\otimes e'+y_{2}\xi'$.
\end{prop}

\begin{proof}
Suppose the first condition holds. Using $\delta=y_{2}\tau+\idop$
and $e_{1}-e_{2}=(x-y)e'$, we can rearrange the first equality: 
\[
\_\otimes e_{1}+y_{2}\xi_{1}=y_{1}\xi_{2}+y_{2}\tau(\_\otimes e_{2})+\_\otimes e_{2},
\]
from which
\[
y_{2}\Bigl(\xi_{1}-\tau(\_\otimes e_{2})\Bigr)=y_{1}\Bigl(\xi_{2}-\_\otimes e'\Bigr).
\]
By Lemma \ref{lem:kernels-lemma}, the image of $\xi_{1}-\tau(\_\otimes e_{2})$
is in $y_{1}y_{2}E^{2}[y]$. Recall that $y_{1}$ is injective. We
can then make the following definition: 
\[
\xi'=y_{1}^{-1}(\xi_{1}-\tau(\_\otimes e_{2})).
\]
 The second condition and the final claim follow from this.

Starting now with the second condition, plugging the second equation
into the first, we find:
\begin{align*}
\xi & =\_\otimes e_{1}+y_{2}\bigl(\tau(\_\otimes e_{2})+y_{1}\xi'\bigr)\\
 & =\delta(\_\otimes e_{2})+\_\otimes(e_{1}-e_{2})+y_{2}y_{1}\xi'\\
 & =\delta(\_\otimes e_{2})+y_{1}\bigl(\_\otimes e'+y_{2}\xi'\bigr).
\end{align*}
 This is the second line of the first condition, and it establishes
the final claim.

The uniqueness claims are clear.
\end{proof}
\begin{prop}
\label{prop:L-conds} Given $(ee_{1},ee_{2},ee_{3},\chi)\in E^{2}[y]^{\oplus3}\oplus\Hom_{A}(_{A}E,E^{3})[y]$
with 
\begin{align}
ee_{1}-ee_{2} & =y_{2}ee'\label{eq:Relns_ee^k}\\
ee_{3}-ee_{2} & =y_{1}ee''\\
\delta(ee_{3})-ee_{1} & =y_{1}ee''',\label{eq:Relns_ee^k-last}
\end{align}
 the following conditions are equivalent:
\begin{align*}
\chi & =\_\otimes ee_{1}+y_{3}\chi_{1}\\
 & =\delta E(\_\otimes ee_{2})+y_{2}\chi_{2}\\
 & =E\delta\circ\delta E(\_\otimes ee_{3})+y_{1}\chi_{3}\\
 & \text{ for some }\chi_{\ell}\in\Hom_{A}(_{A}E,E^{3})[y]
\end{align*}
and 
\begin{align*}
\chi & =\_\otimes ee_{1}+y_{3}\chi_{1}\\
\chi_{1} & =\tau E(\_\otimes ee_{2})+y_{2}\chi_{1}'\\
\chi_{1}' & =E\tau\circ\tau E(\_\otimes ee_{3})+y_{1}\chi''\\
 & \text{ for some }\chi''\in\Hom_{A}(_{A}E,E^{3})[y].
\end{align*}
 When the conditions hold, $\chi_{1}$, $\chi_{2}$, $\chi_{3}$,
$\chi_{1}'$, and $\chi''$ are uniquely determined by the data $(ee_{1},ee_{2},ee_{3},\chi)$,
and there is a unique $\overline{ee}\in E^{2}[y]$ such that 
\begin{align*}
\tau(ee_{3})-ee' & =y_{1}\overline{ee}\\
ee'''-ee'' & =y_{2}\overline{ee}.
\end{align*}
 Define the map $\chi_{2}'=-\_\otimes\overline{ee}+y_{3}\chi''$.
Then we also have 
\[
\chi_{2}=E\tau\circ\delta E(\_\otimes ee_{3})+y_{1}\chi_{2}'
\]
 and 
\[
\chi_{3}=-\delta E(\_\otimes ee'')+y_{2}\chi_{2}'.
\]
Assuming $\chi=\_\otimes ee_{1}+y_{3}\chi_{1}$, the other two conditions
together are equivalent to a single condition on $\chi_{1}$:
\[
\chi_{1}=-\tau Ey_{1}(\_\otimes ee'')+E\delta\circ\tau E(\_\otimes ee_{3})+y_{2}y_{1}\chi''.
\]
\end{prop}

\begin{proof}
Suppose the first condition holds. Equating the first two formulas
for $\chi$ in the first condition and using $\delta E=y_{3}\tau E+\idop$
gives:
\begin{align*}
\_\otimes ee_{1}+y_{3}\chi_{1} & =y_{3}\tau E(\_\otimes ee_{2})+\_\otimes ee_{2}+y_{2}\chi_{2}
\end{align*}
 thus 
\[
y_{3}\bigl(\chi_{1}-\tau E(\_\otimes ee_{2})\bigr)=y_{2}\bigl(\chi_{2}-\_\otimes ee'\bigr).
\]
 By Lemma \ref{lem:kernels-lemma} again, the image of this function
lies in $y_{2}y_{3}E^{3}[y]$, and since each $y_{i}$ is injective,
we can define a new function $\chi_{1}'$ such that: 
\begin{align*}
\chi_{1} & =\tau E(\_\otimes ee_{2})+y_{2}\chi_{1}'\\
\chi_{2} & =\_\otimes ee'+y_{3}\chi_{1}'.
\end{align*}
Equating now the second and third formulas, we have:
\begin{align*}
y_{2}E\tau\circ\delta E(\_\otimes ee_{3})+\delta E(\_\otimes ee_{3})+y_{1}\chi_{3} & =\delta E(\_\otimes ee_{2})+y_{2}\chi_{2}
\end{align*}
 so 
\[
y_{2}\bigl(\chi_{2}-E\tau\circ\delta E(\_\otimes ee_{3})\bigr)=y_{1}\bigl(\chi_{3}+\delta E(\_\otimes ee'')\bigr),
\]
so for some $\chi_{2}'$ we can write:
\begin{align*}
\chi_{2} & =E\tau\circ\delta E(\_\otimes ee_{3})+y_{1}\chi_{2}'\\
\chi_{3} & =-\delta E(\_\otimes ee'')+y_{2}\chi_{2}'.
\end{align*}
We will need a fact derived from the relations (\ref{eq:Relns_ee^k})--(\ref{eq:Relns_ee^k-last})
of the $ee^{k}$. Adding the first and third relations and subtracting
the second yields 
\[
y_{1}\bigl(ee'''-ee''\bigr)=y_{2}\bigl(\tau(ee_{3})-ee'\bigr),
\]

from which we see there must be a (unique) $\overline{ee}$ with 
\begin{align*}
\tau(ee_{3})-ee' & =y_{1}\overline{ee}\\
ee'''-ee'' & =y_{2}\overline{ee}.
\end{align*}
 This gives the third claim of the proposition.

Equating now the two formulas we derived for $\chi_{2}$:
\begin{align*}
y_{3}E\tau\circ\tau E(\_\otimes ee_{3})+E\tau(\_\otimes ee_{3})+y_{1}\chi_{2}' & =\_\otimes ee'+y_{3}\chi_{1}'
\end{align*}
 so 
\[
y_{3}\bigl(\chi_{1}'-E\tau\circ\tau E(\_\otimes ee_{3})\bigr)=y_{1}\bigl(\chi_{2}'+\_\otimes\overline{ee}\bigr).
\]
 Therefore
\begin{align*}
\chi_{1}' & =E\tau\circ\tau E(\_\otimes ee_{3})+y_{1}\chi''\\
\chi_{2}' & =-\_\otimes\overline{ee}+y_{3}\chi''
\end{align*}
 for some $\chi''$, as desired.

In the reverse direction, starting with the second condition, plugging
the $\chi_{1}$ and $\chi_{1}'$ formulas into the first $\chi$ formula
gives:
\[
\chi=\_\otimes ee_{1}+y_{3}\Bigl(\tau E(\_\otimes ee_{2})+y_{2}\bigl(E\tau\circ\tau E(\_\otimes ee_{3})+y_{1}\chi''\bigr)\Bigr),
\]
so
\begin{align*}
\chi-\delta E(\_\otimes ee_{2}) & =\_\otimes\left(ee_{1}-ee_{2}\right)\\
 & \qquad+y_{2}\bigl(E\tau\circ\tau E(\_\otimes ee_{3})+y_{1}\chi''\bigr)\\
 & =y_{2}\Bigl(\_\otimes ee'+E\tau\circ\tau E(\_\otimes ee_{3})+y_{1}\chi''\Bigr),
\end{align*}
 as desired. Similarly:
\begin{align*}
\chi-E\delta\circ\delta E(\_\otimes ee_{3}) & =\chi-y_{3}y_{2}E\tau\circ\tau E(\_\otimes ee_{3})\\
 & \qquad-y_{3}\tau E(\_\otimes ee_{3})-E\delta(\_\otimes ee_{3})\\
 & =\_\otimes ee_{1}+y_{3}\bigl(\tau E(\_\otimes ee_{2})+y_{1}y_{2}\chi''\bigr)\\
 & \qquad-y_{3}\tau E(\_\otimes ee_{3})-E\delta(\_\otimes ee_{3})\\
 & =\_\otimes\bigl(ee_{1}-\delta(ee_{3})\bigr)+y_{1}\Bigl(-y_{3}\tau E(\_\otimes ee'')+y_{2}y_{3}\chi''\Bigr)\\
 & =y_{1}\Bigl(-\_\otimes ee'''-y_{3}\tau E(\_\otimes ee'')+y_{2}y_{3}\chi''\Bigr).
\end{align*}
 The final statement of the proposition is a rearrangement of the
second and third equalities of the second condition.
\end{proof}
\begin{rem}
We will not need to use alternative conditions for $G_{n}$ for $n\geq4$.
\end{rem}

\subsection{Product category $C\textrm{-mod}$}

Let $C=\End_{\text{per }B}(X)^{\op}$. We \textquoteleft change basis\textquoteright{}
from $Be_{1}\oplus Be_{2}$ to $X_{1}\oplus X_{2}$, i.e.~from complexes
of modules over $B$ to complexes of modules over $C$. This is performed
by $\mathscr{H}om_{B}(X,-)$: 
\[
\text{per }B\xrightarrow[\mathscr{H}om_{B}(X,-)]{\sim}\text{per }C,
\]
which is a restricted Rickard (derived Morita) equivalence. It has
an inverse given by $X\otimes_{C}-$. Under this equivalence, the
action of $_{B}E'\otimes_{B}-$ on $\text{per }B$ translates to $_{C}\tilde{E}\otimes_{C}-$
on $\text{per }C$, where $\tilde{E}$ is a $(C,C)$-bimodule that
is finitely generated and projective on the left. Our main theorem
says that $\mathsf{Bim}_{k}(C)$ has the structure of 2-representation
of $\mathcal{U}^{+}$ using $\tilde{E}$. In this section we describe
$C$ and the derived equivalence in more detail.

\subsubsection{New algebra $C$}

Let $\mathscr{C}=\mathscr{E}nd_{B}(X_{1}\oplus X_{2})^{\op}$ be the
$\dg$-algebra of endomorphisms of $X$ (with left-to-right composition).
\begin{defn}
Define two $(A[y],A[y])$-bimodules: 
\[
G_{1}'=A[y]\oplus\Hom_{A[y]}(_{A[y]}E[y],E[y])
\]
 and 
\[
G_{1}''=\Hom_{A[y]}(_{A[y]}E[y],E_{y}).
\]
\end{defn}

The complex $\mathscr{E}nd_{B}(X_{2})$ is given in degrees $0$ and
$1$ by 
\[
G_{1}'\overset{d^{0}}{\to}G_{1}''
\]
 where 
\[
d^{0}\bigl((\theta(y),\varphi)\bigr)=\pi\circ\varphi-\pi(-).\theta(x).
\]

The direct sum decomposition $X_{1}\oplus X_{2}$ provides a matrix
presentation for $\mathscr{C}$ with $\mathscr{C}_{ij}=\mathscr{H}om_{B}(X_{i},X_{j})$.
\begin{defn}
Let $F$ denote the $(A,A)$-bimodule 
\[
F=\Hom_{A}(_{A}E,A).
\]
\end{defn}

Note the canonical isomorphism 
\[
\Hom_{A}(_{A}E,A)[y]\iso\Hom_{A[y]}(_{A[y]}E[y],A[y])
\]
 that exists because $_{A}E$ is finitely generated. Since $_{A}E$
and $_{A[y]}E[y]$ are both finitely generated projective, we also
have canonical isomorphisms of functors: 
\begin{align*}
\Hom_{A}(_{A}E,-) & \iso\Hom_{A}(_{A}E,A)\otimes_{A}-\\
\Hom_{A[y]}(_{A[y]}E[y],-) & \iso\Hom_{A[y]}(_{A[y]}E[y],A[y])\otimes_{A[y]}-.
\end{align*}

\begin{prop}
\label{prop:desc-of-scrC} The algebra $\mathscr{C}$ is isomorphic
to a generalized matrix algebra of complexes concentrated in degrees
$0$ and $1$: 
\[
\begin{pmatrix}A[y] & E[y]\xrightarrow{\pi}E_{y}\\
F[y] & G_{1}'^{\op}\xrightarrow{d^{0}}G_{1}''^{\op}
\end{pmatrix}\iso\begin{pmatrix}\mathscr{C}_{11} & \mathscr{C}_{12}\\
\mathscr{C}_{21} & \mathscr{C}_{22}
\end{pmatrix}.
\]

The map is given on components by: 
\begin{itemize}
\item for $\mathscr{C}_{11}$:
\[
A[y]\ni a\mapsto\left(\begin{pmatrix}1\\
0
\end{pmatrix}\mapsto\begin{pmatrix}a\\
0
\end{pmatrix}\right)
\]
\item for $\mathscr{C}_{12}$: 
\[
\left(E[y]\to E_{y}\right)\ni(e,e')\mapsto\left(\begin{pmatrix}1\\
0
\end{pmatrix}\mapsto\begin{pmatrix}(e,e')\\
0
\end{pmatrix}\right)
\]
\item for $\mathscr{C}_{21}$: 
\[
F[y]\ni f\mapsto\left(\begin{pmatrix}(e,0)\\
1
\end{pmatrix}\mapsto\begin{pmatrix}f(e)\\
0
\end{pmatrix}\right)
\]
\item for $\mathscr{C}_{22}$: 
\begin{align*}
\left(G_{1}'^{\op}\to G_{1}''^{\op}\right)\ni\bigl((\theta,\varphi'),\varphi''\bigr) & \mapsto\left(\begin{pmatrix}(e,0)\\
1
\end{pmatrix}\mapsto\begin{pmatrix}\bigl(\varphi'(e),(\pi\circ\varphi'')(e)\bigr)\\
\theta
\end{pmatrix}\right).
\end{align*}
\end{itemize}
\end{prop}

\begin{proof}
Computation.
\end{proof}
\begin{defn}
\label{def:C} Let $C$ denote the $k$-algebra $\End_{K^{b}(B)}(X)^{\op}$.
\end{defn}

Sometimes we consider $C$ to be a $\dg$-algebra concentrated in
degree $0$.
\begin{lem}
The projection $Z^{0}(\mathscr{C})\to H^{0}(\mathscr{C})=C$ is an
isomorphism. Its inverse gives an injection $C\hookrightarrow\mathscr{C}$
which is a quasi-isomorphism of $\dg$-algebras.
\end{lem}

\begin{proof}
The first claim follows because $\mathscr{C}$ lies in degrees $0$
and $1$. For the second claim we just need that $H^{1}(\mathscr{C})=0$.
It is clear that the map $\pi:E[y]\to E_{y}$ is surjective. We can
see that $d^{0}$ is surjective as well: since $_{A[y]}E[y]$ is projective,
$\Hom_{A[y]}(_{A[y]}E[y],-)$ is exact, so 
\[
\Hom_{A[y]}(_{A[y]}E[y],\pi):\Hom_{A[y]}(_{A[y]}E[y],E[y])\to\Hom_{A[y]}(_{A[y]}E[y],E_{y})
\]
 is surjective.
\end{proof}
The injection of the lemma gives a right action of $C$ on $X$.
\begin{lem}
The algebra $C$ is isomorphic to a generalized matrix algebra: 
\[
\begin{pmatrix}A[y] & y_{1}E[y]\\
F[y] & G_{1}^{\op}
\end{pmatrix}\iso\begin{pmatrix}C_{11} & C_{12}\\
C_{21} & C_{22}
\end{pmatrix},
\]
 with component maps given by (restrictions of) those in Proposition
\ref{prop:desc-of-scrC}.
\end{lem}

\begin{proof}
We have $d^{0}\bigl((\theta,\varphi)\bigr)=0$ exactly when $\varphi=\_.\theta+y_{1}\varphi'$
for some $\varphi'\in\Hom_{A[y]}(_{A[y]}E[y],E[y])$, and it follows
that the map to $C_{22}$ is an isomorphism.
\end{proof}

\subsubsection{Derived equivalence}

Since $X$ is strictly perfect, the triangulated functor 
\[
\mathscr{H}om_{B}(X,-):K^{b}(B)\to K^{b}(C)
\]
 descends to the derived categories and resolutions are not needed:
\[
\mathscr{H}om_{B}(X,-):D^{b}(B)\to D^{b}(C).
\]
 Since $X$ generates $\text{per }B$, it is perfect as a right $\mathscr{C}$-$\dg$-module,
and then also as a complex of $C$-modules because the inclusion $C\overset{}{\hookrightarrow}\mathscr{C}$
is a quasi-isomorphism. It follows that the functor restricts to a
functor 
\[
\mathscr{H}om_{B}(X,-):\text{per }B\to\text{per }C,
\]
 and this is essentially surjective because $C$ is in the essential
image. To show that the functor is fully faithful, it is enough to
check endomorphisms of $X$ and its translates, since $X$ generates
$\text{per }B$. The induced map: 
\[
\Hom_{D^{b}(B)}(X,X[i])\to\Hom_{D^{b}(C)}(\mathscr{E}nd_{B}(X),\mathscr{E}nd_{B}(X)[i])
\]
 is an isomorphism for all $i$: with $i=0$ both sides are canonically
isomorphic to $C$, and the map induces the identity on $C$; with
$i\neq0$ both sides are $0$.

The endofunctor $E'\otimes_{B}-$ on $\text{per }B$ induces an endofunctor
on $\text{per }C$ using this equivalence: first apply $X\otimes_{C}-$,
then $E'\otimes_{B}-$, then $\mathscr{H}om_{B}(X,-)$. Since $X$
is finitely generated and strictly perfect, this induced endofunctor
is isomorphic to $\mathscr{H}om_{B}(X,E'X)\otimes_{C}-$.
\begin{rem}
In the above context a theorem of Rickard shows that $\mathscr{H}om_{B}(X,-):D^{b}(B)\to D^{b}(C)$
is also an equivalence of categories. We do not know $E'\otimes_{B}-$
to be exact, however, so we use the restricted equivalence of perfect
complexes, and the full version of Rickard\textquoteright s theorem
is not needed.
\end{rem}

\begin{defn}
In §\ref{sec:Product-category}, let $\mathscr{E}$ denote the $(C,C)$-bimodule
complex $\mathscr{H}om_{B}(X,E'X)$.
\end{defn}

Then we have the following:
\begin{prop}
\label{prop:^n-scrHom-E'-to-scrHom-E'^n} For each $n$, the morphism
of $(C,C)$ bimodule complexes 
\[
\overbrace{\mathscr{E}\otimes_{C}\dots\otimes_{C}\mathscr{E}}^{n\text{-times}}\to\mathscr{H}om_{B}(X,E'^{n}X)
\]
 given by 
\[
f_{1}\otimes\dots\otimes f_{n}\mapsto E'^{n-1}(f_{n})\circ E'^{n-2}(f_{n-1})\circ\dots\circ f_{1}
\]
 is a quasi-isomorphism. These maps give the vertical maps in diagrams
of the following form, which commute: 
\end{prop}

\noindent \scalebox{0.9}{%
\noindent\begin{minipage}[t][0.7\totalheight]{1\columnwidth}%
\[\begin{tikzcd}[ampersand replacement=\&] 	{\mathscr{H}om_B(X,E'X)^{\otimes\, n}\otimes_C\mathscr{H}om_B(X,E'X)^{\otimes\, m}} \& {} \& {\mathscr{H}om_B(X,E'X)^{\otimes\, n+m}} \\ 	{\mathscr{H}om_B(X,E'^{n}X)\otimes_C\mathscr{H}om_B(X,E'^{m}X)} \& {} \& {\mathscr{H}om_B(X,E'^{n+m}X).} \\ 	\\ 	{} \& {} \& {} 	\arrow[from=1-1, to=2-1] 	\arrow["{f\otimes g\mapsto E'^n(g)\circ f}", from=2-1, to=2-3] 	\arrow[from=1-1, to=1-3] 	\arrow[from=1-3, to=2-3] \end{tikzcd}\]%
\end{minipage}}
\begin{proof}
All diagrams contained in the following diagram commute, up to canonical
isomorphisms in $\text{per }B$ and $\text{per }C$: 

\noindent %
\noindent\begin{minipage}[t]{1\columnwidth}%
\[\begin{tikzcd}[ampersand replacement=\&] 	{\text{per }B} \&\& {\text{per }C} \\ 	{\text{per }B} \&\& {\text{per }C} \\ 	{\text{per }B} \&\& {\text{per }C.} 	\arrow["{\mathscr{H}om_B(X,-)}", shift left=1, from=1-1, to=1-3] 	\arrow["{\mathscr{H}om_B(X,-)}", shift left=1, from=2-1, to=2-3] 	\arrow["{\mathscr{H}om_B(X,-)}", shift left=1, from=3-1, to=3-3] 	\arrow["{E'\otimes_B -}"', from=1-1, to=2-1] 	\arrow["{E'\otimes_B -}"', from=2-1, to=3-1] 	\arrow["{\mathscr{E}\otimes_C-}", from=1-3, to=2-3] 	\arrow["{\mathscr{E}\otimes_C-}", from=2-3, to=3-3] 	\arrow["{X\otimes_C -}", shift left=1, from=1-3, to=1-1] 	\arrow["{X\otimes_C -}", shift left=1, from=2-3, to=2-1] 	\arrow["{X\otimes_C -}", shift left=1, from=3-3, to=3-1] \end{tikzcd}\]%
\end{minipage}

This gives the first statement of the proposition. The diagrams commute
by functoriality of $E'$.
\end{proof}

\subsection{New bimodule $\tilde{E}$}

\subsubsection{Definition of $\tilde{E}$}

Now we define the lead actor of this paper.
\begin{defn}
\label{def:tilde-E} Define a $(C,C)$-bimodule:
\[
\tilde{E}=\Hom_{K^{b}(B)}(X,E'X),
\]
with left $C$ action given by precomposition with $\varphi\in C$,
and right $C$ action given by post-composition with $E'(\varphi)$
for $\varphi\in C$.
\end{defn}

\begin{lem}
\label{lem:scrHom-E'^n-in-nonneg-degree} For each $n$, the complex
$\mathscr{H}om_{B}(X,E'^{n}X)$ of $(C,C)$-bimodules is concentrated
in nonnegative degree.
\end{lem}

\begin{proof}
The lower row of $E'^{n}X$ has components in degrees at least $1$,
and the upper row has components in degrees at least $0$. This is
shown by a simple inductive argument using the formulas for $X$ and
$E'$ in §\ref{subsec:endofunctor-E'}. It follows that there are
no nonzero morphisms in $\mathscr{H}om_{B}(X,E'^{n}X)$ of negative
degree.
\end{proof}
\begin{prop}
The complex $\mathscr{E}=\mathscr{H}om_{B}(X,E'X)$ of $(C,C)$-bimodules
has cohomology concentrated in degree $0$.
\end{prop}

\begin{proof}
We consider separately the matrix components $\mathscr{H}om_{B}(X_{i},E'X_{j}):$ 
\begin{itemize}
\item $\mathscr{H}om_{B}(X_{1},E'X_{1})$: since $X_{1}=Be_{1}$ this is
isomorphic to $e_{1}E'X_{1}$ which is $E[y]\xrightarrow{\pi}E_{y}$,
and $\pi$ is surjective.
\item $\mathscr{H}om_{B}(X_{1},E'X_{2})$: this is isomorphic to $e_{1}E'^{2}X_{1}$,
which is 
\[
E^{2}[y]\xrightarrow{\left(\begin{smallmatrix}\pi_{2}\\
\pi_{1}
\end{smallmatrix}\right)}E_{y}E\oplus EE_{y}\xrightarrow{\left(\begin{smallmatrix}-\pi_{1} & \pi_{2}\end{smallmatrix}\right)}E_{y}E_{y}.
\]
 The second map is clearly surjective. Its kernel consists of pairs
$(ee_{1},ee_{2})\in E^{2}$ such that $ee_{1}-ee_{2}=(Ex-xE)ee^{\circ}$
for some $ee^{\circ}\in E^{2}$. Such a pair is the image of $ee_{2}+(Ex-y)ee^{\circ}$
in $E^{2}[y]$.
\item $\mathscr{H}om_{B}(X_{2},E'X_{1})$: this is isomorphic to $\mathscr{C}_{22}$,
and we saw that $d^{0}$ is surjective.
\item $\mathscr{H}om_{B}(X_{2},E'X_{2})$: this is isomorphic to $G_{2}'\xrightarrow{d^{0}}G_{2}''\xrightarrow{d^{1}}G_{2}'''$,
where 
\begin{align*}
G_{2}' & =E[y]^{\oplus2}\oplus\Hom_{A[y]}(_{A[y]}E[y],E^{2}[y])\\
G_{2}'' & =E_{y}\oplus\Hom_{A[y]}(_{A[y]}E[y],E_{y}E\oplus EE_{y})\\
G_{2}''' & =\Hom_{A[y]}(_{A[y]}E[y],E_{y}E_{y}),
\end{align*}
 with 
\begin{align*}
d^{0} & :(e_{1},e_{2},\xi)\mapsto\bigl(\pi(e_{2}-e_{1}),(\pi_{2}\circ\xi;\pi_{1}\circ\xi)\bigr)\\
d^{1} & :\bigl(e,(\xi';\xi'')\bigr)\mapsto-\pi_{1}\circ\xi'+\pi_{2}\circ\xi''.
\end{align*}
 It is easy to see that $H^{1}=0$ and $H^{2}=0$ by applying the
exact functor $\Hom_{A[y]}(_{A[y]}E[y],-)$ to the sequence considered
in the second bullet.
\end{itemize}
\end{proof}
\begin{cor}
\label{cor:tilde-E-to-cplx-E} The surjection 
\[
Z^{0}\mathscr{H}om_{B}(X,E'X)\to H^{0}\mathscr{H}om_{B}(X,E'X)=\tilde{E}
\]
is an isomorphism. Its inverse gives an injection 
\[
\tilde{E}\hookrightarrow\mathscr{E}
\]
 which is a quasi-isomorphism of complexes of $(C,C)$-bimodules. 
\end{cor}

\begin{rem}
Whereas $E'$ is a complex of bimodules, $\tilde{E}$ is just a bimodule.
This observation is the starting point for our construction. The basis
$X_{1}\oplus X_{2}$ is designed to be more compatible with the $\mathcal{U}^{+}$
action in this sense.
\end{rem}

\begin{lem}
As a left $C$-module, $\tilde{E}$ is finitely generated and projective.
\end{lem}

\begin{proof}
In Prop.~\ref{prop:E'(X)-perfect} we saw that $E'X$ is quasi-isomorphic
to a finite direct sum of summands of $X$, so $_{C}\tilde{E}$ is
a finite direct sum of summands of $C$.
\end{proof}
\begin{lem}
\label{lem:E-tilde^n-to-Hom-cplx} The map $\tilde{E}^{n}\to\mathscr{H}om_{B}(X,E'^{n}X)$
of complexes of $(C,C)$-bimodules given by 
\[
f_{1}\otimes\dots\otimes f_{n}\mapsto E'^{n-1}(f_{n})\circ E'^{n-2}(f_{n-1})\circ\dots\circ f_{1}
\]
 is a quasi-isomorphism.
\end{lem}

\begin{proof}
Use a copy of the morphism 
\[
\tilde{E}\overset{q.i.}{\hookrightarrow}\mathscr{E}
\]
 from Corollary \ref{cor:tilde-E-to-cplx-E} onto each factor of the
product on the left in Proposition \ref{prop:^n-scrHom-E'-to-scrHom-E'^n},
and the fact that $\tilde{E}$ is finitely generated and projective
on the left.
\end{proof}
\begin{lem}
\label{lem:Nat-Isos_E'^n}The maps of Lemma \ref{lem:E-tilde^n-to-Hom-cplx}
induce isomorphisms of $(C,C)$-bimodules 
\[
\tilde{E}^{n}\iso\Hom_{K^{b}(B)}(X,E'^{n}X)
\]
 making the following diagrams commute: 

\noindent \scalebox{0.95}{%
\noindent\begin{minipage}[t][0.7\totalheight]{1\columnwidth}%
\[\begin{tikzcd}[ampersand replacement=\&] 	{\tilde{E}^{n}\otimes_C \tilde{E}^{m}} \& {\tilde{E}^{n+m}} \\ 	{\mathrm{Hom}_{K(B)}(X,E'^{n}(X))\otimes_C\mathrm{Hom}_{K(B)}(X,E'^{m}(X))} \& {\mathrm{Hom}_{K(B)}(X,E'^{n+m}(X)).} \& {} \\ 	\\ 	{} \& {} \& {} 	\arrow["\sim"', from=1-1, to=2-1] 	\arrow["\sim"', from=1-2, to=2-2] 	\arrow["\sim", from=1-1, to=1-2] 	\arrow["\sim", from=2-1, to=2-2] \end{tikzcd}\]%
\end{minipage}}
\end{lem}

\begin{proof}
By Lemma \ref{lem:E-tilde^n-to-Hom-cplx}, the cohomology of $\mathscr{H}om_{B}(X,E'^{n}X)$
is concentrated in degree $0$. By Lemma \ref{lem:scrHom-E'^n-in-nonneg-degree},
\[
Z^{0}\mathscr{H}om_{B}(X,E'^{n}X)=H^{0}\mathscr{H}om_{B}(X,E'^{n}X).
\]
 So the degree $0$ part of the map of Lemma \ref{lem:E-tilde^n-to-Hom-cplx}
is an isomorphism from $\tilde{E}^{n}$ to $Z^{0}\mathscr{H}om_{B}(X,E'^{n}X)$,
which is $\Hom_{K^{b}(B)}(X,E'^{n}X)$. The diagrams commute because
the morphisms are restrictions of the morphisms of Proposition \ref{prop:^n-scrHom-E'-to-scrHom-E'^n}.
\end{proof}
\begin{defn}
We let $\tilde{E}_{ij}^{n}$ denote $\Hom_{K^{b}(B)}(X_{i},E'^{n}X_{j})$.
\end{defn}

Defined in this way, $\tilde{E}_{ij}^{n}$ lies in $\Hom_{K^{b}(B)}(X,E'^{n}X)$,
not in $\tilde{E}^{n}$, but we consider it also in the latter through
the isomorphism of Lemma \ref{lem:Nat-Isos_E'^n}.

\subsubsection{\label{subsec:Desc-of-E^2}\label{subsec:Desc-of-E^3} Some low powers
of $\protect\ee$}

The bimodule $\tilde{E}$ can be presented as a matrix with $ij$-component
$\tilde{E}_{ij}$ given by $\Hom_{K^{b}(B)}(X_{i},E'X_{j})$. This
component is an $(\End(X_{i})^{\op},\End(X_{j})^{\op})$-bimodule.
Recall that $\End(X_{1})^{\op}\cong A[y]$ and $\End(X_{2})^{\op}\cong G_{1}^{\op}$.
\begin{lem}
\label{lem:E'n(X_1)} We have 
\[
y_{1}\dots y_{n}E^{n}[y]\iso\Hom_{K^{b}(B)}(X_{1},E'^{n}X_{1}),
\]
 where $y_{1}\dots y_{n}e$ is sent to the map in $K^{b}(B)$ determined
by: 
\[
\begin{pmatrix}1\\
0
\end{pmatrix}\mapsto\begin{pmatrix}(y_{1}\dots y_{n}e,0,\dots,0)\\
0
\end{pmatrix}.
\]
\end{lem}

\begin{proof}
Computation. Note that $E'^{n}X_{1}$ has just one term in degree
$0$, which is $E^{n}[y]$ in the upper row. The differential of $E'^{n}X_{1}$
out of this term is the map whose kernel is computed in Lemma \ref{lem:kernels-lemma}.
\end{proof}
\begin{prop}
We have: 
\[
\begin{pmatrix}y_{1}\dots y_{n}E^{n}[y] & y_{1}\dots y_{n+1}E^{n+1}[y]\\
G_{n} & G_{n+1}
\end{pmatrix}\iso\begin{pmatrix}\tilde{E}_{11}^{n} & \tilde{E}_{12}^{n}\\
\tilde{E}_{21}^{n} & \tilde{E}_{22}^{n}
\end{pmatrix},
\]
where the maps on the upper row are from Lemma \ref{lem:E'n(X_1)},
and on the lower they are from the definition of $G_{n}$.
\end{prop}

Together with Lemma \ref{lem:Nat-Isos_E'^n}, this gives a parametrization
of $\tilde{E}^{n}$. We may record the matrix presentations for the
first few powers:

\[
\begin{pmatrix}y_{1}E[y] & y_{1}y_{2}E^{2}[y]\\
G_{1} & G_{2}
\end{pmatrix}\iso\tilde{E},
\]

\[
\begin{pmatrix}y_{1}y_{2}E^{2}[y] & y_{1}y_{2}y_{3}E^{3}[y]\\
G_{2} & G_{3}
\end{pmatrix}\iso\tilde{E}^{2},
\]

\[
\begin{pmatrix}y_{1}y_{2}y_{3}E^{3}[y] & y_{1}y_{2}y_{3}y_{4}E^{4}[y]\\
G_{3} & G_{4}
\end{pmatrix}\iso\ee^{3}.
\]

\section{Hecke action}

In this section we introduce $(C,C)$-bimodule endomorphisms $\tilde{x}$
of $\tilde{E}$ and $\tilde{\tau}$ of $\tilde{E}^{2}$, and show
that they satisfy the relations of $\mathcal{U}^{+}$.

\subsection{Definition of the action}

In §\ref{subsec:Defs_x_tau} we give formulas for endomorphisms of
the separate components of $\tilde{E}$ and $\tilde{E}^{2}$. A few
lemmas are needed first in order to show that the formulas are well-defined
on components of the form $G_{n}$, $n=1,2,3$. Then in §\ref{subsec:x,tau-C-equiv}
we argue that these componentwise definitions jointly determine a
morphism of $(C,C)$-bimodules.

\subsubsection{\label{subsec:Defs_x_tau} Formulas for $\protect\xx$ and $\protect\tt$}
\begin{lem}
Let $(\theta,\varphi)\in G_{1}\subset A^{\op}[y]\oplus\Hom_{A}(_{A}E,E)[y]$.
Then $(y\theta,x\circ\varphi)\in G_{1}$.
\end{lem}

\begin{proof}
Compute: 
\begin{align*}
x\circ\varphi-y\theta & =x(\_.\theta+y_{1}\varphi_{1})-y\theta\\
 & =y_{1}(\_.\theta+x\varphi_{1}).
\end{align*}
\end{proof}
\begin{lem}
Let $(e_{1},e_{2},\xi)\in G_{2}\subset E[y]^{\oplus2}\oplus\Hom_{A}(_{A}E,E^{2})[y]$.
Then $(ye_{1},xe_{2},xE\circ\xi)\in G_{2}$ and $(e',e',\tau\circ\xi)\in G_{2}$,
where $e'$ is uniquely determined by $e_{1}-e_{2}=y_{1}e'$.
\end{lem}

\begin{proof}
For the first claim, compute: 
\begin{align*}
xE\circ\xi-\_\otimes ye_{1} & =xE\circ(\_\otimes e_{1}+y_{2}\xi_{1})-\_\otimes ye_{1}\\
 & =y_{2}(\_\otimes e_{1}+xE\circ\xi_{1}),
\end{align*}
 and 
\begin{align*}
xE\circ\xi-\delta(\_\otimes xe_{2}) & =xE\circ\bigl(\delta(\_\otimes e_{2})+y_{1}\xi_{2}\bigr)-\delta(\_\otimes xe_{2})\\
 & =\delta\circ Ex(\_\otimes e_{2})-y_{1}(\_\otimes e_{2})\\
 & \qquad+y_{1}xE\circ\xi_{2}-\delta(\_\otimes xe_{2})\\
 & =y_{1}(-\_\otimes e_{2}+xE\circ\xi_{2}).
\end{align*}
 For the second claim, use the alternative characterization of $G_{2}$
as given in Prop.~\ref{prop:K-conds}, and compute: 
\begin{align*}
\tau\circ\xi & =\tau(\_\otimes e_{1})+\tau y_{2}\xi_{1}\\
 & =\tau(\_\otimes e_{1})+y_{1}\tau\xi_{1}-\xi_{1}\\
 & =\tau(\_\otimes e_{1})+y_{1}\tau y_{1}\xi'-\xi_{1}\\
 & =\tau\bigl(\_\otimes(e_{1}-e_{2})\bigr)+y_{1}y_{2}\tau\xi'\\
 & =\tau y_{1}(\_\otimes e')+y_{1}y_{2}\tau\xi'\\
 & =\_\otimes e'+y_{2}\bigl(\tau(\_\otimes e')+y_{1}\tau\xi'\bigr).
\end{align*}
 The last line has the form of an element of $G_{2}$.
\end{proof}
\begin{lem}
\label{lem:tau-on-chi} Let $(ee_{1},ee_{2},ee_{3},\chi)\in G_{3}\subset E^{2}[y]^{\oplus3}\oplus\Hom_{A}(_{A}E,E^{3})[y]$.
Then $(ee',ee',\tau(ee_{3}),\tau E\circ\chi)\in G_{3}$, where $ee'$
is uniquely determined by $ee_{1}-ee_{2}=y_{2}ee'$.
\end{lem}

\begin{proof}
We use the alternative characterization of $G_{3}$ as given in Prop.~\ref{prop:L-conds},
and compute: 
\begin{align*}
\tau E\circ\chi & =\tau E(\_\otimes ee_{1})+\tau Ey_{3}\chi_{1}\\
 & =\tau E(\_\otimes ee_{1})-\chi_{1}+y_{2}\tau E\circ\chi_{1}\\
 & =\tau E(\_\otimes ee_{1})-\chi_{1}+y_{2}\tau Ey_{2}\bigl(E\tau\circ\tau E(\_\otimes ee_{3})+y_{1}\chi''\bigr)\\
 & =\tau E(\_\otimes ee_{1})-\chi_{1}\\
 & \qquad+\bigl(y_{2}y_{3}\tau E+y_{2}\bigr)\cdot\bigl(E\tau\circ\tau E(\_\otimes ee_{3})+y_{1}\chi''\bigr)\\
 & =\tau E(\_\otimes(ee_{1}-ee_{2}))\\
 & \qquad+y_{2}y_{3}\bigl(\tau E\circ E\tau\circ\tau E(\_\otimes ee_{3})+y_{1}\tau E\circ\chi''\bigr)\\
 & =\tau Ey_{2}(\_\otimes ee')\\
 & \qquad+y_{2}y_{3}\Bigl(E\tau\circ\tau E(\_\otimes ee')+E\delta\circ\tau E(\_\otimes\overline{ee})+y_{1}\tau E\circ\chi''\Bigr)\\
 & =\_\otimes ee'+y_{3}\cdot\\
 & \qquad\Bigl(E\delta\circ\tau E(\_\otimes ee')+y_{2}\bigl(E\delta\circ\tau E(\_\otimes\overline{ee})+y_{1}y_{2}\tau E\circ\chi''\bigr)\Bigr)\\
 & =\_\otimes ee'+y_{3}\cdot\\
 & \qquad\Bigl(-\tau Ey_{1}(\_\otimes\overline{ee})+E\delta\circ\tau E(\_\otimes\tau(ee_{3}))+y_{1}y_{2}\tau E\circ\chi''\Bigr).
\end{align*}
The last line has the form of an element of $G_{3}$, namely $(ee',ee',\tau(ee_{3}),\tau E\circ\chi)$.
\end{proof}
The element $(ee_{1},ee_{2},ee_{3},\chi)\in G_{3}$ is associated
(by Prop.~\ref{prop:L-conds}) with further data that has been notated
$ee^{\ell}$, $\overline{ee}$, $\chi_{\ell}$, $\chi_{1}'$, and
$\chi''$. We record the corresponding data associated with $(ee',ee',\tau(ee_{3}),\tau E\circ\chi)$
using the notation $\bar{e}$ and $\bar{\chi}$ for the new versions:
\begin{align*}
\bar{e}\bar{e}' & =0\\
\bar{e}\bar{e}'' & =\overline{ee}\\
\bar{e}\bar{e}''' & =\overline{ee}\\
\overline{\bar{e}\bar{e}} & =0,
\end{align*}
 and 
\begin{align*}
\bar{\chi} & =(ee',ee',\tau(ee_{3}),\tau E\circ\chi)\\
\bar{\chi}_{1} & =-\tau Ey_{1}(\_\otimes\overline{ee})+E\delta\circ\tau E\circ E\tau(\_\otimes ee_{3})+y_{1}y_{2}\tau E\circ\chi''\\
\bar{\chi}_{2} & =E\tau\circ\delta E\circ E\tau(\_\otimes ee_{3})+y_{1}y_{3}\tau E\circ\chi''\\
\bar{\chi}_{3} & =-\delta E(\_\otimes\overline{ee})+y_{2}y_{3}\tau E\circ\chi''\\
\bar{\chi}_{1}' & =E\tau\circ\tau E\circ E\tau(\_\otimes ee_{3})+y_{1}\tau E\circ\chi''\\
\bar{\chi}'' & =\tau E\circ\chi''.
\end{align*}

Now we give componentwise formulas for $\tilde{x}$ and $\tilde{\tau}$.
These formulas are well-defined on $\tilde{E}_{21}$, $\tilde{E}_{22}$,
$\tilde{E}_{21}^{2}$, and $\tilde{E}_{22}^{2}$ by the lemmas above.
\begin{defn}
\noindent \label{def:tilde-x-tilde-tau} We define the action of $\xx$
on $\ee$ as follows:
\begin{itemize}
\item on $\ee_{11}$: $\xx$ acts by $x$
\item on $\ee_{12}$: $\xx$ acts by $xE$
\item on $\ee_{21}$: $\xx$ acts by $(\theta,\varphi)\mapsto(y\theta,x\circ\varphi)$
\item on $\ee_{22}$: $\xx$ acts by $(e_{1},e_{2},\xi)\mapsto(ye_{1},xe_{2},xE\circ\xi)$.
\end{itemize}
We define the action of $\tt$ on $\ee^{2}$ as follows:
\begin{itemize}
\item on $\ee_{11}^{2}$: $\tt$ acts by $\tau$
\item on $\ee_{12}^{2}$: $\tt$ acts by $\tau E$
\item on $\ee_{21}^{2}$: $\tt$ acts by $(e_{1},e_{2},\xi)\mapsto(e',e',\tau\circ\xi)$
\item on $\ee_{22}^{2}$: $\tt$ acts by ${(ee_{1},ee_{2},ee_{3},\chi)\mapsto(ee',ee',\tau(ee_{3}),\tau E\circ\chi)}$,
\end{itemize}
where $e'$ and $ee'$ are uniquely determined by $e_{1}-e_{2}=y_{1}e'$
and $ee_{1}-ee_{2}=y_{2}ee'$.

\end{defn}

\begin{lem}
The formulas for $\tilde{x}$ give a $(C,C)$-bimodule endomorphism
of $\tilde{E}$.
\end{lem}

\begin{proof}
Recall the definition of the complex $E'$ of $(B,B)$-bimodules in
§\ref{subsec:endofunctor-E'}. There is an $(\left(\begin{smallmatrix}A[y] & 0\\
0 & A[y]
\end{smallmatrix}\right),\left(\begin{smallmatrix}A[y] & 0\\
0 & A[y]
\end{smallmatrix}\right))$-bimodule endomorphism $x'$ of $E'$ given componentwise in degrees
$0$ and $1$ by $(A[y],A[y])$-bimodule endomorphisms: 
\[
x'_{0}=\begin{pmatrix}x & xE_{y}\\
0 & x
\end{pmatrix},\quad x'_{1}=\begin{pmatrix}x & xE_{y}\\
y & x
\end{pmatrix}.
\]
The relation $s\circ E_{y}x=xE_{y}\circ s$ may be used to check that
$x_{0}'$ and $x_{1}'$ together give a morphism of complexes of $(B,B)$-bimodules.
This map induces a $(C,C)$-bimodule endomorphism of $\Hom_{K^{b}(B)}(X,E'X)$
that agrees with the definition of $\tilde{x}$.
\end{proof}
It follows that $\tilde{x}$ induces endomorphisms $\xx\ee$ and $\ee\xx$.
For future reference we write the formulas for those: 
\begin{prop}
The formulas for $\xx$ determine the following formulas for $\xx\ee$
and $\ee\xx$ on $\ee^{2}$: 
\end{prop}

\begin{itemize}
\item on $\ee_{11}^{2}$: $\xx\ee$ acts by $xE$ and $\ee\xx$ acts by
$Ex$
\item on $\ee_{12}^{2}$: $\xx\ee$ acts by $xE^{2}$ and $\ee\xx$ by $ExE$
\item on $\ee_{21}^{2}$: $\xx\ee$ acts by 
\[
(e_{1},e_{2},\xi)\mapsto(ye_{1},xe_{2},xE\circ\xi)
\]
 and $\ee\xx$ by 
\[
{(e_{1},e_{2},\xi)\mapsto(xe_{1},ye_{2},Ex\circ\xi)}
\]
\item on $\ee_{22}^{2}$: $\xx\ee$ acts by 
\[
{(ee_{1},ee_{2},ee_{3},\chi)\mapsto(yee_{1},xE(ee_{2}),xE(ee_{3}),xE^{2}\circ\chi)}
\]
 and $\ee\xx$ by 
\[
{(ee_{1},ee_{2},ee_{3},\chi)\mapsto(xE(ee_{1}),yee_{2},Ex(ee_{3}),ExE\circ\chi)}.
\]
\end{itemize}
\begin{proof}
Use Lemma \ref{lem:Nat-Isos_E'^n}, in particular the diagram in the
case $n=m=1$.
\end{proof}

\subsubsection{\label{subsec:x,tau-C-equiv} Bimodule structure of $\tilde{E}^{2}$
and equivariance of $\tilde{\tau}$}
\begin{lem}
\label{lem:tilde-tau-is-endo} The formulas for $\tilde{\tau}$ give
a $(C,C)$-bimodule endomorphism of $\tilde{E}^{2}$.
\end{lem}

For the maps we defined on components of $\tilde{E}^{2}$ to determine
jointly a $(C,C)$-bimodule endomorphism $\tilde{\tau}$, they must
be equivariant with respect to the left and right $C$-actions. In
order to check equivariance, we write formulas for the actions of
the generators in $C$ in the following four lemmas. The reader may
verify these formulas from the various definitions.
\begin{lem}
Generators in $A[y]\subset C$ act on the right on $\tilde{E}^{2}$,
in terms of the separate bimodule structures of $\tilde{E}_{ij}^{2}$,
as follows: 
\begin{itemize}
\item $\tilde{E}_{11}^{2}\otimes A[y]\to\tilde{E}_{11}^{2}$ by 
\begin{align*}
y_{1}y_{2}E^{2}[y]\otimes_{A[y]}A[y] & \lra y_{1}y_{2}E^{2}[y]\\
y_{1}y_{2}ee\otimes\theta & \mapsto y_{1}y_{2}ee.\theta.
\end{align*}
\item $\tilde{E}_{21}^{2}\otimes A[y]\to\tilde{E}_{21}^{2}$ by 
\begin{align*}
G_{2}\otimes_{A[y]}A[y] & \lra G_{2}\\
\left(e_{1},e_{2},\xi\right)\otimes\theta & \mapsto\left(e_{1}.\theta,e_{2}.\theta,\xi(-).\theta\right).
\end{align*}
\end{itemize}
They act on the left as follows: 
\begin{itemize}
\item $A[y]\otimes\tilde{E}_{11}^{2}\to\tilde{E}_{11}^{2}$ by 
\begin{align*}
A[y]\otimes_{A[y]}y_{1}y_{2}E^{2}[y] & \lra y_{1}y_{2}E^{2}[y]\\
\theta\otimes y_{1}y_{2}ee & \mapsto y_{1}y_{2}\theta.ee.
\end{align*}
\item $A[y]\otimes\tilde{E}_{12}^{2}\to\tilde{E}_{12}^{2}$ by 
\begin{align*}
A[y]\otimes_{A[y]}y_{1}y_{2}y_{3}E^{3}[y] & \lra y_{1}y_{2}y_{3}E^{3}[y]\\
\theta\otimes y_{1}y_{2}y_{3}eee & \mapsto y_{1}y_{2}y_{3}\theta.eee.
\end{align*}
\end{itemize}
\end{lem}

\begin{rem*}
We may confirm that the image of the action map $\tilde{E}_{21}^{2}\to\tilde{E}_{21}^{2}$
preserves the conditions for $G_{2}$: 

\begin{align*}
\xi.\theta-\_\otimes e_{1}.\theta & =y_{2}\xi_{1}.\theta,\\
\xi_{1}.\theta & =\delta(\_\otimes e_{2}).\theta+(y_{1}\xi_{2}).\theta\\
 & =\delta(\_\otimes e_{2}.\theta)+y_{1}(\xi_{2}.\theta),
\end{align*}
and the $e_{\ell}$ relation:
\[
e_{1}.\theta-e_{2}.\theta=y_{1}e'.\theta.
\]
\end{rem*}
\begin{lem}
Generators in $G_{1}^{\op}\subset C$ act on the right on $\tilde{E}^{2}$
as follows: 
\begin{itemize}
\item $\tilde{E}_{12}^{2}\otimes G_{1}^{\op}\to\tilde{E}_{12}^{2}$ by 
\begin{align*}
y_{1}y_{2}y_{3}E^{3}[y]\otimes_{G_{1}^{\op}}G_{1}^{\op} & \lra y_{1}y_{2}y_{3}E^{3}[y]\\
y_{1}y_{2}y_{3}eee\otimes(\theta,\varphi) & \mapsto E^{2}\varphi(y_{1}y_{2}y_{3}eee)
\end{align*}
\item $\tilde{E}_{22}^{2}\otimes G_{1}^{\op}\to\tilde{E}_{22}^{2}$ by 
\begin{align*}
G_{3}\otimes_{G_{1}^{\op}}G_{1}^{\op} & \lra G_{3}\\
\left(ee_{1},ee_{2},ee_{3},\chi\right)\otimes(\theta,\varphi) & \mapsto\left(E\varphi(ee_{1}),E\varphi(ee_{2}),ee_{3}.\theta,E^{2}\varphi\circ\chi\right).
\end{align*}
\end{itemize}
They act on the left as follows:
\begin{itemize}
\item $G_{1}^{\op}\otimes\tilde{E}_{21}^{2}\to\tilde{E}_{21}^{2}$ by 
\begin{align*}
G_{1}^{\op}\otimes_{G_{1}^{\op}}G_{2} & \lra G_{2}\\
(\theta,\varphi)\otimes(e_{1},e_{2},\xi) & \mapsto\left(\theta.e_{1},\theta.e_{2},\xi\circ\varphi\right)
\end{align*}
\item $G_{1}^{\op}\otimes\tilde{E}_{22}^{2}\to\tilde{E}_{22}^{2}$ by 
\begin{align*}
G_{1}^{\op}\otimes_{G_{1}^{\op}}G_{3} & \lra G_{3}\\
(\theta,\varphi)\otimes\left(ee_{1},ee_{2},ee_{3},\chi\right) & \mapsto\left(\theta.ee_{1},\theta.ee_{2},\theta.ee_{3},\chi\circ\varphi\right).
\end{align*}
\end{itemize}
\end{lem}

\begin{rem}
We may confirm that the image of the right action map $\tilde{E}_{22}^{2}\otimes G_{1}^{\op}\to\tilde{E}_{22}^{2}$
preserves the conditions for $G_{3}$: 
\begin{align*}
E^{2}\varphi\circ\chi & =\_\otimes E\varphi(ee_{1})+E^{2}\varphi(\chi-\_\otimes ee_{1})\\
 & =\_\otimes E\varphi(ee_{1})+y_{3}\bigl(E^{2}\varphi\circ\chi_{1}\bigr),\\
E^{2}\varphi\circ\chi_{1} & =\tau E(\_\otimes E\varphi(ee_{2}))+y_{2}E^{2}\varphi\circ\chi_{1}'\\
 & =\tau E\circ E^{2}(\_.\theta+y_{1}\varphi_{1})\circ(\_\otimes ee_{2})+y_{2}E^{2}\varphi\circ\chi_{1}',\\
E^{2}\varphi\circ\chi_{1}' & =E^{2}(\_.\theta)\circ E\tau\circ\tau E(\_\otimes ee_{3})+y_{1}E^{2}\varphi_{1}\circ\chi_{1}'+y_{1}\chi''.\theta\\
 & =E\tau\circ\tau E(\_\otimes ee_{3}.\theta)+y_{1}\bigl(\chi''.\theta+E^{2}\varphi_{1}\circ\chi_{1}'\bigr).
\end{align*}
 And the $ee_{\ell}$ relations: 
\begin{align*}
E\varphi(ee_{1})-E\varphi(ee_{2}) & =y_{2}E\varphi(ee'),\\
ee_{3}.\theta-E\varphi(ee_{2}) & =(ee_{3}-ee_{2}).\theta-y_{1}E\varphi_{1}(ee_{2})\\
 & =y_{1}\left(ee''.\theta-E\varphi_{1}(ee_{2})\right),\\
\delta(ee_{3}.\theta)-E\varphi(ee_{1}) & =y_{2}\tau(ee_{3}).\theta+(ee_{3}-ee_{1}).\theta-y_{1}E\varphi_{1}(ee_{1})\\
 & =y_{2}\tau(ee_{3}).\theta+y_{1}ee''.\theta-y_{2}ee'.\theta-y_{1}E\varphi_{1}(ee_{1})\\
 & =y_{1}\left(y_{2}\overline{ee}.\theta+ee''.\theta-E\varphi_{1}(ee_{1})\right).
\end{align*}
 Similarly we may confirm that the image of the left action map $G_{1}^{\op}\otimes\tilde{E}_{21}^{2}\to\tilde{E}_{21}^{2}$
lies in $G_{2}$: 
\begin{align*}
\xi\circ\varphi & =\varphi(-)\otimes e_{1}+y_{2}\xi_{1}\circ\varphi\\
 & =\_\otimes\theta.e_{1}+y_{2}\bigl(\varphi_{1}(-)\otimes e_{1}+\xi_{1}\circ\varphi\bigr),\\
\xi_{1}\circ\varphi+\varphi_{1}(-)\otimes e_{1} & =\tau(\_\otimes e_{2})\circ\varphi+y_{1}\xi'\circ\varphi+\varphi_{1}(-)\otimes e_{1}\\
 & =\tau(\_\otimes\theta.e_{2})+\tau y_{2}(\varphi_{1}(-)\otimes e_{2})+y_{1}\xi'\circ\varphi+\varphi_{1}(-)\otimes e_{1}\\
 & =\tau(\_\otimes\theta.e_{2})+y_{1}\left(\tau(\varphi_{1}(-)\otimes e_{2})+\varphi_{1}(-)\otimes e'+\xi'\circ\varphi\right).
\end{align*}
And the $e_{\ell}$ relation: 
\[
\theta.e_{1}-\theta.e_{2}=y_{1}\theta.e'.
\]
 And the image of the left action map $G_{1}^{\op}\otimes\tilde{E}_{22}^{2}\to\tilde{E}_{22}^{2}$
lies in $G_{3}$: 
\begin{align*}
\chi\circ\varphi & =\varphi(-)\otimes ee_{1}+y_{3}\chi_{1}\circ\varphi\\
 & =\_\otimes\theta.ee_{1}+y_{3}\bigl(\varphi_{1}\otimes ee_{1}+\chi_{1}\circ\varphi\bigr),\\
\chi_{1}\circ\varphi & =\tau E(\_\otimes\theta.ee_{2})+\tau Ey_{3}(\varphi_{1}\otimes ee_{2})+y_{2}\chi_{1}'\circ\varphi,\\
\varphi_{1}\otimes ee_{1}+\chi_{1}\circ\varphi & =\tau E(\_\otimes\theta.ee_{2})+y_{2}\Bigl(\tau E(\varphi_{1}\otimes ee_{2})+\varphi_{1}\otimes ee'+\chi_{1}'\circ\varphi\Bigr),\\
\chi_{1}'\circ\varphi & =E\tau\circ\tau E(\_\otimes\theta.ee_{3})+E\tau\circ\tau E\circ y_{3}(\varphi_{1}\otimes ee_{3})+y_{1}\chi''\circ\varphi\\
 & =E\tau\circ\tau E(\_\otimes\theta.ee_{3})+y_{1}(E\tau\circ\tau E)(\varphi_{1}\otimes ee_{3})\\
 & \qquad-\tau E(\varphi_{1}\otimes ee_{3})-E\tau(\varphi_{1}\otimes ee_{3})+y_{1}\chi''\circ\varphi,
\end{align*}
 
\begin{multline*}
\tau E(\varphi_{1}\otimes ee_{2})+\varphi_{1}\otimes ee'+\chi_{1}'\circ\varphi=\\
E\tau\circ\tau E(\_\otimes\theta.ee_{3})+y_{1}\Bigl((E\tau\circ\tau E)(\varphi_{1}\otimes ee_{3})-\tau E(\varphi_{1}\otimes ee'')-\varphi_{1}\otimes\overline{ee}+\chi''\circ\varphi\Bigr).
\end{multline*}
And the $ee_{\ell}$ relations: 
\begin{align*}
\theta.ee_{1}-\theta.ee_{2} & =y_{2}\theta.ee'\\
\theta.ee_{3}-\theta.ee_{2} & =y_{1}\theta.ee''\\
\delta(\theta.ee_{3})-\theta.ee_{1} & =y_{1}\theta.ee'''.
\end{align*}
\end{rem}

\begin{lem}
Generators in $y_{1}E[y]\subset C$ act on the right on $\tilde{E}^{2}$
as follows: 
\begin{itemize}
\item $\tilde{E}_{11}^{2}\otimes y_{1}E[y]\to\tilde{E}_{12}^{2}$ by 
\begin{align*}
y_{1}y_{2}E^{2}[y]\otimes_{A[y]}y_{1}E[y] & \lra y_{1}y_{2}y_{3}E^{3}[y]\\
y_{1}y_{2}ee\otimes y_{1}e & \mapsto y_{1}y_{2}y_{3}(ee\otimes e)
\end{align*}
\item $\tilde{E}_{21}^{2}\otimes y_{1}E[y]\to\tilde{E}_{22}^{2}$ by 
\begin{align*}
G_{2}\otimes_{A[y]}y_{1}E[y] & \lra G_{3}\\
\left(e_{1},e_{2},\xi\right)\otimes y_{1}e & \mapsto\left(e_{1}\otimes y_{1}e,e_{2}\otimes y_{1}e,0,\xi(-)\otimes y_{1}e\right).
\end{align*}
\end{itemize}
They act on the left as follows: 
\begin{itemize}
\item $y_{1}E[y]\otimes\tilde{E}_{21}^{2}\to\tilde{E}_{11}^{2}$ by 
\begin{align*}
y_{1}E[y]\otimes_{G_{1}^{\op}}G_{2} & \lra y_{1}y_{2}E^{2}[y]\\
y_{1}e\otimes(e_{1},e_{2},\xi) & \mapsto\xi(y_{1}e)
\end{align*}
\item $y_{1}E[y]\otimes\tilde{E}_{22}^{2}\to\tilde{E}_{12}^{2}$ by 
\begin{align*}
y_{1}E[y]\otimes_{G_{1}^{\op}}G_{3} & \lra y_{1}y_{2}y_{3}E^{3}[y]\\
y_{1}e\otimes(ee_{1},ee_{2},ee_{3},\chi) & \mapsto\chi(y_{1}e).
\end{align*}
\end{itemize}
\end{lem}

\begin{rem*}
We may confirm that the image of the right action map $\tilde{E}_{21}^{2}\otimes y_{1}E[y]\to\tilde{E}_{22}^{2}$
preserves the conditions for $G_{3}$: 
\begin{align*}
\chi & =\xi\otimes y_{1}e,\\
\chi-\_\otimes e_{1}\otimes y_{1}e & =y_{1}y_{3}(\xi_{1}\otimes e),\\
\chi-\delta E(\_\otimes e_{2}\otimes y_{1}e) & =\bigl(\xi-\delta(\_\otimes e_{2})\bigr)\otimes y_{1}e\\
 & =y_{1}y_{2}(\xi_{2}\otimes e).
\end{align*}
 Similarly we may confirm that the image of the left action map $y_{1}E[y]\otimes\tilde{E}_{21}^{2}\to\tilde{E}_{11}^{2}$
lies in $y_{1}y_{2}E^{2}[y]$: 
\begin{align*}
\xi\circ y_{1} & =y_{2}\bigl(\_\otimes e_{1}+\xi_{1}\circ y_{1}\bigr),\\
\xi_{1}\circ y_{1} & =\tau y_{2}(\_\otimes e_{2})+y_{1}\xi'\circ y_{1}\\
 & =y_{1}\bigl(\tau(\_\otimes e_{2})+\xi'\circ y_{1}\bigr)-\_\otimes e_{2},\\
\xi\circ y_{1} & =y_{2}\Bigl(y_{1}\bigl(\tau(\_\otimes e_{2})+\xi'\circ y_{1}\bigr)+\_\otimes(e_{1}-e_{2})\Bigr)\\
 & =y_{1}y_{2}\Bigl(\tau(\_\otimes e_{2})+\_\otimes e'+\xi'\circ y_{1}\Bigr).
\end{align*}
 And the image of the left action map $y_{1}E[y]\otimes\tilde{E}_{22}^{2}\to\tilde{E}_{12}^{2}$
lies in $y_{1}y_{2}y_{3}E^{3}[y]$: 
\begin{align*}
\chi\circ y_{1} & =y_{3}\bigl(\_\otimes ee_{1}+\chi_{1}\circ y_{1}\bigr)\\
\chi_{1}\circ y_{1} & =-\tau Ey_{3}y_{1}(\_\otimes ee'')\\
 & \qquad+E\delta\circ\tau Ey_{3}(\_\otimes ee_{3})+y_{1}y_{2}\chi''\circ y_{1}\\
 & =-\tau Ey_{3}y_{1}(\_\otimes ee'')+E\delta\circ y_{2}\tau E(\_\otimes ee_{3})\\
 & \qquad-E\delta(\_\otimes ee_{3})+y_{1}y_{2}\chi''\circ y_{1}\\
 & =-y_{2}\tau Ey_{1}(\_\otimes ee'')+y_{1}(\_\otimes ee'')+y_{1}y_{2}E\tau\circ\tau E(\_\otimes ee_{3})\\
 & \qquad-y_{1}(\_\otimes ee''')-\_\otimes ee_{1}+y_{1}y_{2}\chi''\circ y_{1}\\
\chi\circ y_{1} & =y_{3}y_{2}y_{1}\Bigl(-\tau E(\_\otimes ee'')+E\tau\circ\tau E(\_\otimes ee_{3})+\chi''\circ y_{1}\Bigr)\\
 & \qquad+y_{3}y_{1}\bigl(\_\otimes ee''-\_\otimes ee'''\bigr)\\
 & =y_{3}y_{2}y_{1}\Bigl(-\tau E(\_\otimes ee'')+E\tau\circ\tau E(\_\otimes ee_{3})-\_\otimes\overline{ee}+\chi''\circ y_{1}\Bigr).
\end{align*}
\end{rem*}
\begin{lem}
Generators in $F[y]\subset C$ act on the right on $\tilde{E}^{2}$
as follows: 
\begin{itemize}
\item $\tilde{E}_{12}^{2}\otimes F[y]\to\tilde{E}_{11}^{2}$ by 
\begin{align*}
y_{1}y_{2}y_{3}E^{3}[y]\otimes_{G_{1}^{\op}}F[y] & \lra y_{1}y_{2}E^{2}[y]\\
y_{1}y_{2}y_{3}eee\otimes f & \mapsto y_{1}y_{2}E^{2}f(y_{1}eee)
\end{align*}
\item $\tilde{E}_{22}^{2}\otimes F[y]\to\tilde{E}_{21}^{2}$ by 
\begin{align*}
G_{3}\otimes_{G_{1}^{\op}}F[y] & \lra G_{2}\\
\left(ee_{1},ee_{2},ee_{3},\chi\right)\otimes f & \mapsto\left(Ef(ee_{1}),Ef(ee_{2}),E^{2}f\circ\chi\right).
\end{align*}
\end{itemize}
They act on the left as follows: 
\begin{itemize}
\item $F[y]\otimes\tilde{E}_{11}^{2}\to\tilde{E}_{21}^{2}$ by 
\begin{align*}
F[y]\otimes_{A[y]}y_{1}y_{2}E^{2}[y] & \lra G_{2}\\
f\otimes y_{1}y_{2}ee & \mapsto\left(0,0,f(-).y_{1}y_{2}ee\right)
\end{align*}
\item $F[y]\otimes\tilde{E}_{12}^{2}\to\tilde{E}_{22}^{2}$ by 
\begin{align*}
F[y]\otimes_{A[y]}y_{1}y_{2}y_{3}E^{3}[y] & \lra G_{3}\\
f\otimes y_{1}y_{2}y_{3}eee & \mapsto\left(0,0,0,f(-).y_{1}y_{2}y_{3}eee\right).
\end{align*}
\end{itemize}
\end{lem}

\begin{rem*}
We may observe that the image of the right action map $\tilde{E}_{22}^{2}\otimes F[y]\to\tilde{E}_{21}^{2}$
preserves the conditions for $G_{3}$: 
\begin{align*}
E^{2}f\circ\chi-\_\otimes Ef(ee_{1}) & =E^{2}f\circ(\chi-\_\otimes ee_{1})\\
 & =y_{2}E^{2}f\circ\chi_{1},\\
E^{2}f\circ\chi-\delta(\_\otimes Ef(ee_{2})) & =E^{2}f\circ(\chi-\delta E(\_\otimes ee_{2}))\\
 & =E^{2}f\circ y_{2}\chi_{2}\\
 & =y_{1}E^{2}f\circ\chi_{2},
\end{align*}
 and the $ee_{\ell}$ relation: 

\begin{align*}
Ef(ee_{1}-ee_{2}) & =Ef(y_{2}ee'')\\
 & =y_{1}Ef(ee'').
\end{align*}
 It is trivial to check the conditions for the images of the left
action maps $F[y]\otimes\tilde{E}_{11}^{2}\to\tilde{E}_{21}^{2}$
and $F[y]\otimes\tilde{E}_{12}^{2}\to\tilde{E}_{22}^{2}$.
\end{rem*}
\begin{proof}[Proof of Lemma \ref{lem:tilde-tau-is-endo}]
 The reader may now check that $\tilde{\tau}$ defined in §\ref{subsec:Defs_x_tau}
is equivariant over the left and right \emph{$C$} actions. These
checks are completely mechanical using the formulas just given.
\end{proof}

\subsection{Hecke relations}

\subsubsection{$\tilde{x}$ and $\tilde{\tau}$ satisfy Hecke relations}

These checks are also mechanical, but we write them out because they
are important.
\begin{prop}
On each component $\tilde{E}_{ij}^{2}$, the maps $\tilde{x}$ and
$\tilde{\tau}$ defined in §\ref{subsec:Defs_x_tau} satisfy 
\begin{align*}
\ee\xx\circ\tt-\tt\circ\xx\ee & =\idop\\
\tilde{\tau}\circ\tilde{E}\tilde{x}-\tilde{x}\tilde{E}\circ\tilde{\tau} & =\idop.
\end{align*}
\end{prop}

\begin{proof}
On the first row, $\tilde{E}_{11}^{2}$ and $\tilde{E}_{12}^{2}$,
the relations follow from the corresponding relations between $x$
and $\tau$.

On $\tilde{E}_{21}^{2}$ presented as $G_{2}$, we have:
\begin{align*}
\ee\xx\circ\tt:(e_{1},e_{2},\xi) & \mapsto(xe',ye',Ex\circ\tau\circ\xi)\\
\tt\circ\xx\ee:(e_{1},e_{2},\xi) & \mapsto(ye'-e_{2},ye'-e_{2},\tau\circ xE\circ\xi)\\
\tt\circ\ee\xx:(e_{1},e_{2},\xi) & \mapsto(e_{2}+xe',e_{2}+xe',\tau\circ Ex\circ\xi)\\
\xx\ee\circ\tt:(e_{1},e_{2},\xi) & \mapsto(ye',xe',xE\circ\tau\circ\xi),
\end{align*}
from which 
\begin{align*}
\ee\xx\circ\tt-\tt\circ\xx\ee:(e_{1},e_{2},\xi) & \mapsto\bigl(y_{1}e'+e_{2},e_{2},(Ex\circ\tau-\tau\circ xE)\circ\xi\bigr)\\
 & \qquad=(e_{1},e_{2},\xi),
\end{align*}
and similarly for the other relation.

On $\tilde{E}_{22}^{2}$ presented as $G_{3}$, we have:
\begin{align*}
\ee\xx\circ\tt:(ee_{1},ee_{2},ee_{3},\chi) & \mapsto(xE(ee'),yee',Ex\circ\tau(ee_{3}),ExE\circ\tau E\circ\chi)\\
\tt\circ\xx\ee:(ee_{1},ee_{2},ee_{3},\chi) & \mapsto(yee'-ee_{2},yee'-ee_{2},\tau\circ xE(ee_{3}),\tau E\circ xE^{2}\circ\chi)\\
\tt\circ\ee\xx:(ee_{1},ee_{2},ee_{3},\chi) & \mapsto(ee_{1}+yee',ee_{1}+yee',\tau\circ Ex(ee_{3}),\tau E\circ ExE\circ\chi)\\
\xx\ee\circ\tt:(ee_{1},ee_{2},ee_{3},\chi) & \mapsto(yee',xE(ee'),xE\circ\tau(ee_{3}),xE^{2}\circ\tau E\circ\chi),
\end{align*}
and so
\begin{multline*}
\ee\xx\circ\tt-\tt\circ\xx\ee:(ee_{1},ee_{2},ee_{3},\chi)\mapsto\\
\bigl(y_{2}ee'+ee_{2},ee_{2},(Ex\circ\tau-\tau\circ xE)(ee_{3}),(ExE\circ\tau E-\tau E\circ xE^{2})\circ\chi\bigr)\\
=(ee_{1},ee_{2},ee_{3},\chi),\qquad\qquad\qquad
\end{multline*}
and similarly for the other relation.
\end{proof}

\subsubsection{$\tilde{\tau}^{2}=0$}

This is clear.

\subsubsection{$\tilde{\tau}$ satisfies the braid relation}

In this section we give formulas defining $k$-module endomorphisms
$\tilde{\tau}_{1}$ and $\tilde{\tau}_{2}$ of the components of the
matrix parametrization of $\tilde{E}^{3}$. We show that these endomorphisms
satisfy the braid relations. Then we argue that they correspond to
the maps $\tilde{E}\tilde{\tau}$ and $\tilde{\tau}\tilde{E}$ induced
on the same bimodule components. This will complete our proof that
$\tilde{x}$ and $\tilde{\tau}$ satisfy the nil affine Hecke relations
in $\mathcal{U}^{+}$.
\begin{lem}
\label{lem:within-G3} Let us be given $(ee_{1},ee_{2},ee_{3},\chi)\in G_{3}$
with $ee''$ defined as in Def.~\ref{def:G_3}. Then 
\[
(\tau(ee_{1}),-ee'',-ee'',E\tau\circ\chi)\in E^{2}[y]^{\oplus3}\oplus\Hom_{A}(_{A}E,E^{3})[y]
\]
 also lies in $G_{3}$.
\end{lem}

\begin{proof}
The reader may check this directly. In Prop.~\ref{prop:Correspondence-tau_i}
we will interpret this element as the image of $(ee_{1},ee_{2},ee_{3},\chi)$
under $\tilde{E}\tilde{\tau}$, and it must therefore lie in $G_{3}$.
\end{proof}
\begin{lem}
\label{lem:within-G4} Let us be given $(eee_{1},eee_{2},eee_{3},eee_{4},\psi)\in G_{4}$
with $eee^{(\ell)}$ defined as in §\ref{def:G_4}. Then the following
elements of $E^{3}[y]^{\oplus4}\oplus\Hom_{A}(_{A}E,E^{4})[y]$ also
lie in $G_{4}$: 
\begin{gather*}
(\tau E(eee_{1}),eee^{(2)},eee^{(2)},E\tau(eee_{4}),E\tau E\circ\psi),\\
(eee^{(4)},eee^{(4)},\tau E(eee_{3}),\tau E(eee_{4}),\tau E^{2}\circ\psi).
\end{gather*}
\end{lem}

\begin{proof}
The reader may check this directly. In Prop.~\ref{prop:Correspondence-tau_i}
we will interpret these elements as the images of $(eee_{1},eee_{2},eee_{3},eee_{4},\psi)$
under $\tilde{E}\tilde{\tau}$ and $\tilde{\tau}\tilde{E}$ respectively,
and they must therefore lie in $G_{4}$.
\end{proof}
\begin{defn}
\label{def:tau_1=002013tau_2} Let $\tilde{\tau}_{1},\tilde{\tau}_{2}$
be $k$-module maps defined on $\tilde{E}_{ij}^{3}$, presented as
in §\ref{subsec:Desc-of-E^3}, as follows: 
\begin{itemize}
\item on $\ee_{11}^{3}$:
\begin{itemize}
\item $\tilde{\tau}_{1}$ acts by $E\tau$
\item $\tilde{\tau}_{2}$ by $\tau E$
\end{itemize}
\item on $\ee_{12}^{3}$:
\begin{itemize}
\item $\tilde{\tau}_{1}$ by $E\tau E$
\item $\tilde{\tau}_{2}$ by $\tau E^{2}$
\end{itemize}
\item on $\ee_{21}^{3}$:
\begin{itemize}
\item $\tilde{\tau}_{1}$ by $(ee_{1},ee_{2},ee_{3},\chi)\mapsto(\tau(ee_{1}),-ee'',-ee'',E\tau\circ\chi)$
\item $\tilde{\tau}_{2}$ by $(ee_{1},ee_{2},ee_{3},\chi)\mapsto(ee',ee',\tau(ee_{3}),\tau E\circ\chi)$,
\\
$\phantom{.}$ \quad{} i.e.~$\tilde{\tau}$ as defined above on
$G_{3}$ considered as $\tilde{E}_{22}^{2}$
\end{itemize}
\item on $\ee_{22}^{3}$:
\begin{itemize}
\item $\tilde{\tau}_{1}$ by $(eee_{1},eee_{2},eee_{3},eee_{4},\psi)\mapsto$
\[
(\tau E(eee_{1}),eee^{(2)},eee^{(2)},E\tau(eee_{4}),E\tau E\circ\psi)
\]
\item $\tilde{\tau}_{2}$ by $(eee_{1},eee_{2},eee_{3},eee_{4},\psi)\mapsto$
\[
(eee^{(4)},eee^{(4)},\tau E(eee_{3}),\tau E(eee_{4}),\tau E^{2}\circ\psi).
\]
\end{itemize}
\end{itemize}
\end{defn}

\begin{prop}
\label{prop:tau_i-braid} The $\tilde{\tau}_{i}$ satisfy $\tilde{\tau}_{1}\circ\tilde{\tau}_{2}\circ\tilde{\tau}_{1}=\tilde{\tau}_{2}\circ\tilde{\tau}_{1}\circ\tilde{\tau}_{2}$.
\end{prop}

\begin{proof}
On $\tilde{E}_{1j}^{3}$ the claim follows from the $\tau_{i}$ braid
relation. On $\tilde{E}_{21}^{3}=G_{3}$ we have:
\begin{flalign*}
 & (ee_{1},ee_{2},ee_{3},\chi)\overset{\tilde{\tau}_{1}}{\longmapsto}\\
 & (\tau(ee_{1}),-ee'',-ee'',E\tau\circ\chi)\overset{\tilde{\tau}_{2}}{\longmapsto}\\
 & (-\overline{ee}-\tau(ee'''),-\overline{ee}-\tau(ee'''),-\tau(ee''),\tau E\circ E\tau\circ\chi)\overset{\tilde{\tau}_{1}}{\longmapsto}\\
 & (-\tau(\overline{ee}),-\tau(\overline{ee}),-\tau(\overline{ee}),E\tau\circ\tau E\circ E\tau\circ\chi)
\end{flalign*}

and 
\begin{flalign*}
 & (ee_{1},ee_{2},ee_{3},\chi)\overset{\tilde{\tau}_{2}}{\longmapsto}\\
 & (ee',ee',\tau(ee_{3}),\tau E\circ\chi)\overset{\tilde{\tau}_{1}}{\longmapsto}\\
 & (\tau(ee'),-\overline{ee},-\overline{ee},E\tau\circ\tau E\circ\chi)\overset{\tilde{\tau}_{2}}{\longmapsto}\\
 & (-\tau(\overline{ee}),-\tau(\overline{ee}),-\tau(\overline{ee}),\tau E\circ E\tau\circ\tau E\circ\chi).
\end{flalign*}

On $\tilde{E}_{22}^{3}=G_{4}$ we have: 
\begin{flalign*}
 & (eee_{1},eee_{2},eee_{3},eee_{4},\psi)\overset{\tilde{\tau}_{1}}{\longmapsto}\\
 & (\tau E(eee_{1}),eee^{(2)},eee^{(2)},E\tau(eee_{4}),E\tau E\circ\psi)\overset{\tilde{\tau}_{2}}{\longmapsto}\\
 & (\tau E(eee^{(5)})+\overline{eee},\tau E(eee^{(5)})+\overline{eee},\tau E(eee^{(2)}),\tau E\circ E\tau(eee_{4}),\tau E^{2}\circ E\tau E\circ\psi)\overset{\tilde{\tau}_{1}}{\longmapsto}\\
 & (\tau E(\overline{eee}),\tau E(\overline{eee}),\tau E(\overline{eee}),E\tau\circ\tau E\circ E\tau(eee_{4}),E\tau E\circ\tau E^{2}\circ E\tau E\circ\psi)
\end{flalign*}

and 
\begin{flalign*}
 & (eee_{1},eee_{2},eee_{3},eee_{4},\psi)\overset{\tilde{\tau}_{2}}{\longmapsto}\\
 & (eee^{(4)},eee^{(4)},\tau E(eee_{3}),\tau E(eee_{4}),\tau E^{2}\circ\psi)\overset{\tilde{\tau}_{1}}{\longmapsto}\\
 & (\tau E(eee^{(4)}),\overline{eee},\overline{eee},E\tau\circ\tau E(eee_{4}),E\tau E\circ\tau E^{2}\circ\psi)\overset{\tilde{\tau}_{2}}{\longmapsto}\\
 & (\tau E(\overline{eee}),\tau E(\overline{eee}),\tau E(\overline{eee}),\tau E\circ E\tau\circ\tau E(eee_{4}),\tau E^{2}\circ E\tau E\circ\tau E^{2}\circ\psi).
\end{flalign*}
\end{proof}
The remaining goal of this section is to relate the $\tilde{\tau}_{i}$
just defined to the $\tilde{\tau}$ acting on $\tilde{E}$ as described
in §\ref{subsec:Defs_x_tau}. The latter is known to be a $(C,C)$-bimodule
morphism.
\begin{prop}
\label{prop:Correspondence-tau_i} Under the isomorphism of Lemma
\ref{lem:Nat-Isos_E'^n}, namely 
\[
\tilde{E}^{3}\iso\Hom_{K^{b}(B)}(X,E'^{3}X),
\]
 the maps $\tilde{E}\tilde{\tau}$ and $\tilde{\tau}\tilde{E}$ on
$\tilde{E}^{3}$ correspond to $\tilde{\tau}_{1}$ and $\tilde{\tau}_{2}$
of Definition \ref{def:tau_1=002013tau_2}.
\end{prop}

\begin{cor}
Lemmas \ref{lem:within-G3} and \ref{lem:within-G4} follow.
\end{cor}

\begin{cor}
Proposition \ref{prop:tau_i-braid} implies $\tilde{E}\tilde{\tau}\circ\tilde{\tau}\tilde{E}\circ\tilde{E}\tilde{\tau}=\tilde{\tau}\tilde{E}\circ\tilde{E}\tilde{\tau}\circ\tilde{\tau}\tilde{E}$.
\end{cor}

\begin{proof}[Proof of the proposition.]
We consider the tensor product $\tilde{E}\otimes_{C}\tilde{E}^{2}$
formed according to the procedure of §\ref{subsec:matrix-tensor-products},
and study the endofunctor $\tilde{E}\tilde{\tau}$ as in Lemma \ref{lem:tensor-components-induced},
and similarly for $\tilde{E}^{2}\otimes_{C}\tilde{E}$ and $\tilde{\tau}\tilde{E}$.
From Lemma \ref{lem:Nat-Isos_E'^n}, we have isomorphisms 
\begin{align*}
\Hom_{K^{b}(B)}(X,E'X)\otimes_{C}\Hom_{K^{b}(B)}(X,E'^{2}X) & \iso\Hom_{K^{b}(B)}(X,E'^{3}X)\\
\Hom_{K^{b}(B)}(X,E'^{2}X)\otimes_{C}\Hom_{K^{b}(B)}(X,E'X) & \iso\Hom_{K^{b}(B)}(X,E'^{3}X)
\end{align*}
 associated with the products 
\begin{align*}
\tilde{E}\otimes_{C}\tilde{E}^{2} & =\tilde{E}^{3}\\
\tilde{E}^{2}\otimes_{C}\tilde{E} & =\tilde{E}^{3}.
\end{align*}
The maps are given by 
\begin{align*}
f\otimes g & \mapsto E'g\circ f\\
f\otimes g & \mapsto E'^{2}g\circ f.
\end{align*}
 These isomorphisms determine actions of $\tilde{E}\tilde{\tau}$
and $\tilde{\tau}\tilde{E}$ on $\Hom_{K^{b}(B)}(X,E'^{3}X)$ that
we may compare to the $\tilde{\tau}_{1}$ and $\tilde{\tau}_{2}$
defined there by components.

The components $\tilde{E}_{ij}$ and $\tilde{E}_{ij}^{2}$ are $(\End(X_{i})^{\op},\End(X_{j})^{\op})$-bimodules,
and $\tilde{\tau}$ gives bimodule endomorphisms $\tilde{\tau}_{\mid ij}$
of the latter. These induce endomorphisms $(\tilde{E}\tilde{\tau})_{\mid ijk}^{1\mid2}$
of 
\[
\tilde{E}_{ijk}^{1\mid2}=\tilde{E}_{ij}\otimes_{\End(X_{j})^{\op}}\tilde{E}_{jk}^{2},
\]
 as in §\ref{subsec:matrix-tensor-products}. We know that $\tilde{E}_{ik}^{3}$
is canonically isomorphic to a quotient of $\tilde{E}_{i1k}^{1\mid2}\oplus\tilde{E}_{i2k}^{1\mid2}$,
and that $\left(\begin{smallmatrix}(\tilde{E}\tilde{\tau})_{\mid i1k}^{1\mid2} & 0\\
0 & (\tilde{E}\tilde{\tau})_{\mid i2k}^{1\mid2}
\end{smallmatrix}\right)$ acting on $\tilde{E}_{i1k}^{1\mid2}\oplus\tilde{E}_{i2k}^{1\mid2}$
descends to $\tilde{E}_{ik}^{3}$, where it gives the components of
$\tilde{E}\tilde{\tau}$. Here it may be compared directly with $\tilde{\tau}_{1}$
that we defined on $\tilde{E}_{ik}^{3}$. It therefore suffices for
our objective to check commutativity of the following diagrams labeled
$D_{1\mid2}(i,j,k)$, indexed by triples $(i,j,k)\in\{1,2\}^{3}$:

\noindent\parbox[t]{1\columnwidth}{%
\[\begin{tikzcd}[ampersand replacement=\&] 	\& {\tilde{E}_{ij}\otimes_{\mathrm{End}(X_j)^{\mathsf{op}}}\tilde{E}^2_{jk}} \& {} \& {\tilde{E}^3_{ik}} \\ 	{D_{1\mid 2}(i,j,k):} \\ 	\& {\tilde{E}_{ij}\otimes_{\mathrm{End}(X_j)^{\mathsf{op}}}\tilde{E}^2_{jk}} \&\& {\tilde{E}^3_{ik}.} 	\arrow["{(\tilde{E}\tilde{\tau})_{\mid ijk}^{1\mid 2}}", from=1-2, to=3-2] 	\arrow["{f\otimes g\mapsto E'g\circ f}"', from=1-2, to=1-4] 	\arrow["{\tilde{\tau}_{1{\mid ik}}}"', from=1-4, to=3-4] 	\arrow["{f\otimes g\mapsto E'g\circ f}", from=3-2, to=3-4] \end{tikzcd}\]%
} Exactly parallel considerations apply to the study of $\tilde{\tau}\tilde{E}$,
where the diagrams for $(i,j,k)$, now labeled $D_{2\mid1}(i,j,k)$,
instead involve maps $(\tilde{E}\tilde{\tau})_{\mid ijk}^{2\mid1}$
and $\tilde{\tau}_{2\mid ik}$.

Checking the diagrams will occupy the next three pages.
\begin{lem}
The diagrams $D_{1\mid2}(i,j,k)$ commute.
\end{lem}

\begin{proof}
We consider the diagrams in turn:
\begin{itemize}
\item Diagram $D_{1\mid2}(1,1,1)$:
\end{itemize}
Consider $(\tilde{E}\tilde{\tau})_{\mid111}^{1\mid2}\in\End(\tilde{E}_{11}\otimes\tilde{E}_{11}^{2})$.
Let $y_{1}e\in\tilde{E}_{11}$ and $y_{1}y_{2}ee\in\tilde{E}_{11}^{2}$.
The image of $y_{1}e\otimes y_{1}y_{2}ee$ in the top right corner
of the diagram is 
\[
E'(y_{1}y_{2}ee)\circ y_{1}e=y_{1}y_{2}y_{3}(e\otimes ee)\in\tilde{E}_{11}^{3}.
\]
 Here we can write out $E'(y_{1}y_{2}ee)=(y_{1}y_{2}ee,0,0,\_\otimes y_{1}y_{2}ee)\in G_{3}$.
On the other hand, $\tilde{\tau}(y_{1}y_{2}ee)=y_{1}y_{2}\tau(ee)$,
so the image of $(\tilde{E}\tilde{\tau})_{\mid111}^{1\mid2}\bigl(y_{1}e\otimes y_{1}y_{2}ee\bigr)$
is $y_{1}y_{2}y_{3}(e\otimes\tau(ee))\in\tilde{E}_{11}^{3}$, which
agrees with $\tilde{\tau}_{1}(y_{1}y_{2}y_{3}(e\otimes ee))$.
\begin{itemize}
\item Diagram $D_{1\mid2}(1,2,1)$:
\end{itemize}
Consider $(\tilde{E}\tilde{\tau})_{\mid121}^{1\mid2}\in\End(\tilde{E}_{12}\otimes\tilde{E}_{21}^{2})$.
Let $y_{1}y_{2}ee\in\tilde{E}_{12}$ and $(e_{1},e_{2},\xi)\in\tilde{E}_{21}^{2}$.
We have no established notation for $E'\bigl((e_{1},e_{2},\xi)\bigr)\in\Hom_{K^{b}(B)}(E'X_{2},E'^{2}X_{1})$.
It is nevertheless easy to check, by tracking \textquoteleft leading
terms\textquoteright{} of the upper rows, that 
\[
E'\bigl((e_{1},e_{2},\xi)\bigr)\circ y_{1}y_{2}ee=E\xi(y_{1}y_{2}ee)\in\tilde{E}_{11}^{3}.
\]
 This lies in $y_{1}y_{2}y_{3}E^{3}[y]$. Then $\tilde{\tau}\bigl((e_{1},e_{2},\xi)\bigr)=(e',e',\tau\circ\xi)$,
so $(\tilde{E}\tilde{\tau})_{\mid121}^{1\mid2}$ applied to $y_{1}y_{2}ee\otimes(e_{1},e_{2},\xi)$
and viewed in $\tilde{E}_{11}^{3}$ is $E\tau\circ E\xi(y_{1}y_{2}ee)$.
\begin{itemize}
\item Diagram $D_{1\mid2}(2,1,1)$:
\end{itemize}
Consider $(\tilde{E}\tilde{\tau})_{\mid211}^{1\mid2}\in\End(\tilde{E}_{21}\otimes\tilde{E}_{11}^{2})$.
Let $(\theta,\varphi)\in\tilde{E}_{21}$ and $y_{1}y_{2}ee\in\tilde{E}_{11}^{2}$.
This time we can write $E'(y_{1}y_{2}ee)=(y_{1}y_{2}ee,0,0,\_\otimes y_{1}y_{2}ee)$.
Then 
\[
E'(y_{1}y_{2}ee)\circ(\theta,\varphi)=(\theta y_{1}y_{2}ee,0,0,\varphi\otimes y_{1}y_{2}ee)\in\tilde{E}_{21}^{3}.
\]
Going around the diagram in either direction yields $(\theta y_{1}y_{2}\tau(ee),0,0,\varphi\otimes y_{1}y_{2}\tau(ee))$.
\begin{itemize}
\item Diagram $D_{1\mid2}(2,2,1)$:
\end{itemize}
Consider $(\tilde{E}\tilde{\tau})_{\mid221}^{1\mid2}\in\End(\tilde{E}_{22}\otimes\tilde{E}_{21}^{2})$.
Let $(e_{1},e_{2},\xi)\in\tilde{E}_{22}$ and $(\bar{e}_{1},\bar{e}_{2},\bar{\xi})\in\tilde{E}_{21}^{2}$.
We have no notation for $E'\left((\bar{e}_{1},\bar{e}_{2},\bar{\xi})\right)$.
One computes that 
\[
E'\left((\bar{e}_{1},\bar{e}_{2},\bar{\xi})\right)\circ(e_{1},e_{2},\xi)=(\bar{\xi}(e_{1}),e_{2}\otimes\bar{e}_{1},e_{2}\otimes\bar{e}_{2},E\bar{\xi}\circ\xi)\in\tilde{E}_{21}^{3}.
\]
 Traversing the diagram in either direction gives $(\tau\circ\bar{\xi}(e_{1}),e_{2}\otimes\bar{e}',e_{2}\otimes\bar{e}',E\tau\circ E\bar{\xi}\circ\xi)$.
\begin{itemize}
\item Diagram $D_{1\mid2}(1,1,2)$:
\end{itemize}
Consider $(\tilde{E}\tilde{\tau})_{\mid112}^{1\mid2}\in\End(\tilde{E}_{11}\otimes\tilde{E}_{12}^{2})$.
Let $y_{1}e\in\tilde{E}_{11}$ and $y_{1}y_{2}y_{3}eee\in\tilde{E}_{12}^{2}$.
Again by tracking \textquoteleft leading terms\textquoteright , one
checks that 
\[
E'(y_{1}y_{2}y_{3}eee)\circ y_{1}e=y_{1}\dots y_{4}(e\otimes eee)\in\tilde{E}_{12}^{3}.
\]
 Traversing the diagram in either direction gives $E\tau E(y_{1}\dots y_{4}e\otimes eee)$
which is $y_{1}\dots y_{4}(e\otimes\tau E(eee))$. 
\begin{itemize}
\item Diagram $D_{1\mid2}(1,2,2)$:
\end{itemize}
Consider $(\tilde{E}\tilde{\tau})_{\mid122}^{1\mid2}\in\End(\tilde{E}_{12}\otimes\tilde{E}_{22}^{2})$.
Let $y_{1}y_{2}ee\in\tilde{E}_{12}$ and $(ee_{1},ee_{2},ee_{3},\chi)\in\tilde{E}_{22}^{2}$.
Then check that 
\[
E'\left((ee_{1},ee_{2},ee_{3},\chi)\right)\circ y_{1}y_{2}ee=E\chi(y_{1}y_{2}ee)\in\tilde{E}_{12}^{3}.
\]
 Traversing the diagram in either direction gives $(E\tau E\circ E\chi)(y_{1}y_{2}ee)$.
\begin{itemize}
\item Diagram $D_{1\mid2}(2,1,2)$:
\end{itemize}
Consider $(\tilde{E}\tilde{\tau})_{\mid212}^{1\mid2}\in\End(\tilde{E}_{21}\otimes\tilde{E}_{12}^{2})$.
Let $(\theta,\varphi)\in\tilde{E}_{21}$ and $y_{1}y_{2}y_{3}eee\in\tilde{E}_{12}^{2}$.
Then check that 
\[
E'(y_{1}y_{2}y_{3}eee)\circ(\theta,\varphi)=(\theta y_{1}y_{2}y_{3}eee,0,0,0,\varphi\otimes y_{1}y_{2}y_{3}eee)\in\tilde{E}_{22}^{3}.
\]
 Traversing the diagram in either direction gives 
\[
\left(\tau E(\theta y_{1}y_{2}y_{3}eee),0,0,0,E\tau E\circ(\varphi\otimes y_{1}y_{2}y_{3}eee)\right).
\]
\begin{itemize}
\item Diagram $D_{1\mid2}(2,2,2)$:
\end{itemize}
Consider $(\tilde{E}\tilde{\tau})_{\mid222}^{1\mid2}\in\End(\tilde{E}_{22}\otimes\tilde{E}_{22}^{2})$.
Let $(e_{1},e_{2},\xi)\in\tilde{E}_{22}$ and $(ee_{1},ee_{2},ee_{3},\chi)\in\tilde{E}_{22}^{2}$.
Then check that 
\[
E'\bigl((ee_{1},ee_{2},ee_{3},\chi)\bigr)\circ(e_{1},e_{2},\xi)=\left(\chi(e_{1}),e_{2}\otimes ee_{1},e_{2}\otimes ee_{2},e_{2}\otimes ee_{3},E\chi\circ\xi\right)\in\tilde{E}_{22}^{3}.
\]
 Traversing the diagram in either direction gives 
\[
\bigl(\tau E(\chi(e_{1})),e_{2}\otimes ee',e_{2}\otimes ee',E\tau(e_{2}\otimes ee_{3}),E\tau E\circ E\chi\circ\xi\bigr).
\]
\end{proof}
\begin{lem}
The diagrams $D_{2\mid1}(i,j,k)$ commute.
\end{lem}

\begin{proof}
We consider the diagrams in turn:
\begin{itemize}
\item Diagram $D_{2\mid1}(1,1,1)$:
\end{itemize}
Consider $(\tilde{\tau}\tilde{E})_{\mid111}^{2\mid1}\in\End(\tilde{E}_{11}^{2}\otimes\tilde{E}_{11})$.
Let $y_{1}y_{2}ee\in\tilde{E}_{11}^{2}$ and $y_{1}e\in\tilde{E}_{11}$.
Then check that 
\[
E'^{2}(y_{1}e)\circ y_{1}y_{2}ee=y_{1}y_{2}y_{3}ee\otimes e\in\tilde{E}_{11}^{3}.
\]
 Traversing the diagram in either direction gives 
\[
y_{1}y_{2}y_{3}(\tau(ee)\otimes e).
\]
\begin{itemize}
\item Diagram $D_{2\mid1}(1,2,1)$:
\end{itemize}
Consider $(\tilde{\tau}\tilde{E})_{\mid121}^{2\mid1}\in\End(\tilde{E}_{12}^{2}\otimes\tilde{E}_{21})$.
Let $y_{1}y_{2}y_{3}eee\in\tilde{E}_{12}^{2}$ and $(\theta,\varphi)\in\tilde{E}_{21}$.
Then check that 
\[
E'^{2}\bigl((\theta,\varphi)\bigr)\circ y_{1}y_{2}y_{3}eee=E^{2}\varphi(y_{1}y_{2}y_{3}eee)\in\tilde{E}_{11}^{3}.
\]
 Traversing the diagram in either direction gives 
\[
(\tau E\circ E^{2}\varphi)(y_{1}y_{2}y_{3}eee).
\]
\begin{itemize}
\item Diagram $D_{2\mid1}(2,1,1)$:
\end{itemize}
Consider $(\tilde{\tau}\tilde{E})_{\mid211}^{2\mid1}\in\End(\tilde{E}_{21}^{2}\otimes\tilde{E}_{11})$.
Let $(e_{1},e_{2},\xi)\in\tilde{E}_{21}^{2}$ and $y_{1}e\in\tilde{E}_{11}$.
Then check that 
\[
E'^{2}(y_{1}e)\circ(e_{1},e_{2},\xi)=(e_{1}\otimes y_{1}e,e_{2}\otimes y_{1}e,0,\xi\otimes y_{1}e)\in\tilde{E}_{21}^{3}.
\]
 Traversing the diagram in either direction gives 
\[
(e'\otimes y_{1}e,e'\otimes y_{1}e,0,(\tau\circ\xi)\otimes y_{1}e).
\]
\begin{itemize}
\item Diagram $D_{2\mid1}(2,2,1)$:
\end{itemize}
Consider $(\tilde{\tau}\tilde{E})_{\mid221}^{2\mid1}\in\End(\tilde{E}_{22}^{2}\otimes\tilde{E}_{21})$.
Let $(ee_{1},ee_{2},ee_{3},\chi)\in\tilde{E}_{22}^{2}$ and $(\theta,\varphi)\in\tilde{E}_{21}$.
Then check that 
\[
E'^{2}\bigl((\theta,\varphi)\bigr)\circ(ee_{1},ee_{2},ee_{3},\chi)=(E\varphi(ee_{1}),E\varphi(ee_{2}),\theta ee_{3},E^{2}\varphi\circ\chi)\in\tilde{E}_{21}^{3}.
\]
 Traversing the diagram in either direction gives 
\[
\left(E\varphi(ee'),E\varphi(ee'),\theta\tau(ee_{3}),E^{2}\varphi\circ\tau E\circ\chi\right).
\]
\begin{itemize}
\item Diagram $D_{2\mid1}(1,1,2)$:
\end{itemize}
Consider $(\tilde{\tau}\tilde{E})_{\mid112}^{2\mid1}\in\End(\tilde{E}_{11}^{2}\otimes\tilde{E}_{12})$.
Let $y_{1}y_{2}ee\in\tilde{E}_{11}^{2}$ and $y_{1}y_{2}\bar{e}\bar{e}\in\tilde{E}_{12}$.
Then check that 
\[
E'^{2}(y_{1}y_{2}\bar{e}\bar{e})\circ y_{1}y_{2}ee=(y_{1}y_{2}ee)\otimes(y_{1}y_{2}\bar{e}\bar{e})=y_{1}\dots y_{4}(ee\otimes\bar{e}\bar{e})\in\tilde{E}_{12}^{3}.
\]
 Traversing the diagram in either direction gives 
\[
y_{1}\dots y_{4}(\tau(ee)\otimes\bar{e}\bar{e}).
\]
\begin{itemize}
\item Diagram $D_{2\mid1}(1,2,2)$:
\end{itemize}
Consider $(\tilde{\tau}\tilde{E})_{\mid122}^{2\mid1}\in\End(\tilde{E}_{12}^{2}\otimes\tilde{E}_{22})$.
Let $y_{1}y_{2}y_{3}eee\in\tilde{E}_{12}^{2}$ and $(e_{1},e_{2},\xi)\in\tilde{E}_{22}$.
Then check that 
\[
E'^{2}\bigl((e_{1},e_{2},\xi)\bigr)\circ y_{1}y_{2}y_{3}eee=E^{2}\xi(y_{1}y_{2}y_{3}eee)\in\tilde{E}_{12}^{3}.
\]
 Traversing the diagram in either direction gives 
\[
(\tau E^{2}\circ E^{2}\xi)(y_{1}y_{2}y_{3}eee).
\]
\begin{itemize}
\item Diagram $D_{2\mid1}(2,1,2)$:
\end{itemize}
Consider $(\tilde{\tau}\tilde{E})_{\mid212}^{2\mid1}\in\End(\tilde{E}_{21}^{2}\otimes\tilde{E}_{12})$.
Let $(e_{1},e_{2},\xi)\in\tilde{E}_{21}^{2}$ and $y_{1}y_{2}ee\in\tilde{E}_{12}$.
Then check that 
\[
E'^{2}(y_{1}y_{2}ee)\circ(e_{1},e_{2},\xi)=(e_{1}\otimes y_{1}y_{2}ee,e_{2}\otimes y_{1}y_{2}ee,0,0,\xi\otimes y_{1}y_{2}ee)\in\tilde{E}_{22}^{3}.
\]
 Traversing the diagram in either direction gives 
\[
\left(e'\otimes y_{1}y_{2}ee,e'\otimes y_{1}y_{2}ee,0,0,(\tau\circ\xi)\otimes y_{1}y_{2}ee\right).
\]
\begin{itemize}
\item Diagram $D_{2\mid1}(2,2,2)$:
\end{itemize}
Consider $(\tilde{\tau}\tilde{E})_{\mid222}^{2\mid1}\in\End(\tilde{E}_{22}^{2}\otimes\tilde{E}_{22})$.
Let $(ee_{1},ee_{2},ee_{3},\chi)\in\tilde{E}_{22}^{2}$ and $(e_{1},e_{2},\xi)\in\tilde{E}_{22}$.
Then check that 
\begin{multline*}
E'^{2}\bigl((e_{1},e_{2},\xi)\bigr)\circ(ee_{1},ee_{2},ee_{3},\chi)=\\
\left(E\xi(ee_{1}),E\xi(ee_{2}),ee_{3}\otimes e_{1},ee_{3}\otimes e_{2},E^{2}\xi\circ\chi\right)\in\tilde{E}_{22}^{3}.
\end{multline*}
 Traversing the diagram in either direction gives 
\[
\left(E\xi(ee'),E\xi(ee'),\tau(ee_{3})\otimes e_{1},\tau(ee_{3})\otimes e_{2},\tau E^{2}\circ E^{2}\xi\circ\chi\right).
\]
\end{proof}
The proposition that $\tilde{E}\tilde{\tau}$ and $\tilde{\tau}\tilde{E}$
correspond to $\tilde{\tau}_{1}$ and $\tilde{\tau}_{2}$ is now proved.
\end{proof}

\subsection{Definition of $\mathcal{L}(1)\dotimes\mathcal{V}$}
\begin{defn}
\label{def:Main-product-def} Let $\mathcal{V}$ be a $2$-representation
of $\mathcal{U}^{+}$ given by the data $(E,x,\tau)$ for a $k$-algebra
$A$ such that $_{A}E$ is finitely generated and projective and $E^{n}$
is free as a $P_{n}$-module. We define $\mathcal{L}(1)\ootimes\mathcal{V}$
to be the $2$-representation of $\mathcal{U}^{+}$ given for the
$k$-algebra $C$ by the data $(\tilde{E},\tilde{x},\tilde{\tau})$.
\end{defn}

\begin{prop}
If $E$ is locally nilpotent, then $\tilde{E}$ is locally nilpotent.
\end{prop}

\begin{proof}
Note that in our setting of bimodules, local nilpotence of $E\otimes_{A}-$
is equivalent to nilpotence of $E$, meaning that $E^{n}\cong0$ for
some $n$. This is because local nilpotence implies $E^{n}\otimes_{A}A\cong0$
for some $n$, but that is just $E^{n}$ as a bimodule.

Recall the expression for $\tilde{E}^{n}$ as a matrix of $(A[y],A[y])$-bimodules:
\[
\begin{pmatrix}y_{1}\dots y_{n}E^{n}[y] & y_{1}\dots y_{n+1}E^{n+1}[y]\\
G_{n} & G_{n+1}
\end{pmatrix}\iso\begin{pmatrix}\tilde{E}_{11}^{n} & \tilde{E}_{12}^{n}\\
\tilde{E}_{21}^{n} & \tilde{E}_{22}^{n}
\end{pmatrix}.
\]
 The method we used to compute a model for $G_{n}$ for $n=1,2,3$
also shows that $G_{n}$ for any $n$ can be described as a sub-bimodule
of $E^{n-1}[y]^{\oplus n}\oplus\Hom_{A}(_{A}E,E^{n})[y]$, given by
the elements satisfying a certain set of conditions. It follows that
$G_{n}$ vanishes for large $n$ if $E^{n}$ does. Also $y_{1}\dots y_{n}E^{n}[y]$
vanishes for large $n$ because $E^{n}$ does. It follows that $\tilde{E}$
is nilpotent.
\end{proof}

\subsubsection{Weights and gradings for $\mathcal{L}(1)\dotimes\mathcal{V}$}

It frequently happens that a $2$-representation has additional structure,
and we may ask whether our $2$-product inherits that structure. A
$2$-representation of $\mathcal{U}^{+}$ may have a weight decomposition,
or its algebra may have a grading.
\begin{defn}
A $2$-representation $\mathcal{V}$ of $\mathcal{U}^{+}$ given for
$k$-algebra $A$ by the data $(E,x,\tau)$ is said to have a weight
decomposition when $A$ has the form $A=\prod_{i\in\z}A_{i}$ with
units $e_{i}\in A_{i}$, and $e_{j}Ee_{i}=\delta_{i+2,j}\cdot e_{i+2}Ee_{i}$.
\end{defn}

\begin{prop}[weight decomposition]
 Let $A$ and $(E,x,\tau)$ satisfy the conditions of Def.~\ref{def:Main-product-def},
and let $\mathcal{V}$ be the $2$-representation they determine.
Suppose that $\mathcal{V}$ has a weight decomposition with units
$e_{i}\in A_{i}$. Let $C$ and $(\tilde{E},\tilde{x},\tilde{\tau})$
give the data of $\mathcal{L}(1)\ootimes\mathcal{V}$. Then $C$ has
a weight decomposition $C=\prod_{i\in\z}C_{i}$ with $C_{i}=f_{i}Cf_{i}$
where the units $f_{i}\in C_{i}\subset C$ are given in matrix form
as follows: 
\[
f_{i}=\begin{pmatrix}e_{i+1} & 0\\
0 & (e_{i-1},\_.e_{i-1})
\end{pmatrix}.
\]
\end{prop}

\begin{proof}
The elements $f_{i}$ are clearly idempotent and orthogonal, and they
sum to the identity. We have for the matrix components of $f_{j}\tilde{E}f_{i}$:
\begin{align*}
[f_{j}\tilde{E}f_{i}]_{11} & =e_{j+1}.y_{1}E[y].e_{i+1}\\{}
[f_{j}\tilde{E}f_{i}]_{12} & =e_{j+1}.y_{1}y_{2}E^{2}[y].e_{i-1}\\{}
[f_{j}\tilde{E}f_{i}]_{21} & =G_{1}\bigcap\Bigl(e_{j-1}A[y]e_{i+1}\oplus e_{j-1}.\Hom_{A}(_{A}E,E).e_{i+1}[y]\Bigr)\\{}
[f_{j}\tilde{E}f_{i}]_{22} & =G_{2}\bigcap\Bigl(e_{j-1}.E[y].e_{i-1}\oplus e_{j-1}.E[y].e_{i-1}\oplus e_{j-1}.\Hom_{A}(_{A}E,E^{2}).e_{i-1}[y]\Bigr).
\end{align*}
 These are all zero unless $j=i+2$.
\end{proof}
\begin{defn}[graded case]
 A $2$-representation $\mathcal{V}$ of $\mathcal{U}^{+}$ given
for $k$-algebra $A$ by the data $(E,x,\tau)$ is said to be a $\z$-graded
$2$-representation when $A$ is a $\z$-graded $k$-algebra, $E$
is a graded bimodule, and $x$ and $\tau$ are graded endomorphisms
with $\deg x=+2$ and $\deg\tau=-2$.
\end{defn}

\begin{prop}
Let $A$ and $(E,x,\tau)$ satisfy the conditions of Def.~\ref{def:Main-product-def},
and let $\mathcal{V}$ be the $2$-representation they determine.
Suppose that $\mathcal{V}$ is a $\z$-graded $2$-representation.
Let $C$ and $(\tilde{E},\tilde{x},\tilde{\tau})$ give the data of
$\mathcal{L}(1)\ootimes\mathcal{V}$. Then $\mathcal{L}(1)\ootimes\mathcal{V}$
is a $\z$-graded $2$-representation. The gradings on generators
in $C$ and $\tilde{E}$ are inherited from the gradings on $A$ and
$E$ with the assumption that $\deg y=+2$.
\end{prop}

\begin{proof}
It is trivial to check that $C$ is graded and $\tilde{E}$ is a graded
bimodule. The formulas for $\tilde{x}$ and $\tilde{\tau}$ in Def.~\ref{def:tilde-x-tilde-tau}
show that they have the right degrees.
\end{proof}

\section{\label{sec:Comparison} Comparison: $\mathcal{V}=\mathcal{L}(1)$}

In §\ref{subsec:KLRW-L(1)xL(1)} we describe a well-known $2$-representation
of $\mathcal{U}^{+}$ categorifying $L(1)\otimes L(1)$ using Soergel
bimodules. In §\ref{subsec:L(1)ootimesL(1)} we describe our product
explicitly for $\mathcal{V}=\mathcal{L}(1)$, and in §\ref{subsec:Comparison-L(1)xL(1)}
we show that the result is equivalent to the known one. The reader
is warned that notations in this section will diverge from the previous
sections.

Let $P_{2}=k[y_{1},y_{2}]$. Let $S_{2}$ denote the symmetric group
on $2$ letters, generated by $t_{1}$, and acting on $P_{2}$ by
permutation of the $y_{i}$. Let $P_{2}^{S_{2}}$ be the subalgebra
generated by invariant homogeneous polynomials.

\subsection{\label{subsec:KLRW-L(1)xL(1)}A categorification of $L(1)\otimes L(1)$}
\begin{defn}
We define:
\begin{itemize}
\item a $(P_{2},P_{2})$-bimodule $B_{s_{1}}=P_{2}\otimes_{P_{2}^{S_{2}}}P_{2}$ 
\begin{itemize}
\item and observe that $B_{s_{1}}$ is also a $P_{2}$-algebra with structure
map $P_{2}\to B_{s_{1}}$ given by $f\mapsto1\otimes f$
\item and that $P_{2}$ is a left $B_{s_{1}}$-module by $(f\otimes g).\theta=fg\theta$
\end{itemize}
\item a $P_{2}$-algebra $T=T_{+2}\oplus T_{0}\oplus T_{-2}$ by 
\[
T_{+2}=P_{2},\;T_{0}=\End_{B_{s_{1}}}(P_{2}\oplus B_{s_{1}})^{\op},\;T_{-2}=P_{2}
\]
\item a $(T,T)$-bimodule $\mathscr{E}=_{+2}\mathscr{E}_{0}\oplus\phantom{}_{0}\mathscr{E}_{-2}$
by 
\begin{align*}
_{0}\mathscr{E}_{-2} & =\begin{pmatrix}P_{2}\\
B_{s_{1}}
\end{pmatrix}\cong T_{0}e_{2}\\
_{+2}\mathscr{E}_{0} & =\begin{pmatrix}P_{2} & B_{s_{1}}\end{pmatrix}\cong e_{2}T_{0}
\end{align*}
 for $e_{2}$ the projection onto $B_{s_{1}}$ 
\begin{itemize}
\item and observe the canonical isomorphism 
\[
_{+2}\mathscr{E}_{-2}^{2}=e_{2}T_{0}\otimes_{T_{0}}T_{0}e_{2}\iso B_{s_{1}}
\]
\end{itemize}
\item a bimodule endomorphism $x\in\End(\mathscr{E})$ by 
\[
_{+2}x_{0}=\begin{pmatrix}y_{2} & y_{2}\otimes1\end{pmatrix},\quad{}_{0}x_{-2}=\begin{pmatrix}y_{1}\\
y_{1}\otimes1
\end{pmatrix}
\]
(acting by multiplication)
\item a bimodule endomorphism $\tau\in\End(\mathscr{E}^{2})$ by 
\[
_{+2}\tau_{-2}:f\otimes g\mapsto\partial_{t_{1}}(f)\otimes g
\]
 where $\partial_{t_{1}}\in\End_{k}(P_{2})$ is a Demazure operator:
\[
\partial_{s_{1}}:f\mapsto\frac{f-f^{t_{1}}}{y_{1}-y_{2}}.
\]
\end{itemize}
\end{defn}

The next theorem is well-known. Cf., for example, Lauda \cite{laudaCategorifiedQuantumMathfraksl2009},
Webster \cite[§2.3]{websterTensorProductAlgebras2016}, Stroppel \cite[§5.1.1]{stroppelCategoryQuiversEndomorphism2003},
Sartori-Stroppel \cite{sartoriCategorificationTensorProduct2015}:
\begin{thm}
The $k$-algebra $T$ and triple $(\mathscr{E},x,\tau)$ defined above
gives a $2$-representation of $\mathcal{U}^{+}$, called $\mathscr{T}$
below, that categorifies the tensor product $L(1)\otimes L(1)$ of
fundamental representations of $\mathfrak{sl}_{2}$.
\end{thm}

\subsection{\label{subsec:L(1)ootimesL(1)}$\mathcal{L}(1)\dotimes\mathcal{L}(1)$}

We notate both factors as in §\ref{subsec:Minimal-rep} except that
on the right factor we use $y_{1}$ in place of $y$, and on the left
factor we use $y_{2}$ in place of $y$. We write $E_{i}$, $x_{i}$,
$\tau_{i}$, $i=1,2$ for the $2$-representation data on the right
($i=1$) and on the left $(i=2$).

In the formulas we have given for the product, the algebra $A$, now
$A_{1}$, becomes $k[y_{1}]_{+1}\times k[y_{1}]_{-1}$ (in its weight
decomposition), $E$ becomes $k[y_{1}]$, $x$ becomes $y_{1}$, and
$y$ becomes $y_{2}$. Let $\omega=y_{1}-y_{2}\in P_{2}$. So $\omega$
will take over the role of \textquoteleft $y_{1}=x-y$\textquoteright{}
that was written in previous sections. Write $\pi:P_{2}\to P_{2}/(\omega)$
for the projection.

We let $B$, $X$, $E'$, $C$, $\tilde{E}$, $\tilde{x}$, and $\tilde{\tau}$
be defined as above. The algebra $B$ and complex $X$ have nonzero
elements only in weights $-2$, $0$, $+2$. These are given as follows:
\begin{align*}
B_{-2} & =\begin{pmatrix}P_{2} & 0\\
0 & 0
\end{pmatrix}, & X_{1_{-2}} & =\begin{pmatrix}P_{2}\\
0
\end{pmatrix}, & X_{2_{-2}} & =\begin{pmatrix}0\\
0
\end{pmatrix},\\
B_{0} & =\begin{pmatrix}P_{2} & k[y]\\
0 & P_{2}
\end{pmatrix}, & X_{1_{0}} & =\begin{pmatrix}P_{2}\\
0
\end{pmatrix}, & X_{2_{0}} & =\begin{pmatrix}P_{2}\xrightarrow{\pi}P_{2}/(\omega)\\
0\to P_{2}
\end{pmatrix},\\
B_{+2} & =\begin{pmatrix}0 & 0\\
0 & P_{2}
\end{pmatrix}, & X_{1_{+2}} & =\begin{pmatrix}0\\
0
\end{pmatrix}, & X_{2_{+2}} & =\begin{pmatrix}0\\
0\to P_{2}
\end{pmatrix}.
\end{align*}
 Here the action of $P_{2}/(\omega)$ from the upper right of $B_{0}$
on $X_{2_{0}}$ is $P_{2}/(\omega)\otimes_{P_{2}}P_{2}\to P_{2}/(\omega)$
given by $f\otimes1\mapsto f$. The complexes for $X$ start in degree
$0$ on the left. The matrix coefficients are in each case from the
$-1$ weight space of $A_{2}$ in the upper left corner.

To compute $\tilde{E}$ we will also need $E'X_{2}$, which is: 
\begin{align*}
_{0}E'_{-2}(X_{2_{-2}}) & =\begin{pmatrix}0\\
0
\end{pmatrix},\\
_{+2}E'_{0}(X_{2_{0}}) & =\begin{pmatrix}0\\
0\to P_{2}\oplus P_{2}\xrightarrow{(-\pi,\pi)}P_{2}/(\omega)
\end{pmatrix}.
\end{align*}

Next we compute $C$: 
\[
[C_{+2}]=\begin{pmatrix}0 & 0\\
0 & P_{2}
\end{pmatrix},\;[C_{0}]=\begin{pmatrix}P_{2} & \omega P_{2}\\
P_{2} & Q_{1}^{\op}
\end{pmatrix},\;[C_{-2}]=\begin{pmatrix}P_{2} & 0\\
0 & 0
\end{pmatrix}.
\]
 Here $Q_{1}^{\op}\subset P_{2}\oplus P_{2}$ is the (commutative)
algebra of all $(\theta,\varphi)$ such that $\varphi-\theta\in\omega P_{2}$,
with componentwise multiplication. It is a $P_{2}$-algebra by $P_{2}\ni f\mapsto(f,f)\in Q_{1}$.
The algebra structure of $C_{0}$ (cf.~§\ref{subsec:matrix-tensor-products})
may be described as follows. The upper right term, $\omega P_{2}$,
is a left $P_{2}$-module by multiplication. It is a right $Q_{1}^{\op}$-module
with $(\theta,\varphi)$ acting by multiplication by $\varphi$. The
lower left $P_{2}$ is a left $Q_{1}^{\op}$-module with the same
action. It has a right $P_{2}$ action by multiplication. The remaining
structure maps are: 
\begin{align}
\omega P_{2}\otimes_{P_{2}}P_{2} & \to P_{2}\label{eq:C-struct-1}\\
\text{by }\omega\theta'\otimes\theta & \mapsto\omega\theta\theta'\nonumber 
\end{align}

and 
\begin{align}
P_{2}\otimes_{P_{2}}\omega P_{2} & \to Q_{1}^{\op}\label{eq:C-struct-2}\\
\text{by }\theta\otimes\omega\theta' & \mapsto(0,\omega\theta\theta').\nonumber 
\end{align}

Now compute $\tilde{E}$ and the endomorphisms $\tilde{x}$ by components:
\begin{align*}
_{0}[\tilde{E}]_{-2} & =\begin{pmatrix}\omega P_{2} & 0\\
Q_{1} & 0
\end{pmatrix}, & _{0}[\tilde{x}]_{-2} & =\begin{pmatrix}y_{1} & 0\\
(y_{2},y_{1}) & 0
\end{pmatrix},\\
_{+2}[\tilde{E}]_{0} & =\begin{pmatrix}0 & 0\\
P_{2} & Q_{2}
\end{pmatrix}, & _{+2}[\tilde{x}]_{0} & =\begin{pmatrix}0 & 0\\
y_{2} & (y_{2},y_{1})
\end{pmatrix},
\end{align*}
 where $Q_{2}\subset P_{2}\oplus P_{2}$ is the $(P_{2},Q_{1}^{\op})$-bimodule
containing all $(e_{1},e_{2})$ such that $e_{1}-e_{2}\in\omega P_{2}$;
$Q_{1}^{\op}$ acts on $Q_{2}$ on the right by $(e_{1},e_{2}).(\theta,\varphi)=(e_{1}\varphi,e_{2}\theta)$
(note the swap), and $P_{2}$ on the left by diagonal multiplication.

In the next section it will be useful to view $_{0}\tilde{E}_{-2}$
as $C_{0}q_{2}$ using the idempotent $q_{2}=\left(\begin{smallmatrix}0 & 0\\
0 & 1
\end{smallmatrix}\right)\in[C_{0}]$, and to view $_{+2}\tilde{E}_{0}$ as $q_{2}C_{0}$ using the isomorphism
of $(P_{2},Q_{1}^{\op})$-bimodules $\sigma:Q_{1}\iso Q_{2}$ by $(\theta,\varphi)\mapsto(\varphi,\theta)$.
Viewing them in this way, we may write $_{0}\tilde{x}_{-2}$ as multiplication
on $C_{0}q_{2}$ on the left by $\left(\begin{smallmatrix}y_{1} & 0\\
0 & (y_{2},y_{1})
\end{smallmatrix}\right)\in C_{0}$, and $_{+2}\tilde{x}_{0}$ as multiplication on $q_{2}C_{0}$ on
the right by $\left(\begin{smallmatrix}y_{2} & 0\\
0 & (y_{1},y_{2})
\end{smallmatrix}\right)\in C_{0}$ (note the swap).

Finally, compute $\tilde{E}^{2}$ and $\tilde{\tau}$ by components:
\begin{align*}
_{+2}[\tilde{E}^{2}]_{-2} & =\begin{pmatrix}0 & 0\\
Q_{2} & 0
\end{pmatrix}, & _{+2}[\tilde{\tau}]_{-2} & =\begin{pmatrix}0 & 0\\
t_{21} & 0
\end{pmatrix},
\end{align*}
where 
\[
t_{21}:(e_{1},e_{2})\mapsto(\omega^{-1}(e_{1}-e_{2}),\omega^{-1}(e_{1}-e_{2})).
\]

\subsection{\label{subsec:Comparison-L(1)xL(1)}Comparison}
\begin{thm}
\label{thm:comparison} There is an equivalence $\mathcal{L}(1)\ootimes\mathcal{L}(1)\iso\mathscr{T}$
of $2$-representations.
\end{thm}

We will use a few intermediate steps.

Define a new algebra $R$: 
\[
R=P_{2}[e]\big/\bigl(e^{2}-\omega e\bigr).
\]
 There is a map of $P_{2}$-algebras $R\xrightarrow{\gamma}B_{s_{1}}$
given by $e\mapsto1\otimes y_{1}-y_{1}\otimes1$. There is another
map of $P_{2}$-algebras $R\xrightarrow{\gamma'}Q_{1}^{\op}$ given
by $P_{2}\ni f\mapsto(f,f)\in Q_{1}^{\op}$ and $e\mapsto(\omega,0)$.
\begin{lem}
The maps $\gamma$ and $\gamma'$ are isomorphisms of $P_{2}$-algebras.
\end{lem}

\begin{proof}
Straightforward.
\end{proof}
We will also use the composition $\sigma\circ\gamma'$ to obtain an
isomorphism of $(P_{2},P_{2})$-bimodules $R\iso Q_{2}$ given by
$f\mapsto(f,f)$ and $e\mapsto(0,\omega)$.

Now we translate $\mathscr{T}$ using $\gamma$. The action of $B_{s_{1}}$
on $P_{2}$ induces an action of $R$ on $P_{2}$ through $\gamma$,
according to which $P_{2}\hookrightarrow R$ acts on $P_{2}$ by multiplication,
and $e$ acts by zero. We have an isomorphism of $R$-modules $P_{2}\iso R/(e)$
using this action on $P_{2}$. In the remainder of this section we
assume this isomorphism and write $R$ in place of $B_{s_{1}}$ everywhere
in the $2$-representation $\mathscr{T}$. Under this translation,
and using the decomposition $R\iso P_{2}\oplus P_{2}e$ as $P_{2}$-modules,
we have: 
\[
_{+2}x_{0}=\begin{pmatrix}y_{2} & y_{2}+e\end{pmatrix},\quad{}_{0}x_{-2}=\begin{pmatrix}y_{1}\\
y_{1}-e
\end{pmatrix},
\]
 and 
\[
_{+2}\tau_{-2}=\bigl(p_{1}+p_{2}e\mapsto-p_{2}\bigr).
\]

\begin{lem}
The matrix presentation of $T_{0}$ is given by: 
\[
\begin{pmatrix}P_{2} & P_{2}\\
P_{2} & R
\end{pmatrix}\iso T_{0},
\]
 where: 
\begin{itemize}
\item for $[T_{0}]_{11}$ the map sends $\theta\in P_{2}$ to $(1\mapsto\theta)\in\End_{R}(P_{2})^{\op}$
\item for $[T_{0}]_{21}$ the map sends $\theta\in P_{2}$ to $(1\mapsto\theta)\in\Hom_{R}(R,P_{2})$
\item for $[T_{0}]_{12}$ the map sends $\theta\in P_{2}$ to $(1\mapsto\theta\omega-\theta e)\in\Hom_{R}(P_{2},R)$
\item for $[T_{0}]_{22}$ the map sends $r\in R$ to $(1\mapsto r)\in\End_{R}(R,R)^{\op}.$
\end{itemize}
The algebra structure maps (cf.~§\ref{subsec:matrix-tensor-products})
are given as follows: 
\begin{itemize}
\item $[T_{0}]_{11}\acts[T_{0}]_{12}$ by $\theta.\theta'=\theta\theta'$
\item $[T_{0}]_{21}\acted[T_{0}]_{11}$ by $\theta'.\theta=\theta'\theta$
\item $[T_{0}]_{12}\acted[T_{0}]_{22}$ by $\theta.(p_{1}+p_{2}e)=\theta p_{1}$
\item $[T_{0}]_{22}\acts[T_{0}]_{21}$ by $(p_{1}+p_{2}e).\theta=p_{1}\theta$
\item $[T_{0}]_{12}\otimes[T_{0}]_{21}\to[T_{0}]_{11}$ by $\theta\otimes\theta'\mapsto\omega\theta\theta'$
\item $[T_{0}]_{21}\otimes[T_{0}]_{12}\to[T_{0}]_{22}$ by $\theta'\otimes\theta\mapsto\omega\theta'\theta-\theta'\theta e$.
\end{itemize}
\end{lem}

\begin{proof}
Let us explain the map to $[T_{0}]_{12}$. Recall that $P_{2}\cong R/(e)$.
An element of $\Hom_{R}(R/(e),R)$ is given by the image $r=p_{1}+p_{2}e$
of $1$, which may be anything satisfying $e.r=0$, and that condition
is equivalent to $p_{1}=-p_{2}\omega$. The other morphisms and the
structure maps are easily computed.
\end{proof}
\begin{lem}
Let $\Phi_{0}:T_{0}\to C_{0}$ be given on components by: 
\[
\begin{pmatrix}\idop_{P_{2}} & \omega\\
\idop_{P_{2}} & \gamma'
\end{pmatrix}.
\]
 Then $\Phi_{0}$ is an isomorphism of $P_{2}$-algebras.
\end{lem}

\begin{proof}
The specified maps give algebra isomorphisms on the diagonal components,
and $k$-module isomorphisms on the off-diagonal components. Now we
check equivariance under the bimodule structure maps. The only nonobvious
cases concern maps involving the lower right component.

An element of $Q_{1}^{\op}$ may be written uniquely as a sum $(\omega\theta,0)+(\varphi,\varphi)$.
This is sent by $\gamma'^{-1}$ to $\varphi+\theta e\in R$. So the
action of $(\theta,\varphi)$ by multiplication by $\varphi$ agrees
with the action of $p_{1}+p_{2}e$ by multiplication by $p_{1}$.
The structure map $[T_{0}]_{12}\otimes[T_{0}]_{21}\to[T_{0}]_{11}$
clearly agrees with Eq.~\ref{eq:C-struct-1} through $\Phi_{0}$.
The map $[T_{0}]_{21}\otimes[T_{0}]_{12}\to[T_{0}]_{22}$ agrees with
Eq.~\ref{eq:C-struct-2} through $\Phi_{0}$ because $\gamma':\omega\theta'\theta-\theta'\theta e\mapsto(0,\omega\theta\theta')$.
\end{proof}
\begin{proof}[Proof of Theorem \ref{thm:comparison}]
 Extend $\Phi_{0}$ to an algebra isomorphism $\Phi:T\iso C$ by
$\Phi_{+2}=\idop_{P_{2}}$ and $\Phi_{-2}=\idop_{P_{2}}$. It remains
to check compatibility with the actions of $E$, $x$, and $\tau$
in $\mathcal{U}^{+}$, and this poses no difficulty. We summarize
that now.

We have in $\mathscr{T}$ that $_{0}\mathscr{E}_{-2}\iso T_{0}r_{2}$
for $r_{2}=\left(\begin{smallmatrix}0 & 0\\
0 & 1
\end{smallmatrix}\right)\in[T_{0}]$, and similarly $_{0}\tilde{E}{}_{-2}=C_{0}q_{2}$ in $\mathcal{L}(1)\ootimes\mathcal{L}(1)$;
and we have $q_{2}=\Phi_{0}(r_{2})$. The action of $_{0}x_{-2}$
on $_{0}\mathscr{E}_{-2}$ in $\mathscr{T}$ can be written in $T_{0}r_{2}$
as multiplication on the left by $\left(\begin{smallmatrix}y_{1} & 0\\
0 & y_{1}-e
\end{smallmatrix}\right)\in[T_{0}]$. In $\mathcal{L}(1)\ootimes\mathcal{L}(1)$ it is written as multiplication
on the left by $\left(\begin{smallmatrix}y_{1} & 0\\
0 & (y_{2},y_{1})
\end{smallmatrix}\right)$. These correspond using $\gamma':R\iso Q_{1}^{\op}$. Similarly for
$_{+2}x_{0}$ since $\gamma':R\ni y_{2}+e\mapsto(y_{1},y_{2})\in Q_{1}^{\op}$.
Finally, the action of $\tau$ in $R$ by $_{+2}\tau_{-2}=\bigl(p_{1}+p_{2}e\mapsto-p_{2}\bigr)$
corresponds to $_{+2}\tilde{\tau}_{-2}$, now using $\sigma\circ\gamma':R\iso Q_{2}$.
\end{proof}
\bibliographystyle{amsalpha}
\providecommand{\bysame}{\leavevmode\hbox to3em{\hrulefill}\thinspace}
\providecommand{\MR}{\relax\ifhmode\unskip\space\fi MR }
\providecommand{\MRhref}[2]{%
	\href{http://www.ams.org/mathscinet-getitem?mr=#1}{#2}
}
\providecommand{\href}[2]{#2}

\end{document}